\documentclass[11pt,a4paper]{amsart}

\usepackage{fullpage}
\usepackage{geometry}

\usepackage{amsmath,amssymb,amsthm}

\usepackage{tikz}
\usepackage{pgfplots}
\usepackage{pgfplotstable}
\usetikzlibrary{matrix,fit} 
\pgfplotsset{compat=1.18}

\usepackage{xcolor}

\usepackage[colorlinks=true, linkcolor=blue, citecolor=blue]{hyperref}
\usepackage{cleveref}

\usepackage[T1]{fontenc}
\usepackage[utf8]{inputenc}

\theoremstyle{plain}
\newtheorem{theorem}{Theorem}[section]
\newtheorem{lemma}[theorem]{Lemma}
\newtheorem{proposition}[theorem]{Proposition}
\newtheorem{corollary}[theorem]{Corollary}

\theoremstyle{definition}
\newtheorem{definition}[theorem]{Definition}
\newtheorem{hyp}[theorem]{Hypothesis}

\theoremstyle{remark}
\newtheorem{remark}[theorem]{Remark}

\newtheorem*{theorem*}{Theorem}

\newcommand{\E}{\mathbb{E}}
\renewcommand{\P}{\mathbb{P}}
\newcommand{\R}{\mathbb{R}}

\newcommand{\gmm}{\gamma_{k,l,\theta,\rho}}
\newcommand{\cc}{c_{k,l,\theta}}

\newif\ifdraft
\drafttrue 

\ifdraft
  \newcommand{\laure}[1]{{\color{magenta}\textbf{#1}}}
  \newcommand{\lili}[1]{{\color{red}\textbf{#1}}}
  
\else
  \newcommand{\laure}[1]{}
  \newcommand{\lili}[1]{}
\fi

\title{Quantitative Convergence of the Nodal Length of Random Trigonometric Polynomials}

\author{Laure Coutin}
\address{Institut de Mathématiques de Toulouse, Université Toulouse, 118 route de Narbonne, F-31062 Toulouse Cedex 9, France}

\email{laure.coutin@math.univ-toulouse.fr}

\author{Liliana Peralta}
\address{Departamento de Matemáticas, Facultad de Ciencias, Universidad Nacional Autónoma de México, Circuito Exterior, C.U., C.P. 04510, Ciudad de México, Mexico}
\email{lylyaanaa@ciencias.unam.mx}
\thanks{The second author was partially supported by UNAM-DGAPA-PAPIIT grant IA103826.}

\begin{document}

\maketitle

\begin{abstract}
We establish quantitative convergence results for the distribution of the nodal length of two-dimensional random trigonometric polynomials with independent and identically distributed centered coefficients towards that of the corresponding Gaussian limiting field. More precisely, we obtain an explicit polynomial rate of convergence.

The proof combines probabilistic and geometric arguments. First, we establish quantitative convergence estimates for the underlying random fields. Using a suitable functional representation, this problem is reduced to the quantitative convergence of an empirical field towards a Brownian sheet in a Wasserstein-type distance. Second, Crofton's formula allows us to express the nodal length in terms of the zeros of the restrictions of the random trigonometric polynomials to parametrized affine lines. This reduction makes it possible to adapt one-dimensional arguments based on Kac's counting formula and transfer the convergence estimates for the random fields to the associated nodal lengths.
\end{abstract}



\bigskip

\tableofcontents

\section{Introduction}

Random trigonometric polynomials in higher dimensions constitute a fundamental class of random fields. Their interest stems from the rich geometry of their nodal sets and their connections with a variety of problems in mathematics and physics. We refer the reader to \cite{Bogomo} for an overview of the role of random trigonometric polynomials in the study of quantum chaos.

In this work we consider the two-dimensional setting, where the geometry of the nodal set can be described through its nodal length. The study of geometric characteristics of random level sets dates back to the pioneering works of Longuet-Higgins \cite{Longuet-Higgins} and Swerling \cite{Swerling}, who investigated contour lengths and related geometric features of random surfaces. Modern treatments rely on Kac--Rice formulas and the geometry of random fields, see Adler and Taylor \cite{Adler} and Azaïs and Wschebor \cite{Wse-Aza}. Understanding the asymptotic behavior of nodal lengths and obtaining quantitative convergence estimates towards limiting models are natural problems in the theory of random fields. As far as we know, the convergence in distribution and the qualitative convergence of the expectation of the nodal length of a particular family of random bivariate trigonometric polynomials were established in \cite{AngstPhamPoly}. This result motivates the present work. We seek to quantify the rate of the distributional convergence. Another motivation comes from the one-dimensional setting, where quantitative convergence results for the number of zeros of random trigonometric polynomials were obtained in \cite{lclp}.

The random trigonometric polynomials considered in \cite{AngstPhamPoly} differ from the model studied here. Nevertheless, both works rely on the same general philosophy: convergence of random fields to a Gaussian limit implies convergence of the corresponding nodal lengths. More specifically, in our setting, we establish quantitative convergence estimates for the underlying random polynomials towards the limiting Gaussian field and use them to derive quantitative convergence results for the associated nodal lengths. Our work can be easily adapted to the framework of \cite{AngstPhamPoly}.

More specifically, we establish a polynomial rate of convergence $m^{-\alpha}$ for the distribution of the nodal length towards its Gaussian counterpart in Fortet--Mourier distance (see Appendix C of  \cite{NP12} for  a definition). The exponent $\alpha$ can be chosen explicitly in terms of the convergence rate of the underlying random fields.

A major difficulty in dimension two comes from the fact that the nodal length is a geometric functional of the corresponding random field. Unlike the one-dimensional case, where zeros are isolated and can be counted directly, the nodal set in dimension two is a collection of curves whose total length depends on the global geometry of the field.

Similarly to the one-dimensional setting considered in \cite{lclp}, the strategy of the proof consists in transferring quantitative convergence estimates between the random trigonometric polynomials and the limiting Gaussian field to the associated geometric functional. The main tool is Crofton's formula, which expresses the nodal length as an average number of intersections between the nodal set and a family of affine lines. This representation reduces the problem to the study of the zeros of the restrictions of the random trigonometric polynomials to parametrized affine lines. Although these restricted processes are not random trigonometric polynomials, Kac's counting formula \cite{kac} can still be applied to study their zeros. This makes it possible to adapt the main ingredients of the one-dimensional analysis developed in \cite{lclp} to the present two-dimensional setting. Combining these arguments with quantitative convergence estimates for the underlying random fields, we obtain a polynomial rate of convergence for the distribution of the nodal length; see Theorem~\ref{thm:main-rate} for a precise statement.

An intermediate step of the proof consists in establishing quantitative convergence estimates for the underlying random fields. To this end, we first introduce a functional representation of the random trigonometric polynomials and their Gaussian limiting field as images of the same linear operator applied, respectively, to an empirical field and a Brownian sheet. This reduces the problem to the study of the convergence of the empirical field towards the Brownian sheet in a suitable Wasserstein-type distance whose test functions are linearly locally Lipschitz. This result seems to be new, to the best of our knowledge.

The proof follows the same general strategy as in the one-dimensional setting. First, we establish uniform regularity estimates for the empirical field and the Brownian sheet in suitable Hölder and fractional Sobolev spaces. We then introduce finite-dimensional interpolation operators and control the associated approximation error. This reduction allows us to apply Stein's method in a finite-dimensional setting and obtain quantitative convergence estimates (see Theorem~\ref{thm-Stein-filed-lip11}). These results suggest that in the $d$-dimensional setting, the rate of convergence depends on the dimension.

The paper is organized as follows. In Section~\ref{RTP_GL}, we introduce the random trigonometric polynomials and their Gaussian limiting field, and establish a functional representation linking both objects to an empirical field and a Brownian sheet. Section~\ref{rat_dons} is devoted to a quantitative Donsker-type theorem for random fields, which provides explicit convergence rates towards the Gaussian limit. Section~\ref{crofton:error} contains the Crofton representation and the decomposition of the approximation error, while Section~\ref{main:result} is devoted to the optimization of the parameters and the proof of the final convergence rate. Finally, Appendix~\ref{moment_bounds} contains auxiliary estimates for the number of zeros of the restricted processes, whereas Appendix~\ref{appendix-Stein} recalls the Stein-type arguments used in the proof of the quantitative convergence theorem.

\section{Random trigonometric polynomials and their Gaussian limit}\label{RTP_GL}

In this section, we introduce the random trigonometric polynomials and state the assumptions on their coefficients. We also introduce the limiting Gaussian field and its covariance structure.

\smallskip

Let $\{a_{k,\ell},\, k,\ell \in \mathbb{N}\}$ and $\{b_{k,\ell},\, k,\ell \in \mathbb{N}\}$ be two independent sequences of independent and identically distributed random variables, defined on a common probability space $(\Omega, \mathcal{F}, \P)$, such that $\E[a_{k,\ell}] = 0$ and $\E[a_{k,\ell}^2] = 1$.

We consider the random trigonometric polynomials defined by
\begin{align}\notag
P_m & : \mathbb{R}^2 \longrightarrow \mathbb{R}, \\
\label{def:pol_1}
P_m(t) &= \frac{1}{m} \sum_{k=1}^{m} \sum_{\ell=1}^{m} \operatorname{Re} \left[ \left(a_{k,\ell} + i b_{k,\ell} \right) e^{-i \left( \frac{k t_1 + \ell t_2}{m} \right)} \right], \quad t=(t_1,t_2) \in \mathbb{R}^2.
\end{align}

It is known that, under suitable assumptions, the sequence $(P_m)_{m\in\mathbb{N}}$ converges in distribution to a centered Gaussian process $G$ on $\mathbb{R}^2$ (see \cite{AngstPhamPoly} for the non-stationary case. Note that, although the result therein is stated without the sine part, the independence of $(a_{k,\ell})$ and $(b_{k,\ell})$ ensures that it extends to our setting). Its covariance function is given by
\begin{align*}
R_G(t,s)= \frac{\cos (t_1+s_1+t_2+s_2) -\cos(t_1+s_1) -\cos (t_2+s_2) + 1}{(t_1+s_1)(t_2+s_2)}, \quad t,s \in \mathbb{R}^2. 
\end{align*}

Our main result is

\begin{theorem}\label{thm:main-rate}
Let $(P_m)_{m\geq 1}$ be the sequence of random trigonometric polynomials defined in \eqref{def:pol_1}. Assume that the family $
\left\{(a_{k,\ell},b_{k,\ell}) : k,\ell\in\mathbb N\right\}$
consists of independent and identically distributed centered random vectors in $\mathbb R^2$ with identity covariance matrix and
\begin{equation*}
\E\left[\|(a_{1,1},b_{1,1})\|^p\right]<\infty
\end{equation*}
for some $p>4$. Let $a>0$ and let
\begin{equation*}
0<\gamma<\frac{1-\frac{4}{p}}{5-\frac{4}{p}}.
\end{equation*}
Let $G$ be the Gaussian field defined in Section~\ref{RTP_GL}. Then, for every
\begin{equation*}
0<\alpha<\frac{\gamma}{4},
\end{equation*}
there exists a constant $C=C(a,p,\gamma,\alpha)>0$ such that
\begin{equation*}
\sup_{\|F\|_{\mathrm{Lip}}+\|F\|_{\infty}\leq 1}
\left|
\E\!\left[F\!\left(\mathcal L(P_m,B(0,a))\right)\right]
-
\E\!\left[F\!\left(\mathcal L(G,B(0,a))\right)\right]
\right|
\leq
C\,m^{-\alpha}.
\end{equation*}
\end{theorem}

\begin{remark}
The restriction to the Fortet--Mourier distance in Theorem~\ref{thm:main-rate} is only due to the event on which the nodal set of $P_m$ is degenerate, that is, on which $P_m$ has a critical point on its nodal set. On this event, Crofton's formula may fail for $P_m$ (see Section~\ref{crofton:error} and Appendix~\ref{app-nondeg}). This restriction can be removed under additional assumptions on the law of the coefficients. Assume that the random variables $a_{k,\ell}$ and $b_{k,\ell}$ have a common density $\rho=e^{-\Phi}$, where $\Phi\in C^\infty(\mathbb R)$ satisfies
\begin{equation*}
\Phi^{(j)}\in\bigcap_{q\geq1}L^q\big(e^{-\Phi}\,dx\big),\qquad j\geq1,
\end{equation*}
where $\Phi^{(j)}$ denotes the $j$-th derivative of $\Phi$. For $t\in\mathbb R^2$, let $Z_m(t)=(P_m(t),\partial_{t_1}P_m(t),\partial_{t_2}P_m(t))$. The Gram matrix of $Z_m(t)$ with respect to the coefficients is non-degenerate: it does not depend on $t$, and its determinant equals
\begin{equation*}
\frac{(m+1)^2(m-1)(7m+5)}{144m^4}>0, \qquad m\geq 2.
\end{equation*}
Using integration by parts with respect to the coefficients, as in \cite[Corollary~10]{ADLNP18}, one can show that the law of $Z_m(t)$ admits a smooth density, bounded uniformly in $m\geq2$ and $t\in\mathbb R^2$. As a consequence, the nodal set of $P_m$ is almost surely non-degenerate. Consequently, Lemma~\ref{crofton-lemma} applies almost surely to the paths of $P_m$, the proof of Theorem~\ref{thm:main-rate} holds for all $1$-Lipschitz functions $F$, and the same rate of convergence $m^{-\alpha}$ is obtained in the Wasserstein-1 distance.
\end{remark}

Both $P_m$ and its Gaussian limit $G$ can be expressed as functionals of an empirical field and of the Brownian sheet, respectively. This will be made precise in the next subsections.

\subsection{Link with the Donsker theorem for fields}\label{link:donsk}
Define
$$z_{k,\ell}= a_{k,\ell} + i b_{k,\ell}, \quad k,\ell \in \mathbb{N}.$$

For $m \in \mathbb{N}$, we define
\begin{align*}
\tilde{S}_m(u,v)
&=\frac{1}{m}\sum_{k=1}^{m}\sum_{l=1}^{m}
z_{k,l}\,
\mathbf{1}_{[\frac{k}{m},1]}(u)\,
\mathbf{1}_{[\frac{l}{m},1]}(v),
\quad (u,v) \in [0,1]^2,
\\
B(u,v) &= B^1(u,v) + i B^2(u,v),
\end{align*}
where $(B^1,B^2)$ are two independent Wiener fields on $[0,1]^2$,  whose precise definition will be given in Section~\ref{FandR}. We now introduce a representation that allows us to express the random polynomials defined in \eqref{def:pol_1} as functionals of the random field $\tilde{S}_m$. This can be seen as a two-dimensional analogue of \cite[Lemma 3.1]{lclp}, where a similar representation is established in the one-dimensional case.

\begin{lemma}\label{lem-def-theta}
Let $a>0$. Let $\Theta:L^{\infty}([0,1]^2,\R^2)\rightarrow C^{\infty}([-a,a]^2, \R)$ be the linear map defined by
\begin{align*}
\Theta(f)(t_1,t_2) = 
&\ \operatorname{Re} \left\{
e^{-i t_1} \left[
f(1,1)e^{-i t_2} - f(1,0) + i t_2 \int_0^1 f(1,v) e^{-i v t_2} \, dv
\right] \right\} \\
& - \operatorname{Re} \left\{
f(0,1)e^{-i t_2} - f(0,0) + i t_2 \int_0^1 f(0,v) e^{-i v t_2} \, dv
\right\} \\
& + \operatorname{Re} \left\{
i t_1 \int_0^1 e^{-i u t_1} \left[
f(u,1)e^{-i t_2} - f(u,0) + i t_2 \int_0^1 f(u,v) e^{-i v t_2} \, dv
\right] du \right\}.
\end{align*}
\begin{itemize}
\item[(i)]
There exist a positive constant $K_0$ and, for all $n\geq 1$ and all multi-indices $j\in \{1,2\}^{n}$, positive constants $K_j$ such that
\begin{align*}
\|\Theta(f)\|_{\infty}&\leq K_0 \|f\|_{\infty},\\
\|\partial_{t_{j_1},...t_{j_n}}\Theta(f)\|_{\infty}&\leq K_j \|f\|_{\infty}.
\end{align*}

\item[(ii)]
If $f_{x_0,y_0}(t_1,t_2)={\mathbf 1}_{[x_0,1]}(t_1)\, {\mathbf 1}_{[y_0,1]}(t_2)$, where $x_0,y_0 \in (0,1]$, then 
\begin{align*}
\Theta(f_{x_0,y_0})(t_1,t_2)
= \operatorname{Re}\big(e^{-i(t_1 x_0 + t_2 y_0)}\big).
\end{align*}
\end{itemize}
\end{lemma}

\begin{proof}
We first prove (ii). Let $x_0,y_0 \in (0,1]$ and set $f=f_{x_0,y_0}$. We note that
$$
f(1,1)=1, \quad f(1,0)=f(0,1)=f(0,0)=0, \quad f(0,v)=0,\; f(u,0)=0,
$$
for all $u,v \in [0,1]$. Hence,
\begin{align*}
\Theta(f)(t_1,t_2) =
&\ \operatorname{Re} \Bigg\{
e^{-i t_1} \Big[e^{-i t_2} + i t_2 \int_{y_0}^1 e^{-i v t_2} \, dv \Big] \\
&\qquad + i t_1 \int_{x_0}^1 e^{-i u t_1} \Big[e^{-i t_2} + i t_2 \int_{y_0}^1 e^{-i v t_2} \, dv \Big] du
\Bigg\}.
\end{align*}
Since
$$
i t_2 \int_{y_0}^1 e^{-i v t_2} \, dv = e^{-i y_0 t_2} - e^{-i t_2},
$$
we obtain
\begin{align*}
\Theta(f)(t_1,t_2)
&= \operatorname{Re} \Bigg\{
e^{-i t_1} e^{-i y_0 t_2}
+ i t_1 \int_{x_0}^1 e^{-i u t_1} e^{-i y_0 t_2} \, du
\Bigg\} \\
&= \operatorname{Re} \Bigg\{
e^{-i y_0 t_2} \left[ e^{-i t_1}
+ i t_1 \int_{x_0}^1 e^{-i u t_1} \, du \right]
\Bigg\}.
\end{align*}
Using that
$$
i t_1 \int_{x_0}^1 e^{-i u t_1} \, du = e^{-i x_0 t_1} - e^{-i t_1},
$$
we conclude that
$$
\Theta(f_{x_0,y_0})(t_1,t_2)
= \operatorname{Re}\big( e^{-i(t_1 x_0 + t_2 y_0)} \big).
$$

\medskip

We now prove (i). Let $f \in L^\infty([0,1]^2,\R^2)$ and $(t_1,t_2)\in[-a,a]^2$. Using the definition of $\Theta$ and the bound $|f|\leq \|f\|_\infty$, we obtain
$$
|\Theta(f)(t_1,t_2)|
\leq C\, \|f\|_{\infty},
$$
where $C>0$ is a constant depending only on $a$. Similarly, differentiating under the integral sign and estimating each term, we obtain
\begin{align*}
|\partial_{t_1}\Theta(f)(t_1,t_2)| &\leq C_1\, \|f\|_{\infty},\\
|\partial_{t_2}\Theta(f)(t_1,t_2)| &\leq C_2\, \|f\|_{\infty}.
\end{align*}
Taking the supremum over $(t_1,t_2)\in[-a,a]^2$, we conclude that
\begin{align*}
\|\Theta(f)\|_{\infty} &\leq K_0 \|f\|_{\infty},\\
\|\partial_{t_i}\Theta(f)\|_{\infty} &\leq K_i \|f\|_{\infty}, \quad i=1,2.
\end{align*}
The proof for the other cases is left to the reader.
\end{proof}

From the expression of $\tilde{S}_m$ and the integration by parts formula for Wiener fields, we derive the following corollary.

\begin{corollary}\label{P-m-Gwiththeta}
The following identities hold:
\begin{align*}
G &= \Theta(B),\\
P_m &= \Theta(\tilde{S}_m), \quad \forall m \in \mathbb{N}.
\end{align*}
\end{corollary}

\begin{proof}
Using the definition of $\Theta$, the fact that $B$ vanishes on the axes, and the integration by parts formula for Wiener fields, we obtain $G=\Theta(B)$.

Using that $\tilde{S}_m$ can be written as a linear combination of indicator functions
$\mathbf{1}_{[k m^{-1},1]\times [l m^{-1},1]}$
and point~(ii) of Lemma~\ref{lem-def-theta}, we conclude that
$P_m=\Theta(\tilde{S}_m)$.
\end{proof}

The representations introduced above allow us to relate the random trigonometric polynomials and their Gaussian limit to a functional version of a Donsker-type theorem for random fields.

\section{Rate of convergence in Donsker like theorem for fields}\label{rat_dons}

In this section, we state and prove a quantitative version of the Donsker-type convergence for random fields introduced above.

\subsection{Framework and results}\label{FandR}
In order to describe the limiting object in the Donsker-type theorem, we recall the definition of a Wiener field. Let $a>0$. Following Lemma~2 in~\cite{CDeo}, a Wiener process on $[0,a]^2$
vanishing on the coordinate axes is a centered Gaussian process
$(W_t)_{t\in[0,a]^2}$
with covariance function
$$
\mathbb{E}[W_s W_t]
=
\min(s_1,t_1)\min(s_2,t_2),
\quad
s=(s_1,s_2),\ t=(t_1,t_2)\in[0,a]^2.
$$
In the sequel, we set $B=(B_1,B_2)$, where $B_i$, $i=1,2$, are two independent
copies of a Wiener process on $[0,1]^2$ vanishing on the coordinate axes.

\medskip

We now introduce the empirical field. Let $\{X_{k,l} : k,l\in\mathbb{N}\}$ be a family of centered, independent, and
identically distributed random vectors in $\mathbb{R}^2$, defined on a probability
space $(\Omega,\mathcal{F},\mathbb{P})$. We consider the empirical field $\tilde{S}_m$ introduced in Section~\ref{link:donsk}. In order to emphasize its partial-sum structure, we write, for $t=(t_1,t_2)\in[0,1]^2$,
\begin{equation}
\tilde{S}_m(t)
=
\frac{1}{m}
\sum_{\substack{1\le k\le m t_1\\ 1\le l\le m t_2}}
X_{k,l},
\end{equation}
where the sums are taken over integers $k,l$. Since the sample paths of $\tilde{S}_m$ are piecewise constant, we define
$S_m$ as its linear interpolation:
$$
S_m(t)
:=
\sum_{k=0}^{m-1}
\sum_{l=0}^{m-1}
X_{k+1,l+1}\,h_k^m(t_1)\,h_l^m(t_2),
$$
where
$$
h_k^m(y)
=
\sqrt{m}
\int_{0}^{y}\mathbf{1}_{\{y\ge 0\}}
\mathbf{1}_{\left[\frac{k}{m},\frac{k+1}{m}\right]}(u)\,du,
\qquad
y\in\mathbb{R}.
$$

\medskip

We now introduce the function spaces that will be used throughout the paper.

Let $D\subset\mathbb{R}^2$. We denote
$$
L^{\infty}(D)
=
\{f:D\to\mathbb{R}^2,\ \text{measurable functions such that }\|f\|_{\infty}<\infty\},
$$
where $
\|f\|_{\infty}=\sup_{t\in D}|f(t)|.$

The set of continuous functions $f:D\to\mathbb{R}^2$ is denoted by $C^0(D)$.

Let $\alpha\in(0,1]$. Let $(E,|\cdot|_E)$ and $(F,|\cdot|_F)$ be two metric spaces. The set of $\alpha$-H\"older continuous functions is defined by
$$
C^\alpha(E,F)
:=
\left\{
f:E\to F\ \big|\ \|f\|_{\alpha,E,F}<\infty
\right\}.
$$
Here, the norm $\|\cdot\|_{\alpha,E,F}$ is given by
$$
\|f\|_{\alpha,E,F}
=
\|f\|_{\infty}
+
\sup_{\substack{s,t\in E\\ s\neq t}}
\frac{|f(s)-f(t)|_F}{|s-t|_E^\alpha}.
$$

In the case $E\subset\mathbb{R}^2$, we endow $E$ with the norm
$$
|s|:=|s_1|+|s_2|,\quad s=(s_1,s_2)\in E.
$$

When the context is clear, we omit the dependence on $E$ and $F$.

\medskip

It is well known that $(\tilde{S}_m)_m$ converges in distribution to $B$ in suitable function spaces, see, for instance, \cite{CDeo}. We are interested in quantifying the rate of this convergence. To this end, we work in the framework of the Wasserstein-1 (or Kantorovich--Rubinstein) distance (see, e.g., \cite{Decreus20}).

\medskip

Let $(E,|\cdot|_E)$ be a normed vector space.

\begin{hyp}\label{hypxp-1}
Let $p>2$. The family $\{X_{k,l}: k,l\in\mathbb{N}\}$ consists of i.i.d.\ centered random vectors in $\mathbb{R}^2$ with unit covariance matrix such that
$$
\mathbb{E}\big[\|X_{1,1}\|^p\big]<\infty.
$$
\end{hyp}

\medskip

In order to apply these convergence results to the study of level sets of random trigonometric polynomials, we introduce a slight variant of the Wasserstein-1 distance.

\begin{definition}\label{def-local-lip-1}
A function $f:E\to\mathbb{R}$ is said to be $C$ linearly locally Lipschitz continuous if
$$
|f(x+h)-f(x)|\le C\,(|x|_E+|h|_E+1)\,|h|_E, \quad \forall x,h\in E.
$$
The set of such functions with constant $1$ is denoted by $\mathrm{Lip}_{1,1}(E)$.
\end{definition}

We now state the quantitative Donsker-type convergence result that will play a key role in the proof of the main theorem.

\begin{theorem}\label{thm-Stein-filed-lip11}
Let $p>4$, $0 \leq \alpha < \frac{1}{2} - \frac{2}{p}$, and 
$\gamma < \frac{1-2\alpha-\frac{4}{p}}{5-\frac{4}{p}}$. 
Then there exists a constant $C_{p,\alpha,\gamma}>0$ such that for every family $\{X_{k,l}: k,l\in\mathbb{N}\}$ satisfying Hypothesis~\ref{hypxp-1},
\begin{align*}
\sup_{f\in\mathrm{Lip}_{1,1}(C^{\alpha}([0,1]^2,\mathbb{R}^2))}
\big|\mathbb{E}[f(S_m)]-\mathbb{E}[f(B)]\big|
\leq C_{p,\alpha,\gamma}
\frac{1}{m^\gamma}\,\mathbb{E}\big[\|X_{1,1}\|^p\big]^{4/p},
\quad \forall m\in\mathbb{N}.
\end{align*}
\end{theorem}

As a direct consequence of Theorem~\ref{thm-Stein-filed-lip11}, we obtain the following corollary.

\begin{corollary}\label{cor:stein-bounded}
Let $p>4$ and $0<\gamma < \frac{1-\frac{4}{p}}{5-\frac{4}{p}}$.
Then there exists a constant $C_{p,\gamma}>0$ such that for every family $\{X_{k,l}: k,l\in\mathbb{N}\}$ satisfying Hypothesis~\ref{hypxp-1},
\begin{align*}
\sup_{f\in\mathrm{Lip}_{1,1}(L^{\infty}([0,1]^2,\mathbb{R}^2))}
\left|
\mathbb{E}[f(\tilde{S}_m)]-\mathbb{E}[f(B)]
\right|
\leq C_{p,\gamma} \frac{1}{m^\gamma}\,\mathbb{E}\big[\|X_{1,1}\|^p\big]^{4/p},
\quad \forall m\in\mathbb{N}.
\end{align*}
\end{corollary}

Using Corollary~\ref{cor:stein-bounded}, we obtain the following result.

\begin{corollary}\label{cor-P-m-cv-G}
Let $p>4$, $0<\gamma < \frac{1-\frac{4}{p}}{5-\frac{4}{p}}$, and let $a>0$.
Then there exists a constant $C_{p,\gamma}>0$ such that for every family $\{X_{k,l}: k,l\in\mathbb{N}\}$ satisfying Hypothesis~\ref{hypxp-1},
\begin{align*}
\sup_{f\in\mathrm{Lip}_{1,1}(C^2([-a,a]^2,\mathbb{R}))}
\left|
\mathbb{E}[f(P_m)]-\mathbb{E}[f(G)]
\right|
\leq C_{p,\gamma} \frac{1}{m^\gamma}\,\mathbb{E}\big[\|X_{1,1}\|^p\big]^{4/p},
\quad \forall m\in\mathbb{N}.
\end{align*}
\end{corollary}

The proof of Theorem~\ref{thm-Stein-filed-lip11} follows the same strategy as in the one-dimensional setting and in related works such as \cite{Decreus20, lclp}. For completeness, we provide below a sketch of the main arguments.

\subsection{Comments and proof}

The argument relies on a functional approach based on fractional Sobolev spaces and their continuous embedding into H\"older spaces.
\medskip

We first recall that such embeddings hold on suitable domains. According to \cite[Definition~4.2.3]{triebelH}, the embedding 
$$
W^{\eta,p}(D)\hookrightarrow C^\alpha(D)
$$
is continuous when $D\subset\mathbb{R}^2$ is a bounded domain of cone-type. This class includes, for instance, bounded domains with $C^1$ boundary as well as open cubes.

\medskip

A bounded domain $D$ is said to be of cone-type if there exist domains $U_1,\dots,U_m$ and cones $C_1,\dots,C_m$ such that
$$
(U_j\cap D)+C_j\subset D,\quad j=1,\dots,m,
$$
where
$$
A_1+A_2=\{a_1+a_2:\ a_i\in A_i,\ i=1,2\}.
$$

\medskip

We next recall the definition of fractional Sobolev spaces (see \cite{andreev_lang}). For $\eta>0$ and $p\ge 1$, we define
\begin{equation*}
W^{\eta,p}(D)
:=
\left\{
u:D\to\mathbb{R}\ \big|\ \|u\|_{W^{\eta,p}}<\infty
\right\},
\end{equation*}
where
\begin{equation*}
\|u\|_{W^{\eta,p}}^p
=
\int_D |u(x)|^p\,dx
+
\iint_{D\times D}
\frac{|u(t)-u(s)|^p}{|t-s|^{2+p\eta}}\,dt\,ds.
\end{equation*}

The following embedding result will be used.

\begin{theorem}\label{theoem-sob-hol-embd}
Let $D\subset\mathbb{R}^2$ be a bounded domain of cone-type, let $p\in(1,\infty)$, and let $\alpha\ge 0$. 
If $\eta>\alpha+\frac{2}{p}$, then the embedding
$$
W^{\eta,p}(D)\hookrightarrow C^\alpha(D)
$$
is continuous.
\end{theorem}

\medskip

As a key intermediate step, we establish the following convergence result in Sobolev spaces.

\begin{theorem}\label{cv-w-eta-p}
Let $p\ge 4$ and $\eta\in(0,\tfrac12)$. Then there exists a constant $C_{\eta,p}>0$ such that for every family $\{X_{k,l}: k,l\in\mathbb{N}\}$ satisfying Hypothesis~\ref{hypxp-1},
\begin{align*}
\sup_{f\in\mathrm{Lip}_{1,1}(W^{\eta,p}([0,1]^2,\mathbb{R}^2))}
\left|
\mathbb{E}[f(S_m)]-\mathbb{E}[f(B)]
\right|
\leq C_{\eta,p}\,
\frac{\ln m}{m^{\frac{1-2\eta}{5}}}
\,\mathbb{E}\big[\|X_{1,1}\|^p\big]^{4/p},
\quad \forall m\in\mathbb{N}.
\end{align*}
\end{theorem}

This result will be used to derive Theorem~\ref{thm-Stein-filed-lip11} via the Sobolev--Hölder embedding.

\medskip

\begin{remark}
The restriction $p\ge 4$ is due to the use of linearly locally Lipschitz test functions. If one considers globally Lipschitz test functions, the condition $p\ge 3$ is sufficient. 
\end{remark}

\begin{remark}
In Section~\ref{crofton:error}, we will work with the set $[0,1]^2$. Although this set is not itself a domain of cone-type, it is contained in the open set $(-2,2)^2$, which satisfies this property. This allows us to apply the embedding theorem in our setting.
\end{remark}

\paragraph{Sketch of the proof of Theorem~\ref{thm-Stein-filed-lip11}.}

First, we extend the definition of $S_m$ and $B$ to $(-2,2)^2$.

\medskip

Second, we prove that $S_m$ and $B$ belong to $W^{\eta,p}((-2,2)^2,\mathbb{R}^2)$ for suitable $\eta$. The Sobolev embedding then implies that they also belong to $C^\alpha([0,1]^2,\mathbb{R}^2)$ for appropriate $\alpha\in[0,1)$. Moreover, their norms are uniformly bounded in $m$.

\medskip

Third, for $N$ to be fixed later, we introduce, for $v:[0,1]^2\to\mathbb{R}^2$, the linear interpolation $\Pi_N(v)$ on the uniform grid of mesh $N^{-1}$ on $(-2,2)^2$.

\medskip

We then prove that
\begin{align*}
\sup_{N,m}
N^{\frac{1}{2}-\eta}
\left(
\mathbb{E}\big[\|S_m-\Pi_N(S_m)\|_{\eta,p}^p\big]
+
\mathbb{E}\big[\|B-\Pi_N(B)\|_{\eta,p}^p\big]
\right)
<\infty.
\end{align*}

\medskip

With this estimate at hand, for $f\in\mathrm{Lip}_{1,1}(W^{\eta,p}((-2,2)^2,\mathbb{R}^2))$, we write
$$
\mathbb{E}[f(S_m)]-\mathbb{E}[f(B)]
=
I_1+I_2+I_3,
$$
where
\begin{align*}
I_1&=\mathbb{E}[f(S_m)]-\mathbb{E}[f(\Pi_N(S_m))],\\
I_2&=\mathbb{E}[f(\Pi_N(S_m))]-\mathbb{E}[f(\Pi_N(\Pi_m(B)))],\\
I_3&=\mathbb{E}[f(\Pi_N(\Pi_m(B)))]-\mathbb{E}[f(B)].
\end{align*}

\medskip

The terms $I_1$, $I_2$, and $I_3$ have different natures.

Using that $f\in\mathrm{Lip}_{1,1}(C^\alpha([0,1]^2,\mathbb{R}^2))$, we obtain
\begin{align*}
I_1
&\le
\mathbb{E}\Big[
\big(
1+\|S_m\|_{\alpha,[0,1]^2}+\|\Pi_N(S_m)\|_{\alpha,[0,1]^2}
\big)
\,
\|\Pi_N(S_m)-S_m\|_{\alpha,[0,1]^2}
\Big]
\\
&\le
\sqrt{
\mathbb{E}\Big[
\big(
1+\|S_m\|_{\alpha,[0,1]^2}+\|\Pi_N(S_m)\|_{\alpha,[0,1]^2}
\big)^2
\Big]
}
\,
\sqrt{
\mathbb{E}\big[
\|\Pi_N(S_m)-S_m\|_{\alpha,[0,1]^2}^2
\big]
}.
\end{align*}

Similarly,
$$
I_3
\le
\sqrt{
\mathbb{E}\Big[
\big(
1+\|B\|_{\alpha,[0,1]^2}+\|\Pi_N(B)\|_{\alpha,[0,1]^2}
\big)^2
\Big]
}
\,
\sqrt{
\mathbb{E}\big[
\|\Pi_N(B)-B\|_{\alpha,[0,1]^2}^2
\big]
}.
$$

Therefore, $I_1+I_3$ is controlled using the interpolation error estimates.

\medskip

The term $I_2$ is handled using a finite-dimensional Stein method in $\mathbb{R}^N$, recalled in Appendix~\ref{appendix-stein}.
\medskip

The proof of Theorem~\ref{thm-Stein-filed-lip11} then follows by choosing $N$ appropriately as a function of $m$.

\medskip

\subsubsection{Proof of Theorem~\ref{thm-Stein-filed-lip11}}

We first extend the definition of the Wiener field to $(-2,2)^2$.

For every $(\varepsilon_1,\varepsilon_2)\in \{\pm 1\}$, let $(W^{\varepsilon_1,\varepsilon_2}_t)_{t\in [0,2]^2}$ be four independent centered Gaussian processes defined on $(\Omega,\mathcal{F},\mathbb{P})$, all with the same covariance function
\begin{align*}
\mathbb{E}[W^{\varepsilon_1,\varepsilon_2}_s W^{\varepsilon_1,\varepsilon_2}_t]
=
\min(s_1,t_1)\min(s_2,t_2),
\quad
s,t\in[0,2]^2.
\end{align*}

We then define a Gaussian process $W=(W_t)_{t\in(-2,2)^2}$ by
\begin{align*}
W_t
=
W^{\mathrm{sg}(t_1),\mathrm{sg}(t_2)}\big(|t_1|,|t_2|\big),
\quad t=(t_1,t_2)\in(-2,2)^2,
\end{align*}
where, for $a\in\mathbb{R}$,
\begin{align*}
\mathrm{sg}(a)=
\begin{cases}
1,& a\ge0,\\
-1,& a<0.
\end{cases}
\end{align*}

The restriction of $W$ to $[0,1]^2$ is a Wiener process vanishing on the coordinate axes.

In the sequel, we set $B=(B^1,B^2)$, where $B^i$, $i=1,2$, are two independent copies of $W$.

\medskip

We now extend the definition of $S_m$ to $(-2,2)^2$ by
\begin{align}
S_m(t)
:=
\sum_{k=-2m+1}^{2m-1}
\sum_{l=-2m+1}^{2m-1}
X_{k,l}\,h_k^m(t_1)\,h_l^m(t_2),
\quad
t=(t_1,t_2)\in(-2,2)^2,
\label{sum_m}
\end{align}
where
\begin{align}
h_k^m(y)
=
\begin{cases}
\displaystyle
\sqrt{m}\int_{0}^{\max(0,y)}
\mathbf{1}_{\left[\frac{k}{m},\frac{k+1}{m}\right]}(u)\,du,
& \text{if } k\ge0,\\[1em]
\displaystyle
\sqrt{m}\int_{\min(0,y)}^{0}
\mathbf{1}_{\left[\frac{k}{m},\frac{k+1}{m}\right]}(u)\,du,
& \text{if } k<0.
\end{cases}
\label{def:h}
\end{align}

\begin{remark}
With this definition, $S_m$ vanishes on the coordinate axes, that is,
\begin{align*}
S_m(t)=0
\quad
\text{for }
t\in\{0\}\times(-2,2)\ \cup\ (-2,2)\times\{0\}.
\end{align*}
\end{remark}

\medskip

We now prove that $B$ and $S_m$ belong to suitable Sobolev and Hölder spaces.

\begin{lemma}\label{lem-B-inWetap}
Let $p>2$, $0<\eta<\frac12$, and $\alpha \in (0,\frac{1}{2})$. Then $B$ belongs to
$$
W^{\eta,p}((-2,2)^2,\mathbb{R}^2)
\quad \text{and} \quad
C^\alpha((-2,2)^2,\mathbb{R}^2).
$$
Moreover,
\begin{align*}
\mathbb{E}\big[\|B\|_\alpha^p\big]
+
\mathbb{E}\big[\|B\|_{W^{\eta,p}}^p\big]
<\infty.
\end{align*}
\end{lemma}

\begin{proof}
Let $s=(s_1,s_2)$ and $t=(t_1,t_2)$ in $(-2,2)^2$. For $i=1,2$, we have
\begin{align}
\mathbb{E}\big[
|B_s^i-B_t^i|^2
\big]
\le
2|t-s|.
\label{maj-accr-B-22}
\end{align}

Indeed, if $s_j t_j\ge0$ for $j=1,2$, then for $i=1,2$,
\begin{align*}
\mathbb{E}\big[
|B_s^i-B_t^i|^2
\big]
&=
s_1s_2+t_1t_2-2\min(s_1,t_1)\min(s_2,t_2)\\
&\le
2\big(|t_2||t_1-s_1|+|s_1||t_2-s_2|\big)
\le
2|t-s|.
\end{align*}

If $s_1t_1<0$ and $s_2t_2>0$, since $|t_1-s_1|=|t_1|+|s_1|$, we have
\begin{align*}
\mathbb{E}\big[
|B_s^i-B_t^i|^2
\big]
&=
\mathbb{E}\big[
|B_s^i-B_{(0,s_2)}^i|^2
\big]
+
\mathbb{E}\big[
|B_t^i-B_{(0,s_2)}^i|^2
\big]\\
&\le
2\big(|s_1|+|t_1|+|t_2-s_2|\big)
\le
2|t-s|.
\end{align*}

The remaining cases are treated similarly, which yields~\eqref{maj-accr-B-22}.

\medskip

We now observe that
\begin{align*}
\int_{(-2,2)^4}
(|t_1-s_1|+|t_2-s_2|)^{-2-p(\eta-\frac12)}
\,ds_1\,ds_2\,dt_1\,dt_2
<\infty,
\end{align*}
provided that $\eta<\frac12$.

\medskip

Since each component $B^i$ is a Gaussian process, it follows that
$$
B\in W^{\eta,p}((-2,2)^2,\mathbb{R}^2)
\quad \text{and} \quad
\|B\|_{W^{\eta,p}}\in L^p(\Omega).
$$

The Hölder regularity and the bound on $\mathbb{E}[\|B\|_\alpha^p]$ follow from Theorem~\ref{theoem-sob-hol-embd}.
\end{proof}

\begin{proposition}\label{SminWetap}
Let $p>2$ and $0<\eta<\frac12$. Let $\{X_{k,l}: k,l\in\mathbb{Z}\}$ be a family of random
vectors satisfying Hypothesis~\ref{hypxp-1} for $p$. Then $S_m$ belongs to
$$
W^{\eta,p}((-2,2)^2,\mathbb{R}^2),
\quad \text{for all } m\in\mathbb{N},
$$
and there exists a constant $C_{p,\eta}>0$, depending only on $p$ and $\eta$, such that
\begin{align*}
\sup_{m\in\mathbb{N}}
\mathbb{E}\big[\|S_m\|_{W^{\eta,p}}^p\big]
\le
C_{p,\eta}\,\mathbb{E}\big[\|X_{1,1}\|^p\big].
\end{align*}

Moreover, if $p>4$ and $0\le\alpha<\frac12-\frac{2}{p}$, then $S_m$ belongs to
$$
C^\alpha((-2,2)^2,\mathbb{R}^2),
$$
and there exists a constant $C_{p,\alpha}>0$, depending only on $p$ and $\alpha$, such that
\begin{align*}
\sup_{m\in\mathbb{N}}
\mathbb{E}\big[\|S_m\|_\alpha^p\big]
\le
C_{p,\alpha}\,\mathbb{E}\big[\|X_{1,1}\|^p\big].
\end{align*}
\end{proposition}

The Proposition \ref{SminWetap} is a consequence of the following Proposition.

\begin{proposition}\label{smbound}
Let $p>2$. Let $\{X_{k,l}: k,l\in\mathbb{Z}\}$ be a family of random vectors satisfying Hypothesis~\ref{hypxp-1}. Then there exists a constant $C_p>0$ such that
\begin{align*}
\sup_{m\in\mathbb{N}}
\mathbb{E}\big[
|S_m(t)-S_m(s)|^p
\big]
\le
C_p\,\mathbb{E}\big[\|X_{1,1}\|^p\big]\,
|t-s|^{p/2},
\quad
\text{for all } s,t\in(-2,2)^2.
\end{align*}
\end{proposition}

\begin{proof}
Recall that
\begin{align*}
S_m(t)
=
\sum_{k=-2m+1}^{2m-1}
\sum_{l=-2m+1}^{2m-1}
X_{k,l}\, h_k^m(t_1)\, h_l^m(t_2).
\end{align*}

Let $s,t \in (-2,2)^2$. Then
\begin{align*}
S_m(t) - S_m(s)
&=
\sum_{k,l}
X_{k,l} \bigl(h_k^m(t_1) - h_k^m(s_1)\bigr) h_l^m(t_2) \\
&\quad +
\sum_{k,l}
X_{k,l} \bigl(h_l^m(t_2) - h_l^m(s_2)\bigr) h_k^m(s_1).
\end{align*}

We recall that, for $a\in\mathbb{R}$, $(a)_+ := \max(a,0)$, and
\begin{align*}
h_k^m(t)
=
\begin{cases}
\displaystyle
\min\left( \frac{1}{\sqrt{m}}, \; \sqrt{m} \left(t - \frac{k}{m}\right)_+ \right), & k \geq 0, \\[10pt]
\displaystyle
\min\left( \frac{1}{\sqrt{m}}, \; \sqrt{m} \left(-t + \frac{k+1}{m}\right)_+ \right), & k < 0.
\end{cases}
\end{align*}

Moreover,
\begin{align*}
\left| h_k^m(t) - h_k^m(s) \right|
\le
\min \left( \frac{1}{\sqrt{m}}, \; \sqrt{m} |t - s| \right)
\mathbf{1}_{\left\{\frac{k}{m} \leq \max(t,s),\ \frac{k+1}{m} \geq \min(t,s)\right\}}.
\end{align*}

Therefore,
\begin{align*}
\mathbb{E} \left[ \left| S_m(t) - S_m(s) \right|^p \right] 
&\leq 2^{p-1} \mathbb{E} \left[ \left| \sum_{k,l} X_{k,l} \left( h_k^m(t_1) - h_k^m(s_1) \right) h_l^m(t_2) \right|^p \right] \\
&\quad + 2^{p-1} \mathbb{E} \left[ \left| \sum_{k,l} X_{k,l} \left( h_l^m(t_2) - h_l^m(s_2) \right) h_k^m(s_1) \right|^p \right].
\end{align*}

We follow the same lines as in the proof of Theorem~3.2 in~\cite{Decreus20}.

Define
\begin{align*}
\mathcal{B}^1_{s,t}
&:=
\left\{ (k,l) : \left(h_k^m(t_1) - h_k^m(s_1)\right) h_l^m(t_2)\neq 0 \right\}, \\
\mathcal{B}^2_{s,t}
&:=
\left\{ (k,l) : \left(h_l^m(t_2) - h_l^m(s_2)\right) h_k^m(s_1)\neq 0 \right\}.
\end{align*}

Using Lemma~4.1 in~\cite{Decreus20}, we obtain
\begin{align*}
\mathbb{E} \left[ |S_m(t) - S_m(s)|^p \right] 
&\leq 2^{p-1} \mathbb{E} \left[ \|X_{1,1}\|^p \right] \left| \mathcal{B}^1_{s,t} \right|^{\frac{p}{2} - 1} 
\sum_{(k,l) \in \mathcal{B}^1_{s,t}} \left| h_k^m(t_1) - h_k^m(s_1) \right|^p |h_l^m(t_2)|^p \\
&\quad + (\text{same term with } \mathcal{B}^2_{s,t}) \\
&=: \lambda_{s,t}^1 + \lambda_{s,t}^2.
\end{align*}

We estimate the first term on the right-hand side. The second term is handled analogously.

First, assume that there exist $k \in \{-2m+1, \ldots, 2m-1\}$ and $\varepsilon \in \{-1,0,1\}$ such that
$(t_1,s_1) \in I_k^m \times I_{k+\varepsilon}^m$.

In this case, only a finite number of indices contribute. More precisely,
\begin{align*}
|h_l^m(t_2)| &\le \frac{1}{\sqrt{m}}, \\
\left| h_{k+\varepsilon}^m(t_1) - h_{k+\varepsilon}^m(s_1) \right| &\le \sqrt{m} |t_1 - s_1|,
\end{align*}
and
\begin{align*}
h_r^m(t_1) - h_r^m(s_1) = 0
\quad \text{for } r \notin \{k-1,k,k+1\}.
\end{align*}

Therefore,
\begin{align*}
|\mathcal{B}^1_{s,t}|
\le C m,
\end{align*}
and we obtain
\begin{align*}
\lambda_{s,t}^1
\le
C_p \mathbb{E}\big[\|X_{1,1}\|^p\big] |t_1 - s_1|^{p/2}.
\end{align*}

Next, assume that $(t_1,s_1) \in I_k^m \times I_r^m$ with $|k-r|\ge 2$. In this case,
\begin{align*}
\frac{1}{m} \le |t_1 - s_1|,
\end{align*}
and one checks that
\begin{align*}
|\mathcal{B}^1_{s,t}|
\le
C m^2 |t_1 - s_1|.
\end{align*}

Using the bounds $|h_k^m|\le m^{-1/2}$ and $|h_k^m(t)-h_k^m(s)|\le m^{-1/2}$, we obtain
\begin{align*}
\lambda_{s,t}^1
\le
C_p \mathbb{E}\big[\|X_{1,1}\|^p\big] |t_1 - s_1|^{p/2}.
\end{align*}

By applying the same argument to $\lambda_{s,t}^2$, we obtain the same bound. Therefore,
\begin{align*}
\mathbb{E}\big[|S_m(t)-S_m(s)|^p\big]
\le
C_p\,\mathbb{E}\big[\|X_{1,1}\|^p\big]\,|t-s|^{p/2},
\end{align*}
which completes the proof.
\end{proof}

We now define the interpolation operator $\Pi_N$.

Let $N\in\mathbb{N}$. Let $f:(-2,2)^2\to\mathbb{R}$ be a function vanishing on the coordinate axes, that is,
$$
f(t_1,0)=f(0,t_2)=0.
$$

We define
\begin{align}\label{pi1}
\Pi_{N}(f)(s)
:=
\sum_{k=-2N+1}^{2N-1}
\sum_{l=-2N+1}^{2N-1}
N \Delta_N^2 f(k,l) \, h_k^N(s_1) h_l^N(s_2),
\quad s=(s_1,s_2)\in(-2,2)^2,
\end{align}
where
\begin{align}\label{incr1}
\Delta_N^2 f(k,l)
:=
f\left(\tfrac{k+1}{N}, \tfrac{l+1}{N}\right)
+ f\left(\tfrac{k}{N}, \tfrac{l}{N}\right)
- f\left(\tfrac{k+1}{N}, \tfrac{l}{N}\right)
- f\left(\tfrac{k}{N}, \tfrac{l+1}{N}\right).
\end{align}

Using telescoping sums and the structure of the functions $h_k^N$, we obtain
\begin{align*}
\Pi_N(f)\left( \tfrac{k_0}{N}, \tfrac{l_0}{N} \right)
=
f\left( \tfrac{k_0}{N}, \tfrac{l_0}{N} \right).
\end{align*}

Let $D:=(-2,2)^2$. We introduce the space
\begin{align*}
\mathcal{H}^1(D)
:=
\left\{
f \ \big| \ 
f(t_1,t_2)
=
\int_{[\min(0,t_1),\max(0,t_1)] \times [\min(0,t_2),\max(0,t_2)]}
\frac{\partial^2 f}{\partial u \partial v}(u,v)\,du\,dv,
\quad
\frac{\partial^2 f}{\partial u \partial v}\in L^2(D)
\right\}.
\end{align*}

We endow $\mathcal{H}^1(D)$ with the scalar product
\begin{align*}
\langle f, g \rangle_{\mathcal{H}^1(D)}
=
\int_{D}
\frac{\partial^2 f}{\partial u \partial v}(u,v)
\frac{\partial^2 g}{\partial u \partial v}(u,v)
\,du\,dv.
\end{align*}

Let $f \in \mathcal{H}^1(D)$ with $f(0,v)=f(u,0)=0$. Then, for each $k,l\in\mathbb{Z}$,
\begin{align}
N \Delta_N^2 f(k,l)
&=
N \left[
f\left( \tfrac{k+1}{N}, \tfrac{l+1}{N} \right)
- f\left( \tfrac{k}{N}, \tfrac{l+1}{N} \right)
- f\left( \tfrac{k+1}{N}, \tfrac{l}{N} \right)
+ f\left( \tfrac{k}{N}, \tfrac{l}{N} \right)
\right] \notag \\
&=
N \int_D
\frac{\partial^2 f}{\partial u \partial v}(u,v)
\,\mathbf{1}_{[\frac{k}{N}, \frac{k+1}{N}]}(u)
\,\mathbf{1}_{[\frac{l}{N}, \frac{l+1}{N}]}(v)
\,du\,dv \label{incr} \\
&=
\left\langle f, h_k^N \otimes h_l^N \right\rangle_{\mathcal{H}^1(D)}.
\notag
\end{align}

Therefore, the operator $\Pi_N(f)$ defined in \eqref{pi1} can be rewritten as
\begin{align}\label{pi}
\Pi_N(f)(s)
:=
\sum_{k=-2N+1}^{2N-1}
\sum_{l=-2N+1}^{2N-1}
\left\langle f, h_k^N \otimes h_l^N \right\rangle_{\mathcal{H}^1(D)}
\, h_k^N(s_1)\, h_l^N(s_2).
\end{align}

If $f=(f^1,f^2):(-2,2)^2\to\mathbb{R}^2$, then
$$
\Pi_N(f) = (\Pi_N(f^1), \Pi_N(f^2)).
$$

The next step consists in controlling the norm of $\mathrm{Id}-\Pi_N$. For this purpose, we establish the following propositions.

\begin{proposition}\label{prop-id-pin-hol-mb}
Let $p>2$ and $0<\eta<\frac{1}{2}$. There exists a constant $C_{p,\eta}$ such that
\begin{align*}
\sup_{N\in\mathbb{N}} N^{\frac{1}{2}-\eta}\, \mathbb{E}\big[\|B-\Pi_N(B)\|_\alpha^p\big] \leq C_{p,\eta}. 
\end{align*}
Let $p>4$ and $0\le\alpha<\frac12-\frac{2}{p}$. Let
$$
0<\varepsilon<\frac12-\frac{2}{p}-\alpha.
$$
There exists a constant $C_{p,\varepsilon,\alpha}$ such that
$$
\sup_{N\in\mathbb{N}}
N^{p\varepsilon}\,
\mathbb{E}\big[\|B-\Pi_N(B)\|_\alpha^p\big]
\le
C_{p,\varepsilon,\alpha}.
$$
\end{proposition}

\begin{proposition}\label{prop-id-pin-hol}
Let $p>2$ and $0<\eta<\frac{1}{2}$. There exists a constant $C_{p,\eta}$ such that for every family $\{X_{k,l},\,k,l\in\mathbb{Z}\}$ satisfying Hypothesis~\ref{hypxp-1} for $p$, 
\begin{align*}
\sup_{m\in\mathbb{N}} \sup_{N\in\mathbb{N}} N^{\frac{1}{2}-\eta}\, \mathbb{E}\big[\|S_m-\Pi_N(S_m)\|_\alpha^p\big] \leq C_{p,\eta}. 
\end{align*}
Let $p>4$ and $0\le\alpha<\frac12-\frac{2}{p}$. Let
$$
0<\varepsilon<\frac12-\frac{2}{p}-\alpha.
$$
There exists a constant $C_{p,\varepsilon,\alpha}$ such that for every family $\{X_{k,l},\,k,l\in\mathbb{Z}\}$ satisfying Hypothesis~\ref{hypxp-1} for $p$,
$$
\sup_{m\in\mathbb{N}}\sup_{N\in\mathbb{N}}
N^{p\varepsilon}\,
\mathbb{E}\big[\|S_m-\Pi_N(S_m)\|_\alpha^p\big]
\le
C_{p,\varepsilon,\alpha}\,
\mathbb{E}\big[\|X_{1,1}\|^p\big].
$$
\end{proposition}

Propositions~\ref{prop-id-pin-hol-mb} and~\ref{prop-id-pin-hol} are consequences of Theorem~\ref{theoem-sob-hol-embd}, the following lemma, and Proposition~\ref{smbound}.

\begin{lemma}\label{aux-lemma}
Let $p>2$ and $0<\eta<\frac{1}{2}$. Assume that $F$ is a stochastic process such that
\begin{align*}
\E \left[ \left| F(t) - F(s)\right|^p \right] \leq |t-s|^{p/2}, \quad \forall (s,t) \in D^2.
\end{align*}
Then,
\begin{align*}
\mathbb{E}\left[ \| F - \Pi_N(F)\|_{W^{\eta,p}}^p\right] \leq C\, N^{-p\left(\frac{1}{2} - \eta\right)}.
\end{align*}
\end{lemma}

We introduce the following definition, which will be used to localize the analysis of the interpolation error.

\begin{definition}
Let $m\in\mathbb{N}$. Assume that
$$
s\in
\left(\frac{i}{m},\frac{i+1}{m}\right)
\times
\left(\frac{j}{m},\frac{j+1}{m}\right)
:= I_{i,j}^m
\subset (-2,2)^2
$$
for some $i,j\in\{-2m+1,\dots,2m-1\}$. We say that a point
$t\in(-2,2)^2$ belongs to the neighborhood of $s$ if
\begin{align*}
t\in
&I_{i,j}^m
\cup I_{i,j+1}^m
\cup I_{i,j-1}^m
\cup I_{i-1,j}^m
\cup I_{i+1,j}^m.
\end{align*}

If $s=\left(s_1,\frac{k}{m}\right)$ with $s_1\in\left(\frac{l}{m},\frac{l+1}{m}\right)$, then the neighborhood
of $s$ is
$$
V_m(s)
:=
\left[\frac{k-1}{m},\frac{k+2}{m}\right]
\times
\left[\frac{l-1}{m},\frac{l+1}{m}\right].
$$

If $s=\left(\frac{k}{m},\frac{l}{m}\right)$, then
$$
V_m(s)
:=
\left[\frac{k-1}{m},\frac{k+1}{m}\right]
\times
\left[\frac{l-1}{m},\frac{l+1}{m}\right].
$$
\end{definition}

\begin{figure}[h]
\centering
\begin{tikzpicture}[thick]

\draw[->] (0,0) -- (5,0);
\draw[->] (0,0) -- (0,5);

\draw[step=1.0cm, very thin] (1,1) grid (4,4);

\fill[gray!40] (2,2) rectangle (3,3); 
\fill[gray!40] (2,3) rectangle (3,4); 
\fill[gray!40] (2,1) rectangle (3,2); 
\fill[gray!40] (1,2) rectangle (2,3); 
\fill[gray!40] (3,2) rectangle (4,3); 

\fill (2.5,2.5) circle (2pt);
\node at (2.8,2.6) {$s$};

\node at (1,0) [below] {$\frac{i-1}{m}$};
\node at (2,0) [below] {$\frac{i}{m}$};
\node at (3,0) [below] {$\frac{i+1}{m}$};
\node at (4,0) [below] {$\frac{i+2}{m}$};

\node at (0.8,1) [left] {$\frac{j-1}{m}$};
\node at (0.8,2) [left] {$\frac{j}{m}$};
\node at (0.8,3) [left] {$\frac{j+1}{m}$};
\node at (0.8,4) [left] {$\frac{j+2}{m}$};

\end{tikzpicture}
\caption{The gray shaded region represents the neighborhood of $s$ when $s$ belongs to the interior of a cell.}
\label{neighF1}
\end{figure}
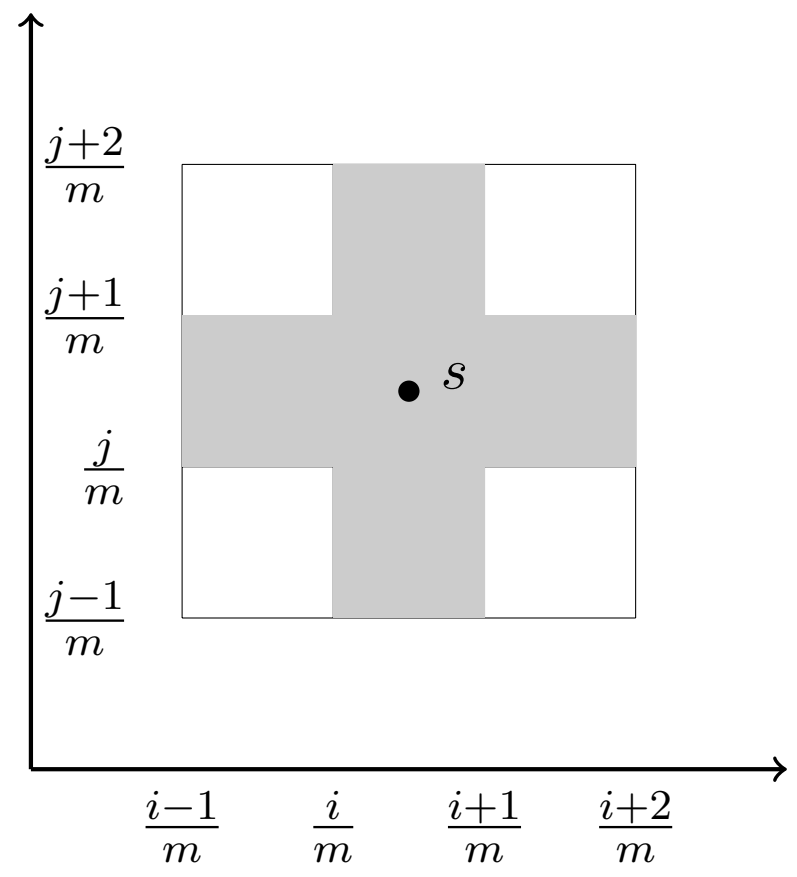

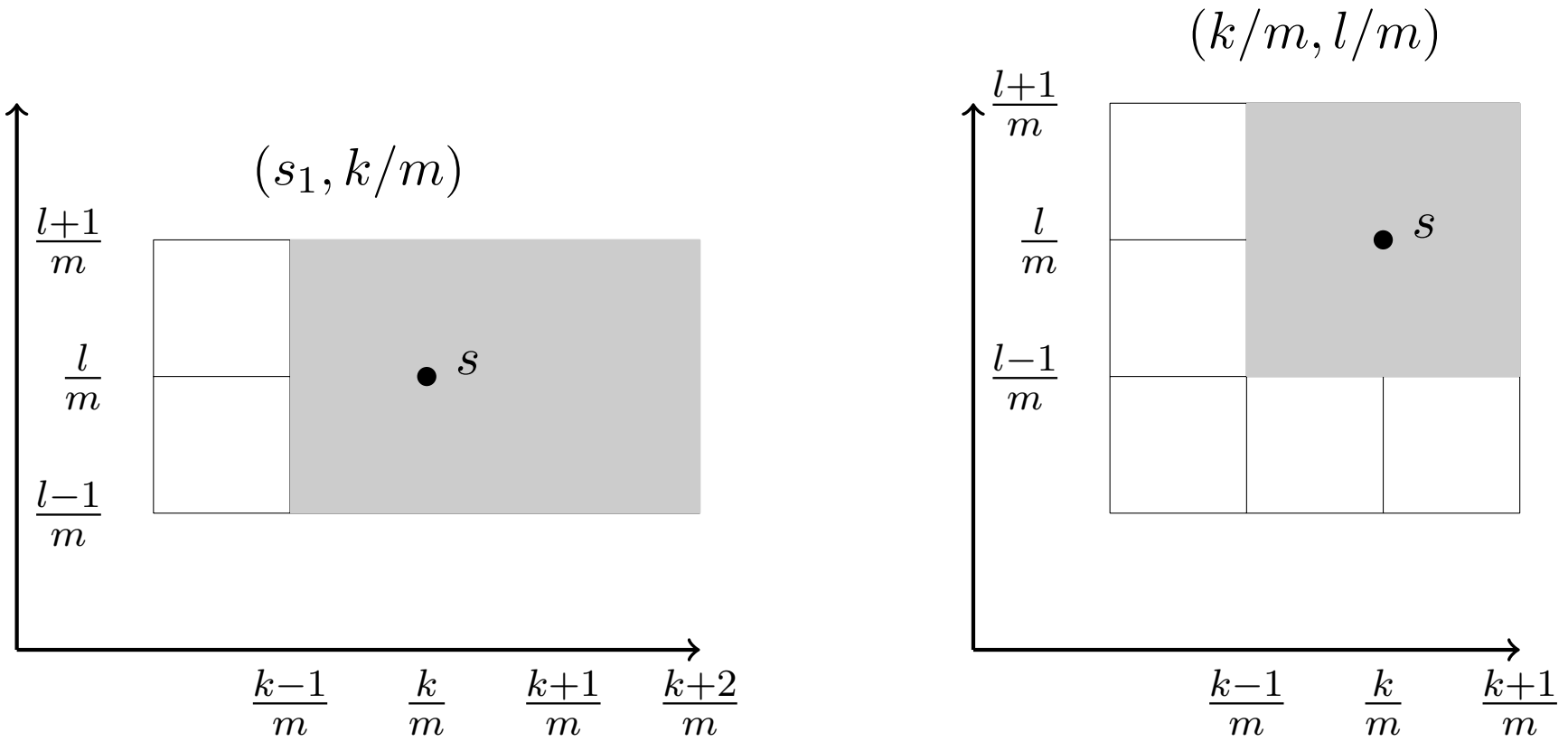
\begin{figure}[h]
\centering
\begin{tikzpicture}[thick]

\begin{scope}[shift={(0,0)}]

\draw[->] (0,0) -- (5,0);
\draw[->] (0,0) -- (0,4);

\draw[step=1.0cm, very thin] (1,1) grid (5,3);

\fill[gray!40] (2,1) rectangle (5,3);

\fill (3,2) circle (2pt);
\node at (3.3,2.1) {$s$};

\node at (2,0) [below] {$\frac{k-1}{m}$};
\node at (3,0) [below] {$\frac{k}{m}$};
\node at (4,0) [below] {$\frac{k+1}{m}$};
\node at (5,0) [below] {$\frac{k+2}{m}$};

\node at (0.8,1) [left] {$\frac{l-1}{m}$};
\node at (0.8,2) [left] {$\frac{l}{m}$};
\node at (0.8,3) [left] {$\frac{l+1}{m}$};

\node at (2.5,3.5) {$(s_1,k/m)$};

\end{scope}

\begin{scope}[shift={(7,0)}]

\draw[->] (0,0) -- (4,0);
\draw[->] (0,0) -- (0,4);

\draw[step=1.0cm, very thin] (1,1) grid (4,4);

\fill[gray!40] (2,2) rectangle (4,4);

\fill (3,3) circle (2pt);
\node at (3.3,3.1) {$s$};

\node at (2,0) [below] {$\frac{k-1}{m}$};
\node at (3,0) [below] {$\frac{k}{m}$};
\node at (4,0) [below] {$\frac{k+1}{m}$};

\node at (0.8,2) [left] {$\frac{l-1}{m}$};
\node at (0.8,3) [left] {$\frac{l}{m}$};
\node at (0.8,4) [left] {$\frac{l+1}{m}$};

\node at (2.5,4.5) {$(k/m,l/m)$};

\end{scope}

\end{tikzpicture}

\caption{Neighborhoods of $s$ in the boundary cases: edge case $s=(s_1,k/m)$ (left) and corner case $s=(k/m,l/m)$ (right).}
\label{fig:boundary-neighborhoods}
\end{figure}

These configurations are illustrated in Figures~\ref{neighF1} and~\ref{fig:boundary-neighborhoods}.

\begin{proof}
First, we prove that there exists a universal constant $C_p>0$ such that
\begin{align}\label{dif-pn-fp}
\mathbb{E} \left[ \left| \Pi_N(F)(t) - \Pi_N(F)(s) + F(s) - F(t) \right|^p \right]
\leq C_p \left( |t - s| \wedge \frac{1}{N} \right)^{\frac{p}{2}},
\quad \forall (s,t) \in D^2.
\end{align}

We split the proof into two cases. In Case~(i), $t$ belongs to a neighborhood of $s$.

In this situation, we decompose the increment along the coordinate directions as follows:
\begin{align}\label{samecoordinate}
&\Pi_N(F)(t) - \Pi_N(F)(s) + F(s) - F(t)\notag\\
&= \big(\Pi_N(F)(t_1,t_2) - \Pi_N(F)(t_1,s_2) + F(t_1,s_2) - F(t_1,t_2)\big)\notag\\
&\quad + \big(\Pi_N(F)(t_1,s_2) - \Pi_N(F)(s_1,s_2) + F(s_1,s_2) - F(t_1,s_2)\big).
\end{align}

Thus, by the triangle inequality, it suffices to consider the case where either $t_1 = s_1$ or $t_2 = s_2$. Assume, for instance, that $t_2 = s_2$.

Then, there exist $k_0,l_0 \in \mathbb{Z}$ and $\varepsilon \in \{-1,0,1\}$ such that  
$s \in I_{k_0,l_0}^N$ and $t_1 \in I_{k_0+\varepsilon}^N$.

We write
\begin{align*}
\Pi_N(F)(t) - \Pi_N(F)(s)
= \sum_{k=-2N+1}^{2N-1} \sum_{l=-2N+1}^{2N-1}
N \Delta^2_N F(k,l)\, \bigl[h_k^N(t_1) - h_k^N(s_1)\bigr] h_l^N(s_2).
\end{align*}

Using the properties of the functions $h_k^N$, only the indices $k \in \{k_0-1,k_0,k_0+1\}$ contribute. Moreover,
\begin{equation*}
\big|h_k^N(t_1) - h_k^N(s_1)\big| \le \sqrt{N}\,|t_1 - s_1|.
\end{equation*}

and therefore,
\begin{align*}
\left| \Pi_N(F)(t) - \Pi_N(F)(s) \right| 
&\leq N \sqrt{N} |t_1 - s_1| 
\sum_{k = k_0 - 1}^{k_0 + 1} 
\left| \sum_{l = -2N+1}^{2N - 1} \Delta^2_N F(k, l)\, h_l^N(s_2) \right|.
\end{align*}

We now describe the structure of the coefficients $h_l^N(s_2)$. Let $l_0$ be such that $s_2 \in I_{l_0}^N$. Then, using the definition of $h_l^N$, we have
\[
|h_l^N(s_2)| \le \frac{1}{\sqrt{N}}, \quad \forall l,
\]
and $h_l^N(s_2)=0$ unless $l$ belongs to a set of indices of cardinality of order $N|s_2|+1$.

More precisely, only the indices $l$ such that the interval $\left[\frac{l}{N},\frac{l+1}{N}\right]$ intersects the segment between $0$ and $s_2$ contribute.

Thus,
\begin{align*}
\left| \Pi_N(F)(t) - \Pi_N(F)(s) \right| 
&\leq N |t_1 - s_1|
\sum_{k = k_0 - 1}^{k_0 + 1} 
\left(
\left| \sum_{l \in \Lambda(s_2)} \Delta_N^2 F(k, l) \right|
+
\sum_{l = l_0 - 1}^{l_0 + 1} \left| \Delta_N^2 F(k, l) \right|
\right),
\end{align*}
where
\[
\Lambda(s_2)
:=
\left\{
l \in \{-2N+1,\dots,2N-1\}
:\ 
0 \leq l \cdot \operatorname{sign}(s_2) \leq |l_0|-2
\right\},
\]

Recall that
\begin{equation*}
\Delta_N^2 F(k, l) = F\left(\frac{k+1}{N}, \frac{l+1}{N}\right) + F\left(\frac{k}{N}, \frac{l}{N}\right) - F\left(\frac{k+1}{N}, \frac{l}{N}\right) - F\left(\frac{k}{N}, \frac{l+1}{N}\right).
\end{equation*}

Using telescoping sums in the second variable and the fact that $F(a,0) = F(0,b) = 0$, we obtain
\begin{equation*}
\sum_{l \in \Lambda(s_2)} \Delta_N^2 F(k, l) 
= F\left(\frac{k+1}{N}, \frac{l_0 - 1}{N}\right) - F\left(\frac{k}{N}, \frac{l_0 - 1}{N}\right).
\end{equation*}

and 
\begin{align*}
\left|\Pi_N(F)(t) - \Pi_N(F)(s)\right|  
&\leq N |t_1-s_1|
\sum_{k=k_0-1}^{k_0 +1}
\Bigg(
\left| F\left(\frac{k+1}{N}, \frac{l_0-1}{N}\right) - F\left(\frac{k}{N}, \frac{l_0-1}{N}\right)\right| \\
&\hspace{3.5cm}
+ \sum_{l=l_0-1}^{l_0+1}\left| \Delta_N^2 F(k,l)\right|
\Bigg).
\end{align*}

We take expectations and use the convexity inequality to obtain
\begin{align*}
\mathbb{E} \left[ \left|\Pi_N(F)(t) - \Pi_N(F)(s)\right| ^p\right]  
&\leq C_p N^p |t_1-s_1|^p
\sum_{k=k_0-1}^{k_0 +1}
\mathbb{E}\left[
\left| F\left(\frac{k+1}{N}, \frac{l_0-1}{N}\right) - F\left(\frac{k}{N}, \frac{l_0-1}{N}\right)\right|^p
\right] \\
&\quad + C_p N^p |t_1-s_1|^p
\sum_{k=k_0-1}^{k_0 +1}
\sum_{l=l_0-1}^{l_0+1}
\mathbb{E}\left[
\left| \Delta_N^2 F (k,l)\right|^p
\right].
\end{align*}

Applying the hypothesis on $F$, and using that $|t_1 - s_1| \leq \frac{1}{N}$, we obtain
\begin{align*}
\mathbb{E} \left[ \left| \Pi_N(F)(t) - \Pi_N(F)(s) \right|^p \right] 
&\leq C_p N^{\frac{p}{2}} |t_1 - s_1|^p \\
&= C_p N^{\frac{p}{2}} |t_1 - s_1|^{\frac{p}{2}} |t_1 - s_1|^{\frac{p}{2}} \\
&\leq C_p |t_1 - s_1|^{\frac{p}{2}},
\end{align*}
where in the last inequality we used that
\[
N^{\frac{p}{2}} |t_1 - s_1|^{\frac{p}{2}}
\leq N^{\frac{p}{2}}\left(\frac{1}{N}\right)^{\frac{p}{2}} = 1.
\]

Consequently,
\[
\mathbb{E} \left[ \left| \Pi_N(F)(t) - \Pi_N(F)(s) \right|^p \right]
\leq C_p \left( |t_1 - s_1| \wedge \frac{1}{N} \right)^{\frac{p}{2}}.
\]

Using once again the assumption on $F$, we have
\[
\mathbb{E} \left[ |F(s) - F(t)|^p \right] 
\leq C_p |t - s|^{\frac{p}{2}}
\leq C_p \left( |t - s| \wedge \frac{1}{N} \right)^{\frac{p}{2}}.
\]

By the convexity inequality, there exists a positive constant $C_p$ such that, for $s$ and $t$ in the same neighborhood with $t_2 = s_2$, 
\[
\mathbb{E} \left[ \left| \Pi_N(F)(t) - \Pi_N(F)(s) + F(s) - F(t) \right|^p \right] 
\leq C_p \left( |t_1 - s_1| \wedge \frac{1}{N} \right)^{\frac{p}{2}}.
\]

Using \eqref{samecoordinate}, there exists a constant $C_p>0$ such that for all $s$ and $t$ in the same neighborhood,
\begin{align}\label{majdiffpiligne1}
\mathbb{E} \left[ \left| \Pi_N(F)(t) - \Pi_N(F)(s) + F(s) - F(t) \right|^p \right] 
\leq C_p \left[
\left( |t_1 - s_1| \wedge \frac{1}{N} \right)^{\frac{p}{2}} 
+ \left( |t_2 - s_2| \wedge \frac{1}{N} \right)^{\frac{p}{2}}
\right].
\end{align}

For Case (ii), we assume that $s$ and $t$ do not belong to the same neighborhood. Again, we may assume that $s_1 = t_1$ or $s_2 = t_2$.

If $s_2 = t_2$ and $s, t$ are not in the same neighborhood, then $|s_1 - t_1| > \frac{1}{N}$. We can find points $t^N$ and $s^N$ on the grid such that $t^N$ and $t$ belong to the same neighborhood, $s^N$ and $t^N$ belong to the same neighborhood, and
\begin{align*}
    &|t_1 - t_1^N| \leq N^{-1} \wedge |t_1 - s_1|, \\
    &|s_1 - s_1^N| \leq N^{-1} \wedge |t_1 - s_1|, \\
    &|t_1^N - s_1^N| \leq |t_1 - s_1|, \\
    &t_2^N = s_2^N.
\end{align*}

Since $F$ and $\Pi_N(F)$ coincide on the grid, we have
\begin{align*}
\mathbb{E} \left[ \left| \Pi_N(F)(t) - \Pi_N(F)(s) + F(s) - F(t) \right|^p \right] 
&\leq 2^{p-1} \mathbb{E} \left[ \left| \Pi_N(F)(t) - \Pi_N(F)(t^N) + F(t^N) - F(t) \right|^p \right] \\
&\quad + 2^{p-1} \mathbb{E} \left[ \left| \Pi_N(F)(s^N) - \Pi_N(F)(s) + F(s) - F(s^N) \right|^p \right].
\end{align*}

Applying estimate \eqref{majdiffpiligne1} to the pairs $(t, t^N)$ and $(s, s^N)$, we obtain
\[
\mathbb{E} \left[ \left| \Pi_N(F)(t) - \Pi_N(F)(s) + F(s) - F(t) \right|^p \right] 
\leq C_p \left( |t - s| \wedge \frac{1}{N} \right)^{\frac{p}{2}},
\]
which corresponds to \eqref{dif-pn-fp}.

Since
\[
\sup_{N \in \mathbb{N}} 
N^{p\left(\frac{1}{2} - \eta\right)} 
\iint_{D^2} \frac{
\left(|t_1 - s_1| \wedge N^{-1}\right)^{\frac{p}{2}} 
+ \left(|t_2 - s_2| \wedge N^{-1}\right)^{\frac{p}{2}} 
}{|t - s|^{2 + p \eta}} \, ds \, dt < \infty,
\]
we conclude that
\[
\mathbb{E} \left[ \| F - \Pi_N(F) \|_{W^{\eta,p}}^p \right] 
\leq C N^{-p\left(\frac{1}{2} - \eta\right)}, 
\quad \forall N \in \mathbb{N}.
\]

This completes the proof of Lemma~\ref{aux-lemma}.

\end{proof}

\medskip
\section{Rates of convergence for the nodal length}

\subsection{Crofton representation and error decomposition}\label{crofton:error}
In this section, we apply Crofton's formula to represent the nodal length through zero counts along one-dimensional sections, thereby reducing the problem to a one-dimensional framework. We then decompose the difference of the expectations into five terms, which will be estimated separately in the subsequent analysis. Crofton's formula requires that a function $f$ has no critical point on its nodal set, that is, no point $t\in B(0,a)$ such that $f(t)=0$ and $\nabla f(t)=0$. In Appendix~\ref{app-nondeg}, we show that this property holds almost surely for $G$, and with a probability tending to one as $m$ goes to infinity for $P_m$. Since the probability of this event is small, it does not affect the rate of convergence, and we may therefore apply Crofton's formula to the paths of both $P_m$ and $G$.

For $a>0$, we denote by
\begin{align*}
{\mathcal B}(f,a):= \{t \in B(0,a)\mid f(t)=0\}
\end{align*}
and by ${\mathcal L}(f,B(0,a))$ the nodal length of $f$ in $B(0,a)$, that is, the one-dimensional Hausdorff measure of ${\mathcal B}(f,a)$, provided that $f$ is sufficiently regular and non-degenerate.

Using Crofton’s formula, ${\mathcal L}(f,B(0,a))$ can be expressed in terms of the number of zeros of the restriction of $f$ to suitable lines.

\begin{lemma}\label{crofton-lemma}
Let $f:\mathbb R^2\to\mathbb R$ be a $C^1$ function such that
\[
\nabla f(t)\neq 0
\qquad\text{whenever}\qquad
f(t)=0.
\]
Then
\begin{align*}
{\mathcal L}(f,B(0,a))
=
\frac{1}{4}\int_0^{2\pi} \int_0^a
N\left(f\circ L_{\theta,\rho},I_\rho\right)
\,d\rho\,d\theta,
\end{align*}
where
\[
L_{\theta,\rho}(\beta)
:=
(\rho\sin\theta+\beta\cos\theta,
-\rho\cos\theta+\beta\sin\theta),
\qquad
\beta\in\mathbb R,
\]
\[
I_\rho
:=
[-\sqrt{a^2-\rho^2},\sqrt{a^2-\rho^2}],
\]
and, for every sufficiently regular function
$g:\mathbb R\to\mathbb R$
and every compact interval $I\subset\mathbb R$,
\[
N(g,I)
=
\#\{t\in I\mid g(t)=0\}.
\]
\end{lemma}

\begin{figure}[h]
\centering
\begin{tikzpicture}[thick]

\draw[->, gray] (-4.5,0) -- (4.5,0);
\draw[->, gray] (0,-4.5) -- (0,4.5);

\node at (4.8,0) {$t_1$};
\node at (0,4.8) {$t_2$};

\draw (0,0) circle(3);

\draw[->, blue] (0,0) -- (2.12,2.12);
\node[blue] at (1.7,2.6) {$v$};

\draw[blue] (-3.1,3.1) -- (3.1,-3.1);
\node[blue] at (2.6,-2.6) {$v^\perp$};

\filldraw[black] (1.5,-1.5) circle (1.5pt);
\node[below right] at (1.5,-1.5) {$y$};

\draw[magenta] (-1.5,-4.5) -- (4.5,1.5);
\node[magenta] at (3.8,1.8) {$L_{\theta,\rho}$};

\draw[dashed,gray] (-3,-3) -- (4,4);

\draw (0.9,0) arc (0:45:0.9);
\node at (0.7,0.35) {$\theta$};

\node at (0.55,-0.8) {$\rho$};

\node[below left] at (0,0) {$O$};

\end{tikzpicture}

\caption{Geometric parametrization used in Crofton’s formula.}
\label{crofF2}
\end{figure}
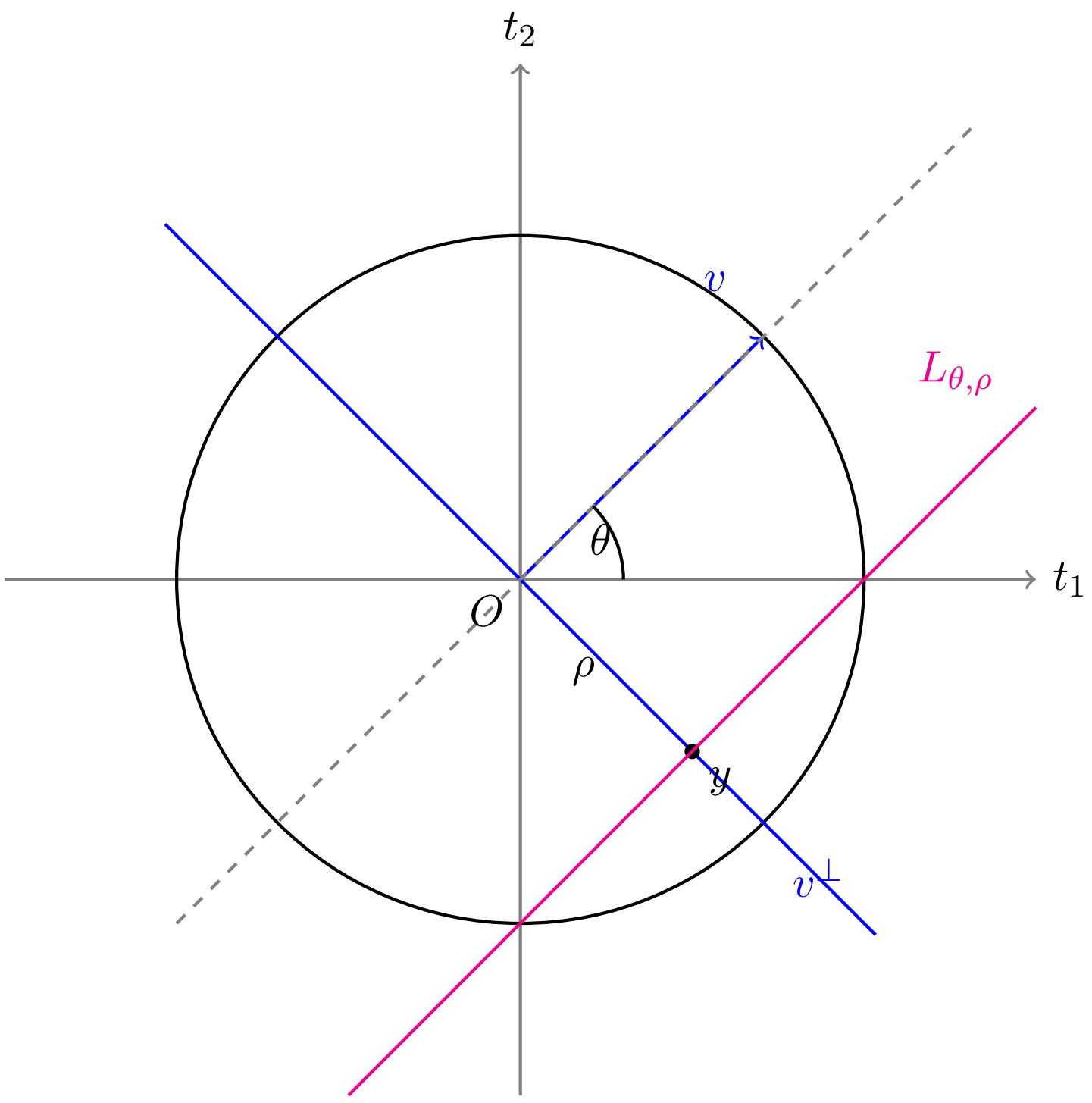

\begin{proof}
We consider the following integral-geometric representation of ${\mathcal B}(P_m,a)$, taken from \cite{leonandmas} and \cite{surfarea}:
\begin{equation}\label{coft_1}
{\mathcal L}(P_m,B(0,a))
:=
\frac{1}{\tilde{\beta}(2)}
\int_{v \in S^{1}}
\int_{y \in v^{\perp}}
\#
\left(
{\mathcal B}(P_m,a)\cap \ell_{v,y}
\right)
\, d\lambda_{1}(y)\, dv,
\end{equation}

where \(S^{1}\) is the unit sphere in \(\mathbb{R}^2\), \(dv\) denotes the arc-length measure on \(S^{1}\), \(\lambda_{1}\) is the one-dimensional Lebesgue measure, and \(\tilde{\beta}(2) = \frac{2}{\pi}\).

Let $v \in S^1$. Then there exists $\theta \in [0,2\pi)$ such that $v=(\cos\theta,\sin\theta)$. We parametrize $v^\perp$ by
\begin{equation*}
y=(\rho\sin\theta,-\rho\cos\theta),
\qquad
\rho\in[0,a]
\end{equation*}
(see Figure~\ref{crofF2}). We now consider the affine line
\begin{align*}
L_{\theta,\rho}
=
\left\{
y+\beta v
\mid
\beta\in\mathbb R
\right\}.
\end{align*}
Then
\begin{align}\label{Lthetarho}
L_{\theta,\rho}
=
\left\{
(\rho\sin\theta+\beta\cos\theta,\;
-\rho\cos\theta+\beta\sin\theta)
\mid
\beta\in\mathbb R
\right\}.
\end{align}

Since ${\mathcal B}(f,a)\subset B(0,a)$, any point
\begin{align*}
x=y+\beta v
\in
L_{\theta,\rho}\cap {\mathcal B}(f,a)
\end{align*}
must satisfy
\begin{equation*}
\|x\|^2=\rho^2+\beta^2\leq a^2.
\end{equation*}
Therefore,
\begin{equation*}
\beta\in
\left[
-\sqrt{a^2-\rho^2},
\sqrt{a^2-\rho^2}
\right]
=:I_\rho.
\end{equation*}

The lemma follows from the identity
\begin{equation*}
\#
\left(
{\mathcal B}(f,a)\cap L_{\theta,\rho}
\right)
=
N\left(
f\circ L_{\theta,\rho},
I_\rho
\right).
\end{equation*}

\end{proof}


We now use the Crofton representation to compare the nodal lengths associated with $P_m$ and $G$. To this end, {let $F:\mathbb{R}\to\mathbb{R}$ be such that $\|F\|_{\mathrm{Lip}}+\|F\|_{\infty}\leq 1$. Our goal is to} estimate
\begin{align*}
\E\left[
F\left(
{\mathcal L}(P_m,B(0,a))
\right)
\right]
-
\E\left[
F\left(
{\mathcal L}(G,B(0,a))
\right)
\right].
\end{align*}

To approximate the number of zeros along one-dimensional sections, we introduce the following functionals. For $f:I\to\mathbb R$ and $\delta>0$, define
\begin{align}
\Phi_{\delta, I}(f)
&=
\frac{1}{2\delta}
\int_I
|f'(u)|
\,\mathbf{1}_{\{|f(u)|\leq\delta\}}
\,du,
\\
\Phi_{\delta,\varepsilon,I}(f)
&=
\frac{1}{2\delta}
\int_I
|f'(u)|
\,H_{\delta,\varepsilon}(f(u))
\,du,
\end{align}
where
\begin{align}\label{def-H-delta-epsi}
H_{\delta,\varepsilon}(u)
:=
\begin{cases}
1,
& \text{if } |u|\leq\delta,
\\
0,
& \text{if } |u|\geq\delta+\varepsilon,
\\
-\dfrac{|u|-\delta}{\varepsilon}+1,
& \text{if } \delta<|u|<\delta+\varepsilon.
\end{cases}
\end{align}

Using the previous approximation functionals together with the Crofton representation, we decompose the difference of expectations into five terms. More precisely,
\begin{align*}
&\mathbb{E}
\left[
F\left(
{\mathcal L}(P_m,B(0,a))
\right)
\right]
-
\mathbb{E}
\left[
F\left(
{\mathcal L}(G,B(0,a))
\right)
\right]
\\
&=
\mathbb{E}
\left[
F\left(
{\mathcal L}(P_m,B(0,a))
\right)
\right]
-
\mathbb{E}
\left[
F\left(
\frac{1}{4}
\int_0^a
\int_0^{2\pi}
\Phi_{\delta,I_\rho}
(P_m\circ L_{\theta,\rho})
\,d\theta\,d\rho
\right)
\right]
\\
&\quad+
\mathbb{E}
\left[
F\left(
\frac{1}{4}
\int_0^a
\int_0^{2\pi}
\Phi_{\delta,I_\rho}
(P_m\circ L_{\theta,\rho})
\,d\theta\,d\rho
\right)
\right]
-
\mathbb{E}
\left[
F\left(
\frac{1}{4}
\int_0^a
\int_0^{2\pi}
\Phi_{\delta,\varepsilon,I_\rho}
(P_m\circ L_{\theta,\rho})
\,d\theta\,d\rho
\right)
\right]
\\
&\quad+
\mathbb{E}
\left[
F\left(
\frac{1}{4}
\int_0^a
\int_0^{2\pi}
\Phi_{\delta,\varepsilon,I_\rho}
(P_m\circ L_{\theta,\rho})
\,d\theta\,d\rho
\right)
\right]
-
\mathbb{E}
\left[
F\left(
\frac{1}{4}
\int_0^a
\int_0^{2\pi}
\Phi_{\delta,\varepsilon,I_\rho}
(G\circ L_{\theta,\rho})
\,d\theta\,d\rho
\right)
\right]
\\
&\quad+
\mathbb{E}
\left[
F\left(
\frac{1}{4}
\int_0^a
\int_0^{2\pi}
\Phi_{\delta,\varepsilon,I_\rho}
(G\circ L_{\theta,\rho})
\,d\theta\,d\rho
\right)
\right]
-
\mathbb{E}
\left[
F\left(
\frac{1}{4}
\int_0^a
\int_0^{2\pi}
\Phi_{\delta,I_\rho}
(G\circ L_{\theta,\rho})
\,d\theta\,d\rho
\right)
\right]
\\
&\quad+
\mathbb{E}
\left[
F\left(
\frac{1}{4}
\int_0^a
\int_0^{2\pi}
\Phi_{\delta,I_\rho}
(G\circ L_{\theta,\rho})
\,d\theta\,d\rho
\right)
\right]
-
\mathbb{E}
\left[
F\left(
{\mathcal L}(G,B(0,a))
\right)
\right]
\\
&=
\sum_{i=1}^5 I_i.
\end{align*}
\

The above decomposition yields several error terms, which we estimate separately in the following subsections.

\subsubsection{Estimation of  $I_3$}
\begin{proposition}\label{I_3}
For every $p>4$ and every $0<\gamma <\frac{1-\frac{4}{p}}{5-\frac{4}{p}}$,
there exists a constant $C_{p,\gamma}$ such that for all
$\{a_{k,l}, b_{k,l} : k,l\in {\mathbb N}\}$ satisfying Hypothesis \ref{hypxp-1} we have
\begin{align*}
I_3 \leq C_{p,\gamma}a^2(1+a)
\frac{\pi(1 + \varepsilon)}{2 C \delta \varepsilon}m^{-\gamma}.
\end{align*}
\end{proposition}

\begin{proof}
Let $F~:{\mathbb R}\rightarrow {\mathbb R}$ be Lipschitz continuous with constant $1$.

Let \( x, y \in C^1(B(0,a),{\mathbb R})\). We endow
\(C^1(B(0,a),{\mathbb R})\) with the norm
\[
\|x\|_{1,a}
=
\sup_{t\in B(0,a)}
\Big(
|x(t)|+ |\partial_1x(t)|+ |\partial_2x(t)|
\Big).
\]
Since \( F \) is Lipschitz continuous, we have
\begin{align*}
&F \left( \frac{1}{4} \int_0^{a} \int_0^{2\pi} \Phi_{\delta, \varepsilon, I_{\rho}}(x\circ L_{\theta,\rho})  \, d\theta \, d\rho \right)
- F \left( \frac{1}{4} \int_0^a  \int_0^{2\pi} \Phi_{\delta, \varepsilon, I_{\rho}}(y\circ L_{\theta,\rho})  \, d\theta \, d\rho \right)\\
&\leq \frac{1}{4} \int_0^a \int_0^{2\pi} \left| \Phi_{\delta, \varepsilon, I_{\rho}}(x\circ L_{\theta,\rho}) - \Phi_{\delta, \varepsilon, I_{\rho}}(y\circ L_{\theta,\rho})\right| \, d\theta \, d\rho.
\end{align*}

On the other hand, the function \( H_{\delta,\varepsilon} \) is bounded by \(1\) and is \(1/\varepsilon\)-Lipschitz on \(\mathbb R\). Then
\[
\frac{C \delta \varepsilon}{(1 + \varepsilon) |I_{\rho}|} \, \Phi_{\delta, \varepsilon, I_{\rho}}
\in \operatorname{Lip}_{1,1}(C^1([0,1], \mathbb{R})),
\]
and therefore,
\begin{align*}
&\left| \Phi_{\delta, \varepsilon, I_{\rho}}(x\circ L_{\theta,\rho}) - \Phi_{\delta, \varepsilon, I_{\rho}}(y\circ L_{\theta,\rho}) \right| \\
&\leq \frac{(1 + \varepsilon) |I_{\rho}|}{C \delta \varepsilon}
\|x\circ L_{\theta,\rho} - y\circ L_{\theta,\rho}\|_{1,I_{\rho}}
\left[
\|x\circ L_{\theta,\rho}\|_{1,I_{\rho}}
+
\|x\circ L_{\theta,\rho} - y\circ L_{\theta,\rho}\|_{1,I_{\rho}}
\right].
\end{align*}

Recall that, for every \(f\in C^1(\mathbb R^2,\mathbb R)\),
\begin{align*}
(f\circ L_{\theta,\rho})(\beta)
:=
f\!\left(
\rho \sin\theta + \beta \cos\theta,\;
-\rho \cos\theta + \beta \sin\theta
\right),
\qquad \forall \beta \in \mathbb{R}.
\end{align*}
Since
\begin{equation*}
\|x\circ L_{\theta,\rho}-y\circ L_{\theta,\rho}\|_{1,I_\rho}
\le
\|x-y\|_{1,B(0,a)},
\end{equation*}
recall that
$$I_{\rho}=\left[-\sqrt{a^2-\rho^2},\,\sqrt{a^2-\rho^2}\right],$$
so that
$$|I_{\rho}| \leq 2a,$$
and
\begin{align*}
    &F \left( \frac{1}{4} \int_0^a \int_0^{2\pi} \Phi_{\delta, \varepsilon, I_\rho} (x\circ L_{\theta,\rho} ) \, d\theta \, d\rho \right)
    - F \left( \frac{1}{4} \int_0^a \int_0^{2\pi} \Phi_{\delta, \varepsilon, I_{\rho}} (y \circ L_{\theta,\rho} ) \, d\theta \, d\rho \right)  \\
   &\leq \frac{\pi(1 + \varepsilon)a^2}{2 C \delta \varepsilon}
   \,  \|x- y\|_{1,B(0,a)}
   \left[
      \|x\|_{1,B(0,a)}
      +
      \|x - y\|_{1,B(0,a)}
   \right].
\end{align*}

Thus,
\[
x \mapsto
\frac{2 C \delta \varepsilon}{\pi(1 + \varepsilon)a^2}
F \left(
\frac{1}{4}
\int_0^{a}
\int_0^{2\pi}
\Phi_{\delta, \varepsilon, I_{\rho}}
(x\circ L_{\theta,\rho})
\, d\theta \, d\rho
\right)
\]
belongs to
\(
\mathrm{Lip}_{1,1}
\bigl(
C^{1}(B(0,a),{\mathbb R}),
{\mathbb R}
\bigr).
\)

According to Corollary \ref{cor-P-m-cv-G} for $p >4$ and $0<\gamma <\frac{1-\frac{4}{p}}{5-\frac{4}{p}}$
there exists a constant $C_{p,\gamma}$ such that
\begin{align*}
I_3 \leq C_{p,\gamma}a^2(1+a)
\frac{\pi(1 + \varepsilon)}{2 C \delta \varepsilon}m^{-\gamma}.
\end{align*}
\end{proof}

\subsubsection{Estimation of  $I_4$}
We follow essentially the same strategy as in the proof of Lemma 4.2 in \cite{lclp}.

Using the representation $G=\Theta(B)$, Euler's formula, and the decomposition
$B=B^1+iB^2$, where $B^1$ and $B^2$ are independent Wiener fields, we may rewrite
$G_{\theta,\rho}$ as follows:
\begin{align*}
G\circ L_{\theta, \rho}(\beta):= G_{\theta, \rho}(\beta)
&= \int_{[0,1]^2} \cos \left( \left[ \rho \sin \theta + \beta \cos \theta \right] u + \left[ -\rho \cos \theta + \beta \sin \theta \right] v \right) \, dB^1(u,v) \\
&\quad - \int_{[0,1]^2} \sin \left( \left[ \rho \sin \theta + \beta \cos \theta \right] u + \left[ -\rho \cos \theta + \beta \sin \theta \right] v \right) \, dB^2(u,v).
\end{align*}

As a consequence, we obtain the following results.

\begin{lemma}\label{lem-cov-g-derg}
The covariance matrix of the Gaussian vector
\begin{equation*}
\left(
G \circ L_{\theta,\rho}(\beta),
(G\circ L_{\theta,\rho})'(\beta)
\right)
\end{equation*}
is given by
\begin{equation*}
\Gamma=
\begin{pmatrix}
    1  &  0  \\
    0  &  \frac{1}{3}+\frac{\sin 2\theta}{4}
\end{pmatrix}.
\end{equation*}
\end{lemma}

\begin{proof}
Using the independence of $B^1$ and $B^2$, we obtain

\begin{align*}
\E[G_{\theta, \rho}^2(\beta)]
&=\E\left[\left[\int_{[0,1]^2}\cos \left(\left[\rho\sin \theta + \beta \cos \theta\right]u+\left[-\rho \cos \theta +\beta \sin \theta\right]v\right)dB^1\right]^2\right]\\
&\quad+\E\left[\left[\int_{[0,1]^2}\sin \left(\left[\rho\sin \theta + \beta \cos \theta\right]u+\left[-\rho \cos \theta +\beta \sin \theta\right]v\right)dB^2\right]^2\right]\\
&=\int_{[0,1]^2}\cos^2 \left(\left[\rho\sin \theta + \beta \cos \theta\right]u+\left[-\rho \cos \theta +\beta \sin \theta\right]v\right)\,du\,dv\\
&\quad+\int_{[0,1]^2}\sin^2 \left(\left[\rho\sin \theta + \beta \cos \theta\right]u+\left[-\rho \cos \theta +\beta \sin \theta\right]v\right)\,du\,dv\\
&=\int_{[0,1]^2}du\,dv\\
&=1.
\end{align*}

Similarly,
\begin{align*}
 (G\circ L_{\theta,\rho})'(\beta)
&= -\int_{[0,1]^2}(u \cos \theta + v\sin \theta)  \sin \left( \left[ \rho \sin \theta + \beta \cos \theta \right] u + \left[ -\rho \cos \theta + \beta \sin \theta \right] v \right) \, dB^1(u,v) \\
&\quad - \int_{[0,1]^2} (u \cos \theta + v\sin \theta) \cos \left( \left[ \rho \sin \theta + \beta \cos \theta \right] u + \left[ -\rho \cos \theta + \beta \sin \theta \right] v \right) \, dB^2(u,v).
\end{align*}

Therefore,
\begin{align*}
\E\!\left[(G\circ L_{\theta,\rho})'(\beta)^2\right]
&=\int_{[0,1]^2}\left(u\cos \theta+v\sin \theta\right)^2\,du\,dv\\
&=\int_{[0,1]^2}(u^2\cos^2\theta +v^2\sin^2 \theta +2uv\cos \theta \sin \theta ) \,du\,dv\\
&=\frac{1}{3}+\frac{\sin 2\theta}{4}.
\end{align*}

Finally,
\begin{align*}
\E\!\left[G_{\theta,\rho}(\beta)\,(G\circ L_{\theta,\rho})'(\beta)\right]=0.
\end{align*}
\end{proof}

\begin{lemma}\label{exp_app}
For the random vector $\left(G_{\theta, \rho}(\beta),(G\circ L_{\theta,\rho})'(\beta)\right)$ we have
\[
\E \left( \left|(G\circ L_{\theta,\rho})'(\beta)\right|
\left| \frac{\left|G_{\theta, \rho}(\beta)\right|-\delta-\varepsilon}{\varepsilon} \right|
1_{\delta\leq \left|G_{\theta, \rho}(\beta)\right|\leq \delta+\varepsilon }\right)
\leq \frac{7\sqrt{3}\varepsilon}{6 \pi}.
\]
\end{lemma}

\begin{proof}
The density of the Gaussian vector
\[
\left(G_{\theta,\rho}(\beta),(G\circ L_{\theta,\rho})'(\beta)\right)
\]
is given by
\begin{align}
h(x,y)
&=
\frac{1}{2\pi |\Gamma|^{1/2}}
\exp\left(
-\frac{1}{2}
(x,y)\Gamma^{-1}(x,y)^{\top}
\right).
\end{align}

Since
\[
\frac{1}{12}
\leq
\frac{1}{3}+\frac{\sin 2\theta}{4}
\leq
\frac{7}{12},
\]
we have
\begin{align*}
h(x,y)
\leq
\frac{\sqrt{3}}{\pi}
e^{- \frac{x^2}{2} - \frac{6y^2}{7}}.
\end{align*}

Then
\begin{align*}
&\E \left(
\left|(G\circ L_{\theta,\rho})'(\beta)\right|
\left(
\frac{\left|G_{\theta, \rho}(\beta)\right|-\delta-\varepsilon}{\varepsilon}
\right)
1_{\delta\leq \left|G_{\theta, \rho}(\beta)\right|\leq \delta+\varepsilon }
\right)\\
&=\int_{\R^2}|y|
\left|\frac{|x|-\delta-\varepsilon}{\varepsilon}\right|
1_{\delta \leq |x|\leq \delta+\varepsilon}
h(x,y)\,dx\,dy\\
&\leq
\int_{\R^2}|y|
\left|\frac{|x|-\delta-\varepsilon}{\varepsilon}\right|
1_{\delta \leq |x|\leq \delta+\varepsilon}
\frac{\sqrt{3}}{\pi}
\exp\left(-\frac{x^2}{2}\right)
\exp\left(-\frac{6y^2}{7}\right)\,dx\,dy\\
&\leq
\frac{2\sqrt{3}}{\varepsilon \pi}
\int_{\R}\int_{\delta}^{\delta + \varepsilon}
|y|\left||x|-\delta-\varepsilon\right|
\exp\left(-\frac{6y^2}{7}\right)\,dx\,dy\\
&\approx \frac{7\sqrt{3}\varepsilon}{ 6\pi}.
\end{align*}
\end{proof}

Therefore, we can estimate the term \(I_4\), defined by

\begin{align*}
I_4=\mathbb{E} \left[ F \left( \frac{1}{4} \int_0^a\int_0^{2\pi} \Phi_{\delta, \varepsilon, I_\rho} (G_{\theta, \rho} (\beta)) d \theta d \rho \right) \right] - \mathbb{E} \left[ F \left( \frac{1}{4}\int_0^a\int_0^{2\pi} \Phi_{\delta, I_{\rho}} (G_{\theta, \rho} (\beta)) d \theta d \rho \right) \right].
\end{align*}

\begin{proposition}\label{I_4}
The following estimate holds:
\begin{equation*}
I_4 \leq \frac{7\sqrt{3}\varepsilon a^2}{24 \delta}.
\end{equation*}
\end{proposition}

\begin{proof}
Since \(F\) is Lipschitz continuous, Lemma \ref{exp_app} yields

\begin{align}\notag
I_4
&\leq \frac{1}{4} \int_0^a \int_0^{2\pi}\E\left[ \left| \Phi_{\delta, \varepsilon, I_{\rho}} (G_{\theta, \rho}(\beta)) - \Phi_{\delta, I_{\rho}} (G_{\theta, \rho}(\beta)) \right|\right] \, d\theta \, d\rho\\\notag
&\leq \frac{1}{4} \int_0^a \int_0^{2\pi} \frac{1}{2 \delta} \int_{I_{\rho}} \E \left[ \left|(G\circ L_{\theta,\rho})'(\beta)\right| \left| H_{\delta, \varepsilon} (G_{\theta, \rho}(\beta)) - 1_{\{ \left |G_{\theta, \rho}(\beta)\right| \leq \delta \}} \right|\right] d\beta \, d\theta \, d\rho\\\notag
&\leq \frac{1}{4} \int_0^a \int_0^{2\pi} \frac{1}{2 \delta} \int_{I_{\rho}} \E \left[ \left|(G\circ L_{\theta,\rho})'(\beta)\right| \left| \frac{\left|G_{\theta, \rho}(\beta)\right|-\delta-\varepsilon}{\varepsilon} \right|1_{\delta\leq \left|G_{\theta, \rho}(\beta)\right|\leq \delta+\varepsilon }\right] d\beta \, d\theta \, d\rho\\\notag
&\leq \frac{1}{4} \int_0^a \int_0^{2\pi} \frac{1}{2 \delta} \int_{I_{\rho}}\frac{7\sqrt{3}\varepsilon}{6 \pi} d\beta \, d\theta \, d\rho = \frac{7\sqrt{3}\varepsilon a^2}{24 \delta}.
\end{align}
\end{proof}


\subsubsection{Estimation of  $I_5$}

We recall that
\begin{align*}
I_5 = \mathbb{E} \left[ 
F \left( \frac{1}{4} \int_0^a \int_0^{2\pi} 
\Phi_{\delta, I_\rho} \big( G \circ L_{\theta, \rho} \big) 
\, d\theta \, d\rho \right) 
\right] 
- 
\mathbb{E} \left[ 
F \left( \mathcal{L} \big( \mathcal{B}(G,a) \big) \right) 
\right].
\end{align*}

\begin{proposition}\label{prop-I5-z}
For all $a > 0$ and $r \in (0,1)$, there exists a constant $C_{r,a} > 0$ such that
\begin{align*}
I_5 \leq C_{r,a} \, \delta^{r}.
\end{align*}
\end{proposition}

Using Lemma \ref{crofton-lemma} for the paths of $G$, together with the fact that $F$ is $1$-Lipschitz, we obtain
\begin{align*}
I_5 
\leq \frac{1}{4} \int_0^a \int_0^{2\pi} 
\mathbb{E}\left[ 
\left| 
\Phi_{\delta, I_\rho} \big( G \circ L_{\theta, \rho} \big) 
- N\big( G \circ L_{\theta, \rho}, I_\rho \big) 
\right| 
\right] 
\, d\theta \, d\rho.
\end{align*}

We now adapt the arguments of Section 4.3 in \cite{lclp}.  
For clarity, we first introduce the necessary notation.

For every $f \in C^2(\mathbb{R}; \mathbb{R})$ and every interval $I = [\alpha,\zeta]$, define
\begin{equation*}
\mathcal{A}_{f,I}
= \min \left\{
|f(\alpha)|,\; |f(\zeta)|,\;
\frac{1}{2} \inf_{t \in I} \big( |f(t)| + |f'(t)| \big)
\right\}.
\end{equation*}

According to \cite[Lemma A.1]{lclp}, we have
\begin{align}\label{est-int-n-phi-delta}
\left| \Phi_{\delta, I}(f) - N(f,I) \right|
\leq 
\left( \left| \Phi_{\delta, I}(f) \right| + \left| N(f,I) \right| \right)
\mathbf{1}_{\mathcal{A}_{f,I} \leq \delta}.
\end{align}

Using \eqref{est-int-n-phi-delta}, we obtain
\begin{align*}
I_5 
\leq \frac{1}{4} \int_0^a \int_0^{2\pi} 
\mathbb{E}\Big[
\Big(
\left| \Phi_{\delta, I_\rho} \big( G \circ L_{\theta, \rho} \big) \right|
+ N\big( G \circ L_{\theta, \rho}, I_\rho \big)
\Big)
\mathbf{1}_{\mathcal{A}_{G \circ L_{\theta,\rho}, I_\rho} \leq \delta}
\Big] 
\, d\theta \, d\rho.
\end{align*}

Thus, Proposition \ref{prop-I5-z} is a consequence of the following result.

\begin{proposition}\label{prop-I5-int}
Let $r \in (0,1)$ and $a > 0$. Then there exists a constant
$C_{r,a}>0$ such that
\begin{align*}
\sup_{\theta \in [0,2\pi], \, \rho \in [-a,a]} 
\mathbb{E}\Big[
\Big(
\left| \Phi_{\delta, I_\rho} \big( G \circ L_{\theta, \rho} \big) \right|
+ N\big( G \circ L_{\theta, \rho}, I_\rho \big)
\Big)
\mathbf{1}_{\mathcal{A}_{G \circ L_{\theta,\rho}, I_\rho} \leq \delta}
\Big]
\leq C_{r,a} \, \delta^{r}.
\end{align*}
\end{proposition}

The proof of Proposition \ref{prop-I5-int} relies on the following lemma, which is a reformulation of Lemmas 4.3 and 4.9 in \cite{lclp}.

\begin{lemma}\label{lem-I5-X}
Let $X$ be a process with sample paths in $C^2(\mathbb{R}, \mathbb{R})$ such that for all $0 < r < 1$ and $p \geq 1$,
\begin{align}\label{hipo}
\mathbb{E} \left[
|X(a)|^{-r} + |X(b)|^{-r} + \|X\|_{\infty,2}^p
\right]
+ \sup_{s \in I} \mathbb{E} \left[ |X'(s)|^{-r} \right] < \infty,
\end{align}
where
\begin{equation*}
\|X\|_{\infty,2}
= \sup_{t \in I} \left( |X(t)| + |X'(t)| + |X''(t)| \right).
\end{equation*}

Then there exists a constant $C > 0$ such that for every compact interval $I$ and every $r,p$ satisfying $rp < 1$ and $1/p + 1/q = 1$,
\begin{align*}
\mathbb{P}\big( \mathcal{A}_{X,I} \leq \delta \big)
\leq 
C \, \delta^{r} \max(1,|I|)\, \varepsilon_{I,X,r,p},
\end{align*}
where
\begin{align}\label{def-erreur-delta}
\varepsilon_{I,X,r,p}
= \mathbb{E}\left[ |X(a)|^{-r} + |X(b)|^{-r} \right]
+ \mathbb{E}\left[ \|X\|_{\infty,2}^p \right]^{1/p}
\cdot
\sup_{s \in I} \mathbb{E} \left[ |X'(s)|^{-rq} \right]^{1/q}.
\end{align}

\end{lemma}

\begin{proof}[Proof of Proposition \ref{prop-I5-int}]
First we estimate $\varepsilon_{I_{\rho},G\circ L_{\theta,\rho},r,p}$. Let $p \geq 1$ and $r \in (0,1)$. There exists a constant $C_{p,r} > 0$ such that
\begin{align}\label{est-uni-epsi}
\varepsilon_{I_{\rho},G\circ L_{\theta,\rho},r,p}\leq C_{p,r}(1+ a^{1+p}).
\end{align}

According to Lemma \ref{lem-cov-g-derg}, for $\beta \in I_{\rho}$,
\begin{align*}
\mathbb{E}\left[\left| (G\circ L_{\theta,\rho})(\beta)\right|^{-r}\right]
= C_r.
\end{align*}

Moreover, for all $\beta \in I_{\rho}$,
\begin{align*}
\mathbb{E}\left[\left| (G\circ L_{\theta,\rho})'(\beta)\right|^{-r}\right]
&=
C_r
\left(
\frac{1}{3} + \frac{\sin(2\theta)}{4}
\right)^{-r/2} \\
&\leq C_r\,12^{r/2}.
\end{align*}

Recall that for $t=(t_1,t_2)\in {\mathbb R}^2$
\begin{align*}
G(t)= \int\int_{[0,1]^2} \cos( t_1u_1 + t_2u_2) \, dB^1(u_1,u_2) 
- \int\int_{[0,1]^2} \sin( t_1u_1 + t_2u_2) \, dB^2(u_1,u_2).
\end{align*}

Thus, the sample paths of $G$ are infinitely differentiable and, for all $j \in {\mathbb N}$ and all $\varepsilon_l \in \{1,2\}$, $l=1,\dots,j$,
\begin{align*}
\partial_{t_{\varepsilon_1},\dots,t_{\varepsilon_j}}^j G(t)
&= \int\int_{[0,1]^2} \prod_{l=1}^j u_{\varepsilon_l}
\cos\left( t_1u_1 + t_2u_2 + \frac{j\pi}{2} \right) \, dB^1(u_1,u_2) \\
&\quad - \int\int_{[0,1]^2} \prod_{l=1}^j u_{\varepsilon_l}
\sin\left( t_1u_1 + t_2u_2 + \frac{j\pi}{2} \right) \, dB^2(u_1,u_2).
\end{align*}
Then, for $X = G$ or $\partial_{t_{\varepsilon_1},\dots,t_{\varepsilon_j}}^j G$,
\begin{align*}
\mathbb{E}\left[ \left( X(t) - X(s) \right)^2 \right]
\leq 2 \left( |t_1 - s_1|^2 + |t_2 - s_2|^2 \right).
\end{align*}

Recall that 
\[
L_{\theta,\rho}(\beta)= (\rho \sin \theta + \beta \cos \theta,\, -\rho \cos \theta + \beta \sin \theta).
\]

Using the Garsia--Rodemich--Rumsey lemma for the Gaussian process $X = (G \circ L_{\theta,\rho})^{(j)}$, we derive
\begin{align*}
&\mathbb{E}\left[ \sup_{\beta \in I_{\rho}} \left| (G \circ L_{\theta,\rho})^{(j)}(\beta) - (G \circ L_{\theta,\rho})^{(j)}(0) \right|^p \right]\\
&\leq |I_{\rho}|^{p/2} \int_{I_{\rho} \times I_{\rho}} 
\frac{\mathbb{E}\left[ \left| (G \circ L_{\theta,\rho})^{(j)}(\beta) - (G \circ L_{\theta,\rho})^{(j)}(\beta') \right|^p \right]}{|\beta - \beta'|^{1 + p/2}} \, d\beta \, d\beta' \\
&\leq C_{p} |I_{\rho}|^{p+1}.
\end{align*}

We now derive the estimate for $\varepsilon_{I,X,r,p}$. Using the previous bounds on the negative moments of $X$ and $X'$, we obtain
\begin{align*}
\varepsilon_{I,X,r,p}
&\leq C_r + C_r \, \mathbb{E}[\|X\|_{\infty,2}^p]^{1/p} \\
&\leq C_{p,r} \left(1 + (1 + |I|^{p+1})^{1/p} \right).
\end{align*}
Using that for $p \geq 1$,
\[
(1 + |I|^{p+1})^{1/p} \leq C (1 + |I|^{p+1}),
\]
we obtain
\begin{align*}
\varepsilon_{I,X,r,p}
\leq C_{p,r} (1 + |I|^{p+1}).
\end{align*}
Since $|I_\rho| \lesssim a$, this yields
\begin{align*}
\varepsilon_{I_{\rho},G\circ L_{\theta,\rho},r,p}
\leq C_{p,r}(1 + a^{1+p}).
\end{align*}

Moreover, for all $p \geq 1$ and $j \in {\mathbb N}$,
\begin{align}\label{est-normsup-derg-p}
\sup_{\theta \in [0,2\pi],\, \rho \in [-a,a]} 
\mathbb{E}\left[\sup_{\beta \in I_{\rho}} \left| (G \circ L_{\theta,\rho})^{(j)}(\beta)\right|^p \right] < \infty.
\end{align}

Now, we turn to the estimation of
\[
\mathbb{E}\left[\left( \left| \Phi_{\delta, I_{\rho}}(G\circ L_{\theta,\rho}) \right| 
+ \left| N(G\circ L_{\theta,\rho}, I_\rho) \right| \right) 
\mathbf{1}_{\mathcal{A}_{G\circ L_{\theta,\rho}, I_\rho}\leq \delta}\right].
\]

We follow the same line of proof as in Lemma 4.10 of \cite{lclp}.  
Let $r \in (0,1)$. Let $(r',p,q)$ be such that $0 < r < r' < 1$, $r' = q r$, and $\frac{1}{p} + \frac{1}{q} = 1$.
Then, Hölder's and Jensen's inequalities imply that
\begin{align*}
\mathbb{E}\Big[
\big( \left| \Phi_{\delta, I_{\rho}}(G\circ L_{\theta,\rho}) \right| 
+ \left| N(G\circ L_{\theta,\rho}, I_\rho) \right| \big)
\mathbf{1}_{\mathcal{A}_{G\circ L_{\theta,\rho}, I_\rho}\leq \delta}
\Big] 
&\leq 2^{1-1/p} \, \mathbb{P}\big( \mathcal{A}_{G\circ L_{\theta,\rho}, I_\rho}\leq \delta \big)^{1/q}
\\
&\quad \times \left[ 
\mathbb{E} \left[ \left| \Phi_{\delta, I_{\rho}}(G\circ L_{\theta,\rho}) \right|^p \right] 
+ \mathbb{E} \left[ \left| N(G\circ L_{\theta,\rho}, I_\rho) \right|^p \right] 
\right]^{1/p}.
\end{align*}

Applying Lemma \ref{lem-I5-X} to $G\circ L_{\theta,\rho}$ and using \eqref{est-uni-epsi}, we obtain
\begin{align}\label{est-aronde-delta-G}
\mathbb{P}\big( \mathcal{A}_{G\circ L_{\theta,\rho}, I_\rho}\leq \delta \big)
\leq C_{r,a}\delta^r.
\end{align}
Using estimate \eqref{est-aronde-delta-G} with $r' = rq$, we obtain 
\begin{align}\label{est-I5-4} 
\mathbb{E}\Big[
\big( \left| \Phi_{\delta, I_{\rho}}(G\circ L_{\theta,\rho}) \right| 
+ \left| N(G\circ L_{\theta,\rho}, I_\rho) \right| \big)
\mathbf{1}_{\mathcal{A}_{G\circ L_{\theta,\rho}, I_\rho}\leq \delta}
\Big] 
&\leq C_{r,a} \delta^r
\\
&\quad \times \left[ 
\mathbb{E} \left[ \left| \Phi_{\delta, I_{\rho}}(G\circ L_{\theta,\rho}) \right|^p\right] 
+ \mathbb{E} \left[ \left| N(G\circ L_{\theta,\rho}, I_\rho )\right|^p\right] 
\right]^{1/p}.
\end{align} 

Recall that $G \circ L_{\theta,\rho}$ is a centered Gaussian process such that 
\[
\mathbb{E}\big[(G \circ L_{\theta,\rho})(\beta)^2\big]=1, \quad \text{for all } \beta \in \mathbb{R}.
\]
Then
\begin{align*}
\sup_{t \in B(0,a)} \sup_{x \in \mathbb{R}} p_{G(t)}(x) 
\leq \frac{1}{\sqrt{2\pi}}.
\end{align*}

According to Theorem 3.6 of \cite{Wse-Aza}, together with estimate \eqref{est-normsup-derg-p}, for all $p \geq 1$,
\begin{align}\label{maj-moment-zero-G}
\sup_{\theta \in [0,2\pi],\, \rho \in [-a,a]} 
\mathbb{E} \left[ N( G\circ L_{\theta,\rho}, [-a,a]) \right] < \infty.
\end{align}

In the same spirit, since 
\[
\frac{1}{12} \leq \mathbb{E}\big[((G \circ L_{\theta,\rho})'(\beta))^2\big] \leq \frac{7}{12}, 
\quad \text{for all } \beta \in \mathbb{R},
\]
we obtain
\begin{align*}
\sup_{\beta \in [-a,a]} \sup_{x \in \mathbb{R}} 
p_{(G\circ L_{\theta,\rho})'(\beta)}(x) 
\leq \frac{\sqrt{12}}{\sqrt{2\pi}},
\end{align*}
and
\begin{align*}
\sup_{\theta \in [0,2\pi],\, \rho \in [-a,a]} 
\mathbb{E} \left[ N\big( (G\circ L_{\theta,\rho})', [-a,a] \big) \right] < \infty.
\end{align*}

Now, using that 
\[
\Phi_{\delta, I_{\rho}}(G\circ L_{\theta,\rho}) 
\leq N(G\circ L_{\theta,\rho}, I_{\rho}) + 1
\]
(see the proof of Proposition 13 in \cite{ADLNP18}) and Corollary 3.7 in \cite{Wse-Aza}, we obtain
\begin{align}\label{maj-moment-phi-delta}
\sup_{\theta \in [0,2\pi],\, \rho \in [-a,a]} 
\mathbb{E} \left[ \left| \Phi_{\delta, I_{\rho}}(G\circ L_{\theta,\rho}) \right|^p \right] < \infty.
\end{align}

Plugging estimates \eqref{maj-moment-zero-G} and \eqref{maj-moment-phi-delta} into \eqref{est-I5-4} completes the proof of Proposition \ref{prop-I5-int}.
\end{proof}


\subsubsection{Estimation of  $I_1$}
We recall that
\begin{align*}
I_1 = \mathbb{E} \left[ F \left( \frac{1}{4} \int_0^a \int_0^{2\pi} 
\Phi_{\delta, I_\rho} (P^m \circ L_{\theta, \rho} ) \, d\theta \, d\rho \right) \right] 
- \mathbb{E} \left[ F \left( \mathcal{L}(\mathcal{B}(P^m,a)) \right) \right].
\end{align*}

\begin{proposition}\label{prop-I1-z}
For all $a > 0$, $r\in(0,1)$, $q \in (1,\frac{1}{r})$ and $0 < \gamma < \frac{1-\frac{4}{p}}{5-\frac{4}{p}}$, there exists a constant $C_{r,a,\gamma} > 0$ such that
\begin{align*}
I_1 \leq C_{r,a,\gamma} \left[\delta^r
+
\frac{1+\varepsilon}{\delta\varepsilon}
m^{-\gamma}\ln m
+
\varepsilon^{-1/q}
\big(m^{-\gamma}\ln m\big)^{1/q}\right].
\end{align*}
\end{proposition}
{Let ${\mathcal D}_0={\mathcal D}_{0,P_m,a}$ (see \eqref{def-vad-crofton}).} On ${\mathcal D}_0^c$, we use Lemma \ref{crofton-lemma} for the paths of $P^m$, together with the fact that $F$ is $1$-Lipschitz. On ${\mathcal D}_0$, we use that $\|F\|_{\infty}\leq 1$ and Lemma \ref{Pm-crofton}. We obtain
{\begin{align*}
I_1 \leq \frac{1}{4} \int_0^a \int_0^{2\pi}  
\mathbb{E}\left[ 
\left| \Phi_{\delta, I_\rho} (P^m \circ L_{\theta, \rho}) 
- N(P^m \circ L_{\theta, \rho}, I_{\rho}) \right|  
\right] 
\, d\theta \, d\rho
+ 2\,{\mathbb P}[{\mathcal D}_0]
\end{align*}}
and
\begin{align}\label{aux:boun_D_0}
{\mathbb P}[{\mathcal D}_0]\leq C\, m^{-\frac{\gamma r}{1+r}}.
\end{align}
According to estimate \eqref{est-int-n-phi-delta},
\begin{align}\label{est-splii1}
I_1 \leq \frac{1}{4} \int_0^a \int_0^{2\pi}  
\mathbb{E}\Big[
\big( \left| \Phi_{\delta, I_\rho} (P^m \circ L_{\theta, \rho}) \right| 
+ N(P^m \circ L_{\theta, \rho}, I_{\rho}) \big)
\mathbf{1}_{\mathcal{A}_{P^m \circ L_{\theta,\rho}, I_{\rho}} \leq \delta}
\Big] 
\, d\theta \, d\rho
+ 2\,{\mathbb P}[{\mathcal D}_0].
\end{align}

We first prove the following lemma, adapted from Lemma 4.9 of \cite{lclp}.

\begin{lemma}\label{lem-maj-padeltapm}
For all $r \in (0,1)$ and $0 < \gamma < \frac{1-\frac{4}{p}}{5-\frac{4}{p}}$, there exists a constant $C_{r,a,\gamma} > 0$ such that for all $m \in \mathbb{N}$, all $\varepsilon < \delta$, and all sequences $(X_{k,l}, Y_{k,l})_{k,l \in \mathbb{N}}$ satisfying Hypothesis \ref{hypxp-1}
\begin{align*}
\sup_{\theta \in [0,2\pi],\, \rho \in [-a,a]} 
\mathbb{P}\left( \mathcal{A}_{P^m \circ L_{\theta,\rho}, I_{\rho}} \leq \delta \right) 
\leq C_{r,a,\gamma} \left[ \delta^{r} + \frac{1}{\varepsilon} m^{-\gamma} \ln m \right] 
\mathbb{E}[ |X_{1,1}|^p + |Y_{1,1}|^p ].
\end{align*}
\end{lemma}

\begin{proof}
Let us denote, for continuous functions on $[-a,a]$,
\begin{equation*}
\overline{\Psi}_{\delta,\varepsilon,I_{\rho}} (f)= H_{\delta,\varepsilon} ( {\mathcal A}_{f,I_{\rho}}),
\end{equation*}
where $H_{\delta,\varepsilon}$ is defined in \eqref{def-H-delta-epsi}. 

Since ${\mathbf 1}_{u\leq \delta} \leq H_{\delta,\varepsilon}(u) \leq {\mathbf 1}_{u\leq 2\delta}$, we have
\begin{align*}
\mathbb{P}\left( \mathcal{A}_{P^m\circ L_{\theta,\rho},I_{\rho}}\leq \delta\right) 
&\leq \mathbb{E} \left[H_{\delta,\varepsilon} ( \mathcal{A}_{P^m\circ L_{\theta,\rho},I_{\rho}}) \right] \\
&\leq \left| \mathbb{E} \left[H_{\delta,\varepsilon} ( \mathcal{A}_{P^m\circ L_{\theta,\rho},I_{\rho}}) \right] 
- \mathbb{E} \left[H_{\delta,\varepsilon} ( \mathcal{A}_{G\circ L_{\theta,\rho},I_{\rho}}) \right] \right| \\
&\quad + \mathbb{E} \left[H_{\delta,\varepsilon} ( \mathcal{A}_{G\circ L_{\theta,\rho},I_{\rho}}) \right] \\
&\leq \left| \mathbb{E}\left[\overline{\Psi}_{\delta,\varepsilon,I_{\rho}}(P^m) \right] 
- \mathbb{E}\left[\overline{\Psi}_{\delta,\varepsilon,I_{\rho}}(G) \right] \right| \\
&\quad + \mathbb{P}\left( \mathcal{A}_{G\circ L_{\theta,\rho},I_{\rho}}\leq 2\delta\right).
\end{align*}

The map $\overline{\Psi}_{\delta,\varepsilon,I_{\rho}}$ is Lipschitz from $C^1([-a,a],\mathbb{R})$ into $\mathbb{R}$ with Lipschitz constant bounded by $3 \varepsilon^{-1}$. 

Thus, according to Corollary \ref{cor:stein-bounded}, for $0< \gamma < \frac{1-\frac{4}{p}}{5-\frac{4}{p}}$, there exists a constant $C_{p,\gamma}$ such that
\begin{align*}
\sup_{\theta \in [0,2\pi],\, \rho \in [-a,a]}  
\left| \mathbb{E} \left[\overline{\Psi}_{\delta,\varepsilon,I_{\rho}}(P^m) \right] 
- \mathbb{E} \left[\overline{\Psi}_{\delta,\varepsilon,I_{\rho}}(G) \right] \right| 
\leq C_{p,\gamma} \frac{1}{\varepsilon} m^{-\gamma} \ln m \, \mathbb{E}[ |X_{1,1}|^p + |Y_{1,1}|^p].
\end{align*}

Applying estimate \eqref{est-aronde-delta-G} with $2\delta$, we obtain
\begin{align*}
\sup_{\theta \in [0,2\pi],\, \rho \in [-a,a]} 
\mathbb{P}\left( \mathcal{A}_{P^m\circ L_{\theta,\rho},I_{\rho}}\leq \delta\right) 
\leq C_{p,\gamma} \frac{1}{\varepsilon} m^{-\gamma} \ln m \, \mathbb{E}[ |X_{1,1}|^p + |Y_{1,1}|^p]
+ 2^r C_{r,a} \delta^r.
\end{align*}
\end{proof}

Using the same ideas, we obtain the following lemma.

\begin{lemma}\label{lem-I1-phidelta}
For all $r \in (0,1)$ and $0< \gamma < \frac{1-\frac{4}{p}}{5-\frac{4}{p}}$, there exists a constant $C_{r,a,\gamma} > 0$ such that for all $m \in \mathbb{N}$, all $\varepsilon < \delta$, and all sequences $(X_{k,l}, Y_{k,l})_{k,l \in \mathbb{N}}$ satisfying Hypothesis \ref{hypxp-1}

\begin{align*}
\sup_{\theta \in [0,2\pi],\, \rho \in [-a,a]} 
\mathbb{E}\left[
\left| \Phi_{\delta, I_\rho} (P^m \circ L_{\theta, \rho} )\right|  
\mathbf{1}_{\mathcal{A}_{P^m \circ L_{\theta,\rho},I_{\rho}}\leq \delta}
\right]
\leq C_{r,a,\gamma} \left[\delta^{r} +  \frac{1+\varepsilon}{\delta\varepsilon} m^{-\gamma} \ln m \right]
\mathbb{E}[ |X_{1,1}|^p + |Y_{1,1}|^p].
\end{align*}
\end{lemma}

\begin{proof}
Recall that $ \Phi_{\delta, I} \leq \Phi_{\delta, I,\varepsilon} \leq \Phi_{2\delta, I}$. Then
\begin{align*}
\mathbb{E}\left[
\left| \Phi_{\delta, I_\rho} (P^m \circ L_{\theta, \rho} )\right|  
\mathbf{1}_{\mathcal{A}_{P^m \circ L_{\theta,\rho},I_{\rho}}\leq \delta}
\right]
&\leq 
\left| 
\mathbb{E}\left[
\left| \Phi_{\delta, I_\rho,\varepsilon } (P^m \circ L_{\theta, \rho} )\right| 
H_{\delta,\varepsilon} ( \mathcal{A}_{P^m \circ L_{\theta,\rho},I_{\rho}})   
\right]\right. \\
&\quad \left.
- \mathbb{E}\left[
\left| \Phi_{\delta, I_\rho,\varepsilon } (G \circ L_{\theta, \rho} )\right| 
H_{\delta,\varepsilon} ( \mathcal{A}_{G \circ L_{\theta,\rho},I_{\rho}})   
\right]
\right| \\
&\quad + \mathbb{E}\left[
\left| \Phi_{2\delta, I_\rho} (G \circ L_{\theta, \rho} )\right|  
\mathbf{1}_{\mathcal{A}_{G \circ L_{\theta,\rho},I_{\rho}}\leq 2\delta}
\right].
\end{align*}

The map 
\[
\frac{\delta \varepsilon}{|I_{\rho}|(1+\varepsilon)} 
\Phi_{\delta, I_{\rho} ,\varepsilon}\,\overline{\Psi}_{\delta,\varepsilon,I_{\rho}}
\]
belongs to ${\mathcal L}_{1,1}(C^1([-a,a],\mathbb{R}),\mathbb{R})$. 

Thus, according to Corollary \ref{cor:stein-bounded}, for $0< \gamma < \frac{1-\frac{4}{p}}{5-\frac{4}{p}}$, there exists a constant $C_{p,\gamma,a}$ such that 
\begin{align*}
\sup_{\theta \in [0,2\pi],\, \rho \in [-a,a]}
\left| 
\mathbb{E}\left[
\left| \Phi_{\delta, I_\rho,\varepsilon } (P^m \circ L_{\theta, \rho} )\right| 
H_{\delta,\varepsilon} ( \mathcal{A}_{P^m \circ L_{\theta,\rho},I_{\rho}})   
\right]\right. \\
\left.
- \mathbb{E}\left[
\left| \Phi_{\delta, I_\rho,\varepsilon } (G \circ L_{\theta, \rho} )\right| 
H_{\delta,\varepsilon} ( \mathcal{A}_{G \circ L_{\theta,\rho},I_{\rho}})   
\right]
\right|
\leq C_{p,\gamma,a}\frac{1}{\delta \varepsilon} m^{-\gamma} \ln m \, 
\mathbb{E}[ |X_{1,1}|^p + |Y_{1,1}|^p].
\end{align*}

Then, using Proposition \ref{prop-I5-int} for the term in $G$, we obtain
\begin{align*}
\mathbb{E}\left[
\left| \Phi_{\delta, I_\rho} (P^m \circ L_{\theta, \rho} )\right|  
\mathbf{1}_{\mathcal{A}_{P^m \circ L_{\theta,\rho},I_{\rho}}\leq \delta}
\right]
&\leq C_{p,\gamma,a}\frac{1}{\delta \varepsilon} m^{-\gamma} \ln m \, 
\mathbb{E}[ |X_{1,1}|^p + |Y_{1,1}|^p]
+ C_{a,r} \delta^r.
\end{align*}
\end{proof}
\begin{proof}[Proof of Proposition \ref{prop-I1-z}]
Recall that, for all $p^{\prime}\geq 1$, by Proposition~\ref{prop-maj-expec-N} applied with
\(I=[-a,a]\),
\begin{equation*}
\sup_m
\sup_{\theta\in[0,2\pi],\,\rho\in[-a,a]}
\mathbb{E}\!\left[
N\!\left(P^m\circ L_{\theta,\rho},[-a,a]\right)^{p^{\prime}}
\right]
=
\sup_m
\sup_{\theta\in[0,2\pi],\,\rho\in[-a,a]}
\mathbb{E}\!\left[
N\!\left(P_{\theta,\rho}^{m},[-a,a]\right)^{p^{\prime}}
\right]
<\infty.
\end{equation*}

Thus, applying Hölder's inequality together with Lemma \ref{lem-maj-padeltapm} for $r' = qr < 1$ and $p^{\prime} = \frac{q}{q-1}$, we obtain
\begin{align*}
\sup_{\theta \in [0,2\pi],\, \rho \in [-a,a]} 
\mathbb{E}\left[
N(P^m \circ L_{\theta, \rho}, I_{\rho}) 
\mathbf{1}_{\mathcal{A}_{P^m \circ L_{\theta,\rho}, I_{\rho}} \leq \delta}
\right] 
\leq C_{r,a,\gamma}\left(\delta^r+\varepsilon^{-1/q}\left(m^{-\gamma}\ln(m)\right)^{1/q}\right).
\end{align*}

{Combining \eqref{est-splii1}, \eqref{aux:boun_D_0}, Lemma \ref{lem-I1-phidelta} we obtain}

\begin{align*}
m^{-\frac{\gamma r}{1+r}}\leq \delta^r+\frac{1+\varepsilon}{\delta\varepsilon}\,m^{-\gamma}\ln m
\end{align*}
and therefore
\begin{align*}
I_1\leq C_{r,a,\gamma}
\Bigg[
\delta^r
+
\frac{1+\varepsilon}{\delta\varepsilon}
m^{-\gamma}\ln m
+
\varepsilon^{-1/q}
\big(m^{-\gamma}\ln m\big)^{1/q}
\Bigg].
\end{align*}
Then, we conclude the proof.
\end{proof}
\subsubsection{Estimation of  $I_2$}
Recall that
\begin{align*}
I_2 := \mathbb{E} \left[ F \left( \frac{1}{4} \int_0^a \int_0^{2\pi} \Phi_{\delta, I_{\rho}} (P^m \circ L_{\theta, \rho} ) \, d\theta \, d\rho \right) \right] 
- \mathbb{E} \left[ F \left( \frac{1}{4} \int_0^a \int_0^{2\pi} \Phi_{\delta, \varepsilon, I_\rho} (P^m \circ L_{\theta,\rho} ) \, d\theta \, d\rho \right) \right].
\end{align*}
\begin{proposition}\label{I_2}
Let $p>4$, $a>0$ and $\gamma < \frac{1-\frac{4}{p}}{5-\frac{4}{p}}$. There exists a constant $C_{p,\gamma,a} > 0$ such that for all $1$-Lipschitz functions $F$ and all sequences $\{X_{k,l},Y_{k,l},\,k,l\in \mathbb{N}\}$ of centered independent identically distributed random variables satisfying Hypothesis \ref{hypxp-1},
\begin{align*}
I_2 \leq C_{p,\gamma,a}
\left(\frac{1}{\varepsilon \delta} m^{-\gamma} \ln (m) + \frac{\varepsilon}{\delta}\right).
\end{align*}
\end{proposition}

\begin{proof}
First, we note that for any function $f \in C^1(I,\mathbb{R})$,
\begin{align*}
0 \leq \Phi_{\delta,\varepsilon,I}(f) - \Phi_{\delta,I}(f)
&\leq \frac{1}{2\delta} \int_I |f'(s)| \frac{|\delta + \varepsilon - |f(s)||}{\varepsilon} \mathbf{1}_{\delta \leq |f(s)| \leq \delta + \varepsilon} \, ds \\
&\leq \tilde{\psi}_{\delta,\varepsilon,I}(f),
\end{align*}
where $\tilde{\psi}_{\delta,\varepsilon,I}$ is defined on $C^{1}( I,\mathbb{R})$ by
\begin{align*}
\tilde{\psi}_{\delta,\varepsilon,I}(f) := \frac{1}{2\delta} \int_I |f'(s)| \tilde{H}_{\delta,\varepsilon}(|f(s)|)\, ds,
\end{align*}
and $\varepsilon \tilde{H}_{\delta,\varepsilon}$ is a $1$-Lipschitz function bounded by $\varepsilon$ on $\mathbb{R}$ given by
\begin{align*}
\tilde{H}_{\delta,\varepsilon}(x)= 
\begin{cases}
0 & \text{if } x > \delta + \varepsilon, \\
\frac{\delta + \varepsilon - x}{\varepsilon} & \text{if } \delta < x < \delta + \varepsilon, \\
\frac{\delta - x}{\varepsilon} & \text{if } \delta - \varepsilon < x < \delta, \\
0 & \text{if } -\delta < x < \delta - \varepsilon, \\
\frac{\delta + \varepsilon + x}{\varepsilon} & \text{if } -\delta - \varepsilon < x < -\delta, \\
0 & \text{if } x < -\delta - \varepsilon.
\end{cases}
\end{align*}

Thus, the map
\[
\frac{\delta \varepsilon}{(1+\varepsilon)a}\, \tilde{\psi}_{\delta,\varepsilon,I_{\rho}}
\]
belongs to ${\mathcal L}ip_{1,1}(C^1([-a,a],\mathbb{R}),\mathbb{R})$. 

Consequently, using the above notation and Corollary \ref{cor-P-m-cv-G}, for $0<\gamma<\frac{1-\frac4p}{5-\frac4p}$, we obtain
\begin{align*}
&\mathbb{E}\left[
\left| \Phi_{\delta, I_{\rho}}(P^{m}\circ L_{\theta,\rho}) 
- \Phi_{\delta, \varepsilon, I_{\rho}}(P^{m}\circ L_{\theta,\rho}) \right|
\right] \\
&\leq 
\left| \mathbb{E}\left[\tilde{\psi}_{\delta,\varepsilon,I_{\rho}}(P^{m}\circ L_{\theta,\rho})\right]
- \mathbb{E}\left[\tilde{\psi}_{\delta,\varepsilon,I_{\rho}}(G\circ L_{\theta,\rho})\right] \right| \\
&\quad + \mathbb{E}\left[\tilde{\psi}_{\delta,\varepsilon,I_{\rho}}(G\circ L_{\theta,\rho})\right] \\
&\leq \frac{C_{p,\gamma}}{\varepsilon \delta}\frac{\ln m}{m^{\gamma}}
+ \mathbb{E}\left[\tilde{\psi}_{\delta,\varepsilon,I_{\rho}}(G\circ L_{\theta,\rho})\right].
\end{align*}

For all $\beta \in [-a,a]$, the vector 
\[
\big(G\circ L_{\theta,\rho}(\beta),\, (G\circ L_{\theta,\rho})'(\beta)\big)
\]
is centered Gaussian with covariance given in Lemma \ref{lem-cov-g-derg}. Its components are independent, with variances $1$ and $\frac{1}{3} + \frac{\sin(2\theta)}{4}$, respectively. 

Thus,
\begin{align*}
\mathbb{E}\left[\tilde{\psi}_{\delta,\varepsilon,I_{\rho}}(G\circ L_{\theta,\rho})\right] 
&\leq \frac{1}{2\delta\varepsilon} \int_{-a}^{a}
\int_{\mathbb{R}^2} |y| \big(\delta + \varepsilon - |x|\big)
\mathbf{1}_{\delta < |x| < \delta+\varepsilon}
\, p_{(G\circ L_{\theta,\rho}(\beta), (G\circ L_{\theta,\rho})'(\beta))}(x,y)\, dx\, dy\, d\beta \\
&\leq \frac{1}{2\delta} \int_{-a}^{a} \int_{\mathbb{R}} |y|\, p_{(G\circ L_{\theta,\rho})'(\beta)}(y)
\frac{2\varepsilon}{\sqrt{2\pi}}\, dy\, d\beta \\
&\leq \frac{a\varepsilon\sqrt{7}}{\delta \sqrt{24 \pi}}.
\end{align*}

Finally,
\begin{align*}
I_2 
&\leq \frac{1}{4}\int_{0}^a \int_{0}^{2\pi}
\mathbb{E}\left(
\left|\Phi_{\delta, I_{\rho}}\left(P^{m}\circ L_{\theta,\rho}\right)
- \Phi_{\delta, \varepsilon, I_{\rho}}\left(P^{m}\circ L_{\theta,\rho}\right)\right|
\right)
\, d\theta\, d\rho \\
&\leq \frac{1}{4}\int_{0}^a \int_{0}^{2\pi}
C_{p,\gamma}\left(\frac{1}{\delta\varepsilon}m^{-\gamma}\ln m + \frac{\varepsilon}{\delta}\right)
\, d\theta\, d\rho \\
&= C_{p,\gamma}
\left(\frac{1}{\delta\varepsilon}m^{-\gamma}\ln m + \frac{\varepsilon}{\delta}\right)
\frac{1}{4} \int_{0}^a \int_{0}^{2\pi} d\theta\, d\rho \\
&= C_{p,\gamma,a}
\left(\frac{1}{\delta\varepsilon}m^{-\gamma}\ln m + \frac{\varepsilon}{\delta}\right).
\end{align*}
\end{proof}


\subsection{Optimization of the Parameters and Final Convergence Rate}\label{main:result}

We are now in a position to combine the estimates obtained for the error terms $I_1,\ldots,I_5$ and derive the quantitative rate of convergence for the nodal length.

Recall that, for every Lipschitz continuous function $F:\mathbb R\to\mathbb R$,
\begin{align*}
\mathbb E\left[
F\left({\mathcal L}(P_m,B(0,a))\right)
\right]
-
\mathbb E\left[
F\left({\mathcal L}(G,B(0,a))\right)
\right]
=
\sum_{i=1}^{5} I_i.
\end{align*}

The estimates established in the previous subsections provide bounds for each of the terms $I_1,\ldots,I_5$ in terms of the smoothing parameters $\delta$ and $\varepsilon$, as well as the degree parameter $m$. Optimizing these parameters yields the final rate of convergence.

We recall the statement of the main result of the paper.

\begin{theorem*}[Theorem~\ref{thm:main-rate}]
Let $(P_m)_{m\geq 1}$ be the sequence of random trigonometric polynomials defined in \eqref{def:pol_1}. Assume that the family $
\left\{(a_{k,\ell},b_{k,\ell}) : k,\ell\in\mathbb N\right\}$
consists of independent and identically distributed centered random vectors in $\mathbb R^2$ with identity covariance matrix and
\begin{equation*}
\E\left[\|(a_{1,1},b_{1,1})\|^p\right]<\infty
\end{equation*}
for some $p>4$. Let $a>0$ and let
\begin{equation*}
0<\gamma<\frac{1-\frac{4}{p}}{5-\frac{4}{p}}.
\end{equation*}
Let $G$ be the Gaussian field defined in Section~\ref{RTP_GL}. Then, for every
\begin{equation*}
0<\alpha<\frac{\gamma}{4},
\end{equation*}
there exists a constant $C=C(a,p,\gamma,\alpha)>0$ such that
\begin{equation*}
\sup_{\|F\|_{\mathrm{Lip}}+\|F\|_{\infty}\leq 1}
\left|
\E\!\left[F\!\left(\mathcal L(P_m,B(0,a))\right)\right]
-
\E\!\left[F\!\left(\mathcal L(G,B(0,a))\right)\right]
\right|
\leq
C\,m^{-\alpha}.
\end{equation*}\end{theorem*}

\begin{proof}
{Let $F:\mathbb R\to\mathbb R$ be such that $\|F\|_{\mathrm{Lip}}+\|F\|_{\infty}\leq 1$, in particular, $F$ is $1$-Lipschitz and $\|F\|_{\infty}\leq 1$.} We now combine the estimates obtained for the terms
$I_1,\ldots,I_5$ in order to derive the final rate of convergence. By Propositions
\ref{prop-I1-z}, \ref{I_2}, \ref{I_3}, \ref{I_4} and \ref{prop-I5-z}, we obtain
\begin{align*}
&\left|
\mathbb E\left[
F\left({\mathcal L}(P_m,B(0,a))\right)
\right]
-
\mathbb E\left[
F\left({\mathcal L}(G,B(0,a))\right)
\right]
\right| \\
&\leq
C_{r,a,\gamma}
\left[
\delta^r
+
\frac{1+\varepsilon}{\delta\varepsilon}
m^{-\gamma}\ln m
+
\varepsilon^{-1/q}
\bigl(m^{-\gamma}\ln m\bigr)^{1/q}
\right] \\
&\quad
+
C_{p,\gamma,a}
\left(
\frac{1}{\varepsilon \delta}
m^{-\gamma}\ln m
+
\frac{\varepsilon}{\delta}
\right)\\
&\quad
+
C_{p,\gamma}a^2(1+a)
\frac{\pi(1+\varepsilon)}
{2C\,\delta\varepsilon}
m^{-\gamma}
+
\frac{7\sqrt{3}\,\varepsilon a^2}{24\,\delta}
+
C_{r,a}\,\delta^r\\
&\leq
C
\left[
\delta^r
+
\frac{\varepsilon}{\delta}
+
\frac{1}{\delta\varepsilon}
m^{-\gamma}\ln m
+
\varepsilon^{-1/q}
\bigl(m^{-\gamma}\ln m\bigr)^{1/q}
\right],
\end{align*}
where in the last inequality we used that $0<\varepsilon<1$ and that
$m^{-\gamma}\leq m^{-\gamma}\ln m$ for $m$ sufficiently large.
Let $
B=m^{-\gamma}\ln m$ and define
\begin{equation*}
f_B(\delta,\varepsilon)
=
\frac{B}{\delta\varepsilon}
+
\frac{\varepsilon}{\delta}
+
\delta^r
+
\left(\frac{B}{\varepsilon}\right)^{1/q}.
\end{equation*}

Balancing the first two terms yields
\begin{equation*}
\frac{B}{\delta\varepsilon}
\sim
\frac{\varepsilon}{\delta},
\end{equation*}
that is,
\begin{equation*}
\varepsilon\sim \sqrt B.
\end{equation*}

Next, balancing the second and third terms gives
\begin{equation*}
\frac{\varepsilon}{\delta}
\sim
\delta^r,
\end{equation*}
hence
\begin{equation*}
\varepsilon\sim \delta^{r+1},
\qquad
\delta\sim \varepsilon^{1/(r+1)}
\sim
B^{\frac{1}{2(r+1)}}.
\end{equation*}

Therefore
\begin{equation*}
\delta^r
\sim
B^{\frac{r}{2(r+1)}}.
\end{equation*}

Moreover,
\begin{equation*}
\left(\frac{B}{\varepsilon}\right)^{1/q}
\sim
\left(\frac{B}{\sqrt B}\right)^{1/q}
=
B^{\frac1{2q}}.
\end{equation*}

Since $1<q<\tfrac{1}{r}$,
\begin{equation*}
\frac1{2q}>\frac{r}{2}>\frac{r}{2(r+1)},
\end{equation*}
and consequently
\begin{equation*}
f_B(\delta,\varepsilon)
\leq
C_r
B^{\frac{r}{2(r+1)}}.
\end{equation*}
On the other hand, define
\begin{equation*}
\phi(r)=\frac{r\gamma}{2(r+1)},
\qquad r\in(0,1).
\end{equation*}

Since
\begin{equation*}
\phi'(r)
=
\frac{\gamma}{2(r+1)^2}>0,
\end{equation*}
the function $\phi$ is increasing on $(0,1)$ and
\begin{equation*}
\sup_{0<r<1}\phi(r)
=
\lim_{r\to 1^-}\phi(r)
=
\frac{\gamma}{4}.
\end{equation*}

Therefore, for every
\begin{equation*}
0<\alpha<\frac{\gamma}{4},
\end{equation*}
there exists $r\in(0,1)$ such that
\begin{equation*}
\alpha<
\frac{r\gamma}{2(r+1)}.
\end{equation*}

It follows that
\begin{equation*}
f_B(\delta,\varepsilon)
\leq
C_{\alpha}
m^{-\alpha}
\end{equation*}
for all sufficiently large $m$.
\end{proof}
\appendix
\section{Non-degeneracy property}\label{app-nondeg}
Crofton's formula (see Lemma \ref{crofton-lemma}) is valid for any $C^1$ random field $Y\colon B(0,a)\rightarrow {\mathbb R}$ such that ${\mathbb P}\left[{\mathcal D}_{0,Y,a}\right]=0$ where for all $\delta \geq 0 $
\begin{align}\label{def-vad-crofton}
{\mathcal D}_{\delta,Y,a}= \{\omega, \min_{t \in B(0,a)}\big(|Y(\omega,t)| + \|\nabla Y(\omega,t)\|\big) \leq \delta\}.
\end{align}

We have the following results. 
\begin{lemma}\label{g-crofton}
For all $a>0$, ${\mathbb P}\left[{\mathcal D}_{0,G,a}\right]=0.$ 
Moreover, for all $r \in (0,1)$, there exists a constant $C_r$ such that, for all $\delta >0$ and $a\geq 0$,\begin{align*}
{\mathbb P}\left[{\mathcal D}_{\delta ,G,a}\right]\leq C_r (1+a^2) \delta^{r}.
\end{align*}
\end{lemma}

\begin{lemma}\label{Pm-crofton}
Let $p>4$, $0<\gamma < \frac{1- \frac{4}{p}}{5- \frac{4}{p}}$, $a>0$ and $r\in(0,1)$. Let $\left\{(a_{k,\ell},b_{k,\ell}) : k,\ell\in\mathbb N\right\}$ be a family of independent and identically distributed centered random vectors in $\mathbb R^2$ with identity covariance matrix such that
\begin{equation*}
\E\left[\|(a_{1,1},b_{1,1})\|^p\right]<\infty.
\end{equation*}
Then there exists a constant $C$, depending only on $p$, $\gamma$, $r$, $a$ and $\E\left[\|(a_{1,1},b_{1,1})\|^p\right]$, such that, for all $\delta\in(0,1]$,
\begin{align*}
\sup_{m\geq \delta^{-\frac{1+r}{\gamma}}} {\mathbb P}\left[{\mathcal D}_{\delta ,P_m,a}\right]\leq C \delta^{r}.
\end{align*}
In particular, ${\mathbb P}\left[{\mathcal D}_{0,P_m,a}\right]\leq C\, m^{-\frac{\gamma r}{1+r}}$ for all $m\in\mathbb N$.
\end{lemma}

\begin{remark} The proof of Theorem 2 of \cite{AngstPhamPoly} relies on Theorem 6 of \cite{AngstPhamPoly}, which requires that the level set be a $C^1$ curve.
However, this property does not seem to be established there, so it appears to us that an argument is missing.
\end{remark}

\subsection{Proofs}

\begin{lemma}\label{lem-gauss-crofton}
Let $Z=(Z(t), t \in [-a,a]^2)$ be a centered Gaussian field with $C^3$ sample paths, and let $\Gamma_t$ denote the covariance matrix of $(Z(t),\partial_1 Z(t), \partial_2 Z(t))$. Assume that there exist two constants $0<\underline{C}<\overline{C}$ such that, for all $v \in {\mathbb R}^3$ and all $t \in [-a,a]^2$,
\begin{align*}
&\underline{C} \|v\|^2 \leq \langle v, \Gamma_t v \rangle \leq \overline{C}\|v\|^2,\\
&{\mathbb E} \left[ |\partial Z(t)|^2 \right] \leq \overline{C}, \quad \text{for } \partial\in\{\partial_{11},\partial_{12},\partial_{22},\partial_{112},\partial_{122}\}.
\end{align*}
Then, for all $r\in(0,1)$, there exists a constant $C_{r}$, depending only on $\underline{C}$, $\overline{C}$ and $r$, such that
\begin{align*}
{\mathbb E} \left[\sup_{t\in [-a,a]^2} \left( \frac{1}{Z(t)^2 + \partial_1 Z(t)^2 + \partial_2 Z(t)^2} \right)^{\frac{r}{2}}\right] \leq C_{r}(1+a+ a^2). 
\end{align*}
\end{lemma}
This lemma is the two-dimensional analogue of Lemma 4.3 in \cite{lclp}.

\begin{proof}
Let $f\colon {\mathbb R}^3 \rightarrow {\mathbb R}$ and $x\colon [-a,a]^2 \rightarrow {\mathbb R}^3$ be $C^2$ functions. Applying the fundamental theorem of calculus and the chain rule repeatedly, we obtain
\begin{align*}
f(x(t))&= f(x(-a,-a)) + \int_{-a}^{t_1} \langle \nabla f(x(u,-a)),\partial_1x(u,-a)\rangle du + \int_{-a}^{t_2} \langle \nabla f(x(-a,v)),\partial_2x(-a,v)\rangle dv\\
&+ \int_{-a}^{t_1} \int_{-a}^{t_2} \left\{D^2 f(x(u,v))\big(\partial_1x(u,v),\partial_2x(u,v)\big)+  \langle \nabla f(x(u,v)), \partial_{12}x(u,v)\rangle\right\} du dv.
\end{align*}
Thus, by the Cauchy--Schwarz inequality and the definition of the norm $\|D^2 f\|$,
\begin{align*}
\sup_{t \in [-a,a]^2}|f(x(t)&)|\leq  | f(x(-a,-a))| + \int_{-a}^{a} \|\nabla f(x(u,-a))\|\,\|\partial_1x(u,-a) \|du\\
& + \int_{-a}^{a} \|\nabla f(x(-a,v))\|\,\|\partial_2x(-a,v)\| dv\\
&+ \int \int_{[-a,a]^2}  \left\{\|D^2 f(x(u,v))\|\,\|\partial_1x(u,v)\|\,\|\partial_2x(u,v)\|+  \|\nabla f(x(u,v))\|\,\|\partial_{12}x(u,v)\|\right\} du dv.
\end{align*}
Using the fact that, by Hölder's and the Cauchy--Schwarz inequalities, for conjugate exponents $p,q$ and nonnegative random variables $U$, $V$ and $W$,
\begin{align*}
{\mathbb E} \left[UVW\right]\leq {\mathbb E} \left[U^p\right]^{1/p} {\mathbb E} \left[V^{2q}\right]^{1/(2q)} {\mathbb E} \left[W^{2q}\right]^{1/(2q)},
\end{align*}
we derive that, for every $C^2$ random field $(X(t), t \in [-a,a]^2)$,
\begin{align}\label{lem-gauss-crofton-1} 
{\mathbb E} \left[ \sup_{t \in [-a,a]^2}|\right. & \left.f(X(t))|\right] \leq {\mathbb E} \left[ | f(X(-a,-a))|\right] \\
&+ 2a \sup_{t \in [-a,a]^2}{\mathbb E} \left[ \|\nabla f(X(t))\|^p\right]^{1/p} \left\{\sup_{t \in [-a,a]^2} {\mathbb E} \left[ \|\partial_1X(t)\|^q\right]^{1/q}+ \sup_{t \in [-a,a]^2} {\mathbb E} \left[ \|\partial_2X(t)\|^q\right]^{1/q}\right\}\nonumber\\
&+4a^2\sup_{t \in [-a,a]^2}{\mathbb E} \left[ \|\nabla f(X(t))\|^p\right]^{1/p} \sup_{t \in [-a,a]^2}{\mathbb E} \left[ \|\partial_{12}X(t)\|^q\right]^{1/q} \nonumber\\
&+4a^2\sup_{t \in [-a,a]^2}{\mathbb E} \left[ \|D^2 f(X(t))\|^p\right]^{1/p} \sup_{t \in [-a,a]^2} {\mathbb E} \left[ \|\partial_1X(t)\|^{2q}\right]^{1/(2q)} \sup_{t \in [-a,a]^2} {\mathbb E} \left[ \|\partial_2X(t)\|^{2q}\right]^{1/(2q)}.\nonumber
\end{align}
Let $X(t)=(Z(t),\partial_1 Z(t),\partial_2 Z(t))\in{\mathbb R}^3$, so that $\|X(t)\|^2= Z(t)^2 + \partial_1 Z(t)^2 + \partial_2 Z(t)^2$. Under the conditions of Lemma \ref{lem-gauss-crofton}, $X$ has $C^2$ sample paths and, since $\partial_i X$ and $\partial_{12}X$ are centered Gaussian vectors, for all $q>0$ there exists a constant $C_q$, depending only on $q$ and $\overline{C}$, such that 
\begin{align}\label{lem-gauss-crofton-2} 
&\sup_{i=1,2} \sup_{t \in [-a,a]^2} {\mathbb E} \left[ \|\partial_iX(t)\|^{q}\right]^{1/q}\leq C_q,\\
&\sup_{t \in [-a,a]^2}{\mathbb E} \left[ \|\partial_{12}X(t)\|^q\right]^{1/q} \leq C_q.\nonumber
\end{align}
Moreover, since $X(t)$ is a centered Gaussian vector with covariance matrix $\Gamma_t$, for $\kappa \in [0,\tfrac{3}{2})$,
\begin{align}\label{lem-gauss-crofton-3} 
\sup_{t \in [-a,a]^2}{\mathbb E} \left[ \|X(t)\|^{-2\kappa}\right] &= \sup_{t \in [-a,a]^2}\int_{{\mathbb R}^3} \frac{1}{(2\pi)^{3/2}} \|\sqrt{\Gamma_t}\, x\|^{-2\kappa} e^{-\frac{\|x\|^2}{2}}dx \nonumber\\
&\leq \frac{1}{\underline{C}^{\kappa}}\int_{{\mathbb R}^3} \frac{1}{(2\pi)^{3/2}} \| x\|^{-2\kappa} e^{-\frac{\|x\|^2}{2}}dx\nonumber\\
&\leq \frac{C_{\kappa}}{\underline{C}^{\kappa}}.
\end{align}
Let $r\in(0,1)$, $\varepsilon>0$ and $f^{\varepsilon} (x)= (\varepsilon +\|x\|^2)^{-\frac{r}{2}}$, $x\in{\mathbb R}^3$. Then
\begin{align}\label{lem-gauss-crofton-4} 
&|f^{\varepsilon} (x)|\leq  \|x\|^{-r}, \\
&\|\nabla f^{\varepsilon}(x)\| \leq  r\,\|x\|^{-(r+1)},\nonumber\\
&\|D^2 f^{\varepsilon}(x)\| \leq  r(r+3)\,\|x\|^{-(r+2)}.\nonumber
\end{align}
Let $p\in\big(1,\tfrac{3}{r+2}\big)$ and let $q$ be its conjugate exponent. Since $\tfrac{r}{2}$, $\tfrac{(r+1)p}{2}$ and $\tfrac{(r+2)p}{2}$ belong to $[0,\tfrac32)$, plugging estimates \eqref{lem-gauss-crofton-2} to \eqref{lem-gauss-crofton-4} into \eqref{lem-gauss-crofton-1}, we derive that there exists a constant $C_{r}$, depending only on $\underline{C}$, $\overline{C}$ and $r$, and not on $\varepsilon$, such that
\begin{align*}
{\mathbb E} \left[\sup_{t\in [-a,a]^2} \left( \frac{1}{\varepsilon+ Z(t)^2 + \partial_1 Z(t)^2 + \partial_2 Z(t)^2} \right)^{\frac{r}{2}}\right] \leq C_{r}(1+a+ a^2). 
\end{align*}
Letting $\varepsilon$ go to $0$, the monotone convergence theorem yields the desired estimate.\end{proof}

\subsubsection{Proof of Lemma \ref{g-crofton}}

Since the covariance matrix of $(G(t),\partial_1 G(t),\partial_2 G(t))$ is given by
\begin{align*}
\left( \begin{array}{ccc}1 & 0 &0\\
0 &\frac{1}{3} &\frac{1}{4}\\
0& \frac{1}{4}& \frac{1}{3}\end{array} \right),
\end{align*}
whose eigenvalues are $1$, $\frac{7}{12}$ and $\frac{1}{12}$, and since, by the Itô isometry, for all $n\in\{2,3\}$ and $j\in\{1,2\}^n$,
\begin{align*}
{\mathbb E}\left[\left(\partial_{t_{j_1},\dots,t_{j_n}} G(t)\right)^2\right]=\iint_{[0,1]^2}\prod_{l=1}^n u_{j_l}^2\,du_1\,du_2\leq \frac{1}{5},
\end{align*}
$G$ satisfies the assumptions of Lemma~\ref{lem-gauss-crofton} with $\underline{C}=\frac{1}{12}$ and $\overline{C}=1$, which are independent of $a$. Thus, by Markov's inequality, for all $r\in(0,1)$ there exists a constant $C_r$ such that 
\begin{align*}
{\mathbb P}\left[ \min_{t\in B(0,a)} \left(|G(t)|+\|\nabla G(t)\|\right) \leq \delta\right] \leq C_r (1+a^2)\, \delta^{r}.
\end{align*}
%
Moreover, since ${\mathcal D}_{0,G,a}\subset{\mathcal D}_{\delta,G,a}$ for all $\delta>0$, letting $\delta$ go to $0$ gives ${\mathbb P}\left[{\mathcal D}_{0,G,a}\right]=0$.

\subsubsection{Proof of Lemma \ref{Pm-crofton}}
Let $0<\varepsilon <\delta$ and let $H_{\delta,\varepsilon}$ be defined in \eqref{def-H-delta-epsi}. Since $\mathbf{1}_{[0,\delta]}\leq H_{\delta,\varepsilon}\leq \mathbf{1}_{[0,\delta+\varepsilon]}$ on $[0,\infty)$, we have
\begin{align*}
{\mathbb P}  \left[ {\mathcal D}_{\delta,P_m,a}
\right] &\leq {\mathbb E}\left[ H_{\delta,\varepsilon}\big( \min_{t \in B(0,a)}(|P_m(t)| +\|\nabla P_m(t)\|)\big) \right]\\
&\leq {\mathbb E}\left[ H_{\delta,\varepsilon}\big( \min_{t \in B(0,a)}(|P_m(t)| +\|\nabla P_m(t)\|)\big) \right]\\&-{\mathbb E}\left[ H_{\delta,\varepsilon}\big( \min_{t \in B(0,a)}(|G(t)| +\|\nabla G(t)\|)\big) \right]\\&+ {\mathbb P}  \left[ \min_{t \in B(0,a)}(|G(t)| +\|\nabla G(t)\|)\leq \varepsilon+ \delta \right].
\end{align*}
Thus, using Corollary \ref{cor-P-m-cv-G} for the $\frac{1}{\varepsilon}$-Lipschitz continuous function 
$g \mapsto H_{\delta,\varepsilon} \big( \min_{t \in B(0,a)}(|g(t)| +\|\nabla g(t)\|)\big)$ and Lemma \ref{g-crofton}, we obtain
\begin{align*}
 {\mathbb P}  \left[ {\mathcal D}_{\delta,P_m,a}
\right] &\leq C\left(\frac{1}{\varepsilon m^{\gamma}} + (\delta + \varepsilon)^{r}\right)\\
&\leq C\left(\frac{1}{\varepsilon m^{\gamma}} + \delta ^{r}\right),
\end{align*}
{where $C$ depends only on $p$, $\gamma$, $r$ and $a$.
Choosing $\varepsilon=\frac{\delta}{2}$, for all $m\geq\delta^{-\frac{1+r}{\gamma}}$, we have $\frac{1}{\varepsilon m^{\gamma}}\leq 2\delta^{r}$, which yields the inequality above.
Finally, let $m\in\mathbb N$ and $\delta=m^{-\frac{\gamma}{1+r}}\in(0,1]$. Then $m=\delta^{-\frac{1+r}{\gamma}}$ and, since ${\mathcal D}_{0,P_m,a}\subset{\mathcal D}_{\delta,P_m,a}$,
\begin{align*}
{\mathbb P}\left[{\mathcal D}_{0,P_m,a}\right]\leq C\,\delta^{r}= C\, m^{-\frac{\gamma r}{1+r}}.
\end{align*}}

\section{Moment bounds for the number of zeros}\label{moment_bounds}

In this appendix, we study the number of zeros of the restrictions of the random polynomial $P_m$ to the lines $L_{\theta,\rho}$. According to the beginning of Section~8 of \cite{NguyenVu}, this number of zeros has a finite expectation. Here, we prove that all its moments are bounded uniformly in $m$, $\theta$ and $\rho$. 

{Throughout this appendix, we use the convention $N(0,I):=0$ for every interval $I$, where $0$ denotes the null function.}

To simplify the notation, we denote by
\[
P_{\theta,\rho}^m:=P^m\circ L_{\theta,\rho}
\]
the restriction of the random polynomial defined in \eqref{def:pol_1} to the parametrized line \eqref{Lthetarho}. That is,
\begin{align}\notag
P_{\theta,\rho}^m(\beta)
=
\frac{1}{m}
\sum_{k,l=1}^{m}
\operatorname{Re}
\left[
(a_{k,l}+ib_{k,l})
\,e^{-i\gamma_{k,l}}
\exp\!\left(-i\cc\,\beta\right)
\right],
\end{align}
where the random variables \(a_{k,l}\) and \(b_{k,l}\) are centered, independent, identically distributed, and have unit variance,
\begin{equation*}
\gamma_{k,l}
=
\rho\left(
\frac{k}{m}\sin\theta
-
\frac{l}{m}\cos\theta
\right),
\end{equation*}
and
\begin{equation}\label{cc}
\cc
:=
\frac{k\cos(\theta)+l\sin(\theta)}{m}.
\end{equation}

\begin{proposition}\label{prop-maj-expec-N}
Let $p\geq 1$ and $a>0$. There exists a constant $C(p,a)>0$, depending only on $p$ and $a$, such that for every interval $I\subset[-a,a]$,
\begin{align*}
\sup_m \sup_{\theta \in [0,2\pi],\,\rho \in [-a,a]}
\E\left(N\left(P_{\theta,\rho}^m, I \right)^p\right)
\leq C(p,a).
\end{align*}
\end{proposition}

The proof of Proposition~\ref{prop-maj-expec-N} relies on a detailed study of the distributions of
\(
P_{\theta,\rho}^{m(j)}(\beta),
\)
\(j\geq0\).
We begin by obtaining estimates on the characteristic function of
\(
P_{\theta,\rho}^{m}(\beta),
\)
which will later be used to derive quantitative information on these distributions.

\subsection{Characteristic function of $P_{\theta,\rho}^m(\beta)$}

Let $j\geq 0$ be an integer. Then the $j$th derivative of \(P_{\theta,\rho}^m\) is given by
\begin{align}\notag
P_{\theta,\rho}^{m(j)}(\beta)
=
\frac{1}{m}
\sum_{k,l=1}^{m}
\operatorname{Re}
\left[
(a_{k,l}+ib_{k,l})
(\cc)^j
e^{-i\gamma_{k,l}}
\exp\left(
-i\left(\cc\,\beta+\frac{j\pi}{2}\right)
\right)
\right].
\end{align}

Consequently,
\begin{align}\label{der-pol}
P_{\theta,\rho}^{m(j)}(\beta)
=&\frac{1}{m}
\sum_{k=1}^{m}
\sum_{l=1}^{m}
a_{k,l}
(\cc)^j
\cos\left(
\gmm+\cc\,\beta+\frac{j\pi}{2}
\right)
\notag\\
&-
\frac{1}{m}
\sum_{k=1}^{m}
\sum_{l=1}^{m}
b_{k,l}
(\cc)^j
\sin\left(
\gmm+\cc\,\beta+\frac{j\pi}{2}
\right).
\end{align}

Let us introduce the following notation:
\begin{align}
A_{k,l,\theta,\rho,\beta}
&:=
\cos\left(
\gmm+\cc\,\beta+\frac{j\pi}{2}
\right),
\label{not-der-arg_1}
\\
B_{k,l,\theta,\rho,\beta}
&:=
\sin\left(
\gmm+\cc\,\beta+\frac{j\pi}{2}
\right).
\label{not-der-arg_2}
\end{align}

Then,
\begin{align*}
P_{\theta,\rho}^{m(j)}(\beta)
=&\frac{1}{m}
\sum_{k=1}^{m}
\sum_{l=1}^{m}
a_{k,l}
(\cc)^j
A_{k,l,\theta,\rho,\beta}
\\
&-
\frac{1}{m}
\sum_{k=1}^{m}
\sum_{l=1}^{m}
b_{k,l}
(\cc)^j
B_{k,l,\theta,\rho,\beta}.
\end{align*}
We now compute the characteristic function of \(P_{\theta,\rho}^{m(j)}(\beta)\).

Since \(a_{k,l}\) and \(b_{k,l}\) are independent and identically distributed, we have
\begin{align*}
&\E\left(\exp\left(ivP_{\theta,\rho}^{m(j)}(\beta)\right)\right)
=
\prod_{k,l=1}^{m}
\phi\left(\frac{v}{m}(\cc)^jA_{k,l,\theta,\rho,\beta}\right)
\phi\left(\frac{v}{m}(\cc)^jB_{k,l,\theta,\rho,\beta}\right),
\end{align*}
where \(\phi\) denotes the common characteristic function of \(a_{k,l}\) and \(b_{k,l}\).

We next derive key properties of the characteristic function, which will be essential in the estimation of Lemma \ref{lem-maj-phi-gamma-r}.

Since
$$\{(a_{k,\ell}, b_{k,\ell}) : k,\ell =1,\ldots,m\}$$
are independent and identically distributed centered random vectors with identity covariance matrix, there exists a constant \(c>0\) such that
\begin{align}\label{phi<gaus}
| \phi(v)| \leq \exp\left(- \frac{v^2}{4}\right),
\qquad \forall\, |v| \leq c.
\end{align}

For \(r\in\{2,\ldots,2m\}\), let
\begin{align*}
C_r
:=
\left\{
(k,l)\in\{1,\ldots,m\}^2
\,:\,
k+l\le r
\right\}.
\end{align*}

\begin{lemma}
For every \(r\in\{2,\ldots,2m\}\) and every
\[
v\leq \frac{cm^{j+1}}{r^j},
\]
we have
\begin{align*}
\left| \E\left[\exp\left(ivP_{\theta,\rho}^{m(j)}(\beta)\right)\right] \right|
\leq
\exp\left(
-\frac{v^2}{4m^2}
\sum_{(k,l)\in C_r}
\cc^{2j}
\right).
\end{align*}
\end{lemma}

\begin{proof}
From the definition of \(\cc\) given in \eqref{cc}, and the expressions for
\(A_{k,l,\theta,\rho,\beta}\) and \(B_{k,l,\theta,\rho,\beta}\) given in
\eqref{not-der-arg_1} and \eqref{not-der-arg_2}, we have
\begin{align*}
|A_{k,l,\theta,\rho,\beta}|^2
+
|B_{k,l,\theta,\rho,\beta}|^2
&=1,
\\
|\cc|
&\leq \frac{k+l}{m}.
\end{align*}

Let
\[
v\leq \frac{cm^{j+1}}{r^j}.
\]
If \((k,l)\in C_r\), then \(k+l\le r\), and therefore
\[
|\cc|
\le \frac{r}{m}.
\]
Hence
\[
\left|
\frac{v}{m}
(\cc)^j
A_{k,l,\theta,\rho,\beta}
\right|
\le
\frac{v}{m}
\left(\frac{r}{m}\right)^j
\le c,
\]
and similarly
\[
\left|
\frac{v}{m}
(\cc)^j
B_{k,l,\theta,\rho,\beta}
\right|
\le c.
\]

Therefore, for \((k,l)\in C_r\), we may apply \eqref{phi<gaus} and obtain
\begin{align*}
&
\left|
\phi\left(
\frac{v}{m}(\cc)^jA_{k,l,\theta,\rho,\beta}
\right)
\phi\left(
\frac{v}{m}(\cc)^jB_{k,l,\theta,\rho,\beta}
\right)
\right|
\\
&\le
\exp\left(
-\frac{v^2}{4m^2}
(\cc)^{2j}
A_{k,l,\theta,\rho,\beta}^{\,2}
\right)
\exp\left(
-\frac{v^2}{4m^2}
(\cc)^{2j}
B_{k,l,\theta,\rho,\beta}^{\,2}
\right)
\\
&=
\exp\left(
-\frac{v^2}{4m^2}
(\cc)^{2j}
\bigl(
A_{k,l,\theta,\rho,\beta}^{\,2}
+
B_{k,l,\theta,\rho,\beta}^{\,2}
\bigr)
\right)
\\
&=
\exp\left(
-\frac{v^2}{4m^2}
(\cc)^{2j}
\right).
\end{align*}

For \((k,l)\notin C_r\), we simply use the trivial bound
\[
|\phi(u)|\le 1.
\]

Combining these estimates with the product representation of the characteristic function yields
\begin{align*}
\left| \E\left[\exp\left(ivP_{\theta,\rho}^{m(j)}(\beta)\right)\right] \right|
&\le
\exp\left(
-\frac{v^2}{4m^2}
\sum_{(k,l)\in C_r}
\cc^{2j}
\right),
\end{align*}
which concludes the proof.
\end{proof}

We now derive a lower bound for
\begin{equation*}
\sum_{(k,l)\in C_r}\cc^{2j}.
\end{equation*}

\begin{lemma}
Let \(j\geq 0\), \(12\leq r\leq 2m\), and
\(\theta\in[0,2\pi]\). Then
\begin{align*}
m^{2j}\sum_{(k,l)\in C_r}\cc^{2j}
\geq
\frac{1}{2^{5j+4}(2j+1)}
\,r^{2j+2}.
\end{align*}
\end{lemma}

\begin{proof}
First, note that
\begin{align*}
m^2\cc^2
=
k^2(\cos\theta)^2
+
l^2(\sin\theta)^2
+
2kl\sin\theta\cos\theta.
\end{align*}
Thus,
\begin{align*}
m^2\cc^2
\geq
\left[
k|\cos\theta|
-
l|\sin\theta|
\right]^2.
\end{align*}

Note that either \( |\cos \theta|^2 \geq \frac{1}{2} \) or \( |\sin \theta|^2 \geq \frac{1}{2} \). We assume, for instance, that \( |\cos \theta|^2 \geq \frac{1}{2} \). Thus,
\begin{align*}
\left(k\cos \theta + l\sin \theta\right)^2
\geq
\left(k|\cos \theta|- l|\sin \theta|\right)^2
\geq
\frac{1}{2}|k-l|^2,
\qquad \forall\, k\geq l.
\end{align*}

Recall that
\begin{equation*}
C_r
=
\{(k,l)\in\{1,\ldots,m\}^2:\ k+l\le r\}.
\end{equation*}
Thus,
\begin{align*}
\sum_{(k,l)\in C_r}|k-l|^{2j}
&\geq
\sum_{k=1}^{[r/2]}
\sum_{l=1}^{k}
|k-l|^{2j}
\\
&\geq
\frac12
\sum_{k=1}^{[r/2]}
\sum_{l=1}^{[r/2]}
|k-l|^{2j}.
\end{align*}

We perform the change of variable $\tilde l=k-l.$

Thus,
\begin{align*}
\sum_{(k,l)\in C_r}|k-l|^{2j}
&\geq
\frac12
\sum_{k=1}^{[r/2]}
\sum_{\tilde l=k-[r/2]}^{k-1}
|\tilde l|^{2j}.
\end{align*}

Using Fubini's theorem,
\begin{align*}
\sum_{(k,l)\in C_r}|k-l|^{2j}
&\geq
\frac12
\sum_{\tilde l=1-[r/2]}^{[r/2]-1}
\sum_{k=\tilde l+1}^{\tilde l+[r/2]}
|\tilde l|^{2j}
\\
&\geq
[r/2]
\sum_{\tilde l=1}^{[r/2]-1}
|\tilde l|^{2j}
\\
&\geq
\frac{1}{2j+1}
([r/2]-2)^{2j+2}.
\end{align*}

We have used in the last inequality that
\begin{equation*}
\sum_{l=0}^{n} l^i
\geq
\frac{1}{i+1}(n-1)^{i+1}.
\end{equation*}

For $r\geq 12$,
\begin{equation*}
[r/2]-2
\geq
\frac r2-3
\geq
\frac r4,
\end{equation*}
thus
\begin{align*}
m^{2j}\sum_{(k,l)\in C_r}\cc^{2j}
\geq
\frac{1}{2^{5j+4}(2j+1)}
\,r^{2j+2}.
\end{align*}

If $ |\cos \theta|^2 \leq \frac12 $, we exchange the roles of $k$ and $l$.
\end{proof}

We summarize the results of this section in the following lemma.

\begin{lemma}\label{lem-maj-phi-gamma-r}
For all \(v\in\mathbb R\),
\begin{align*}
\left|
\E\left(
\exp\left(ivP_{\theta,\rho}^{m(j)}(\beta)\right)
\right)
\right|
\leq 1.
\end{align*}

There exists a constant \(c>0\) such that, for every \(j\geq 0\), there exists a constant \(\gamma_j>0\) (depending only on \(j\) and independent of \(\theta\), \(\beta\), \(\rho\), \(a\), and \(m\)) such that, for all \(12\leq r\leq 2m\) and
\[
v\leq \frac{cm^{j+1}}{r^j},
\]
we have
\begin{align*}
\left|
\E\left(
\exp\left(ivP_{\theta,\rho}^{m(j)}(\beta)\right)
\right)
\right|
\leq
\exp\left(
-\gamma_j
\frac{v^2r^{2j+2}}{m^{2j+2}}
\right).
\end{align*}
\end{lemma}

\subsection{Distribution function of $|P_{\theta,\rho}^{m(j)}(\beta)|$}

We now derive distributional estimates for
\(
P_{\theta,\rho}^{m(j)}(\beta).
\)
The following result is adapted from \cite[Lemma 4]{Flasche}.

\begin{lemma}\label{poly_lem1}
Let $j \in \mathbb{N}_0$ be fixed, and let $I$ be an interval. There exists a constant $C(j) > 0$, depending only on $j$, such that for all
\[
m\in\mathbb N,\qquad
T>0,\qquad
\theta\in[0,2\pi],\qquad
\rho\in[-a,a],
\qquad\text{and}\qquad
\beta\in I,
\]
the random variable \(P_{\theta,\rho}^{m(j)}(\beta)\) satisfies
\begin{align}
\P\left(|P_{\theta,\rho}^{m(j)}(\beta)|\leq T\right)
\leq
C(j)\left(T+T^{-1/2}m^{-(j+1)/2}\right).
\end{align}
\end{lemma}

\begin{proof}
The value of the constant $C(j)$ may change from line to line, but it always depends only on $j$.
Let $\eta$ be a random variable independent of $P_{\theta,\rho}^{m(j)}$, with characteristic function
\[
\psi(v)=\frac{\sin^2(v\lambda)}{v^2\lambda^2},
\qquad \lambda>0.
\]
Define
\[
\tilde{Y}_m^{(j)}(\beta)
=
P_{\theta,\rho}^{m(j)}(\beta)+\eta.
\]
Then, for $T>0$,
\begin{align}\label{den:Yj}
\P\left(|P_{\theta,\rho}^{m(j)}(\beta)| \leq T\right)
\leq
\P\left(|\tilde{Y}_m^{(j)}(\beta)| \leq \frac{3T}{2}\right)
+
\P\left(|\eta| \geq \frac{T}{2}\right).
\end{align}

To estimate the first term on the right-hand side of \eqref{den:Yj}, we consider the characteristic function of \(\tilde{Y}_m^{(j)}(\beta)\), namely
\begin{align}\notag
\tilde{\psi}_m(v)
=
\E\left[
\exp\left(iv P_{\theta,\rho}^{m(j)}(\beta)\right)
\right]
\psi(v).
\end{align}

By Fourier inversion,
\begin{align*}
\P\left(|\tilde{Y}_m^{(j)}(\beta)| \leq y\right)
&=
\frac{2}{\pi}
\int_0^{\infty}
\frac{\sin(yv)}{v}
\,\operatorname{Re}\!\left(\tilde{\psi}_m(v)\right)
\,dv
\\
&\leq
\frac{2y}{\pi}
\int_0^{\infty}
\psi(v)
\left|
\E\left(
\exp\left(iv P_{\theta,\rho}^{m(j)}(\beta)\right)
\right)
\right|
\,dv.
\end{align*}

To estimate this integral, we split the integration domain as follows:
\begin{align*}
\P\left(|\tilde{Y}_m^{(j)}(\beta)|\leq y\right)
&\leq
\frac{2y}{\pi}
\left(
I_0
+
\sum_{r=1}^{2m-12} I_r
+
I_{2m-11}
\right).
\end{align*}

where
\begin{align}\label{aux:int}
I_r
=
\int_{\Gamma_r}
\psi(v)
\left|
\E\left(
\exp\left(ivP_{\theta,\rho}^{m(j)}(\beta)\right)
\right)
\right|
\,dv.
\end{align}
We consider the following partition of \(\mathbb{R}_+\), where \(c>0\) is the constant given by Lemma~\ref{lem-maj-phi-gamma-r}:
\begin{align*}
\Gamma_{0}
:=&
\left\{
v:\;
0\leq v<
\frac{cm^{j+1}}{(2m)^j}
\right\},
\\
\Gamma_{r}
:=&
\left\{
v:\;
\frac{cm^{j+1}}{(2m-r+1)^j}
\leq v<
\frac{cm^{j+1}}{(2m-r)^j}
\right\},
\qquad
r=1,\ldots,2m-12,
\\
\Gamma_{2m-11}
:=&
\left\{
v:\;
v\geq
\frac{cm^{j+1}}{12^j}
\right\}.
\end{align*}

For the integral on \(\Gamma_{2m-11}\), we use that \(|\phi|\leq 1\) and \(\sin^2(v\lambda)\leq 1\). Thus,
\begin{align*}
I_{2m-11}
&\leq
\int_{\frac{cm^{j+1}}{12^j}}^{\infty}
\frac{dv}{v^2\lambda^2}
=
\frac{12^j}{c\lambda^2m^{j+1}}.
\end{align*}

In the case of \(\Gamma_{0}\), we use Lemma~\ref{lem-maj-phi-gamma-r} with \(r=2m\) and derive
\begin{align*}
I_{0}
&\leq
\int_{0}^{\frac{cm^{j+1}}{(2m)^j}}
\exp\left(
-\frac{(2m)^{2j+2}\gamma_j}{m^{2j+2}}\,v^2
\right)\,dv \\
&\leq
\int_{0}^{\infty}
\exp\left(
-\gamma_j2^{2j+2}v^2
\right)\,dv \\
&=
\frac{\sqrt{\pi}}{2}\,\tilde{\gamma}_j,
\end{align*}
where \(\tilde{\gamma}_j\) is a constant depending only on \(j\) and independent of \(m\), \(\rho\), \(\theta\), and \(\beta\).

Lastly, for \(\Gamma_r\), by applying Lemma~\ref{lem-maj-phi-gamma-r}, we obtain
\begin{align*}
I_{r}
&\leq
\int_{\Gamma_{r}}
\psi(v)
\exp\left(
-v^2 \gamma_j
\frac{(2m-r)^{2j+2}}{m^{2j+2}}
\right)
\,dv.
\end{align*}

We perform the change of variable
$$u=\left(\frac{2m-r}{m}\right)^{j+1}v,$$
and thus
\begin{align*}
I_r
&\leq
\frac{1}{\lambda^2}
\left(\frac{2m-r}{m}\right)^{j+1}
\int_{\tilde{\Gamma}_r}
\frac{1}{u^2}
e^{-\gamma_j u^2}
\,du,
\end{align*}
where
$$\tilde{\Gamma}_r
=
\left[
c\frac{(2m-r)^j}{(2m-r+1)^j}(2m-r),
\;
c(2m-r)
\right].$$

The function $u\mapsto u^{-2}e^{-\gamma_j u^2}$
is decreasing, and thus
\begin{align*}
I_r
&\leq
\frac{1}{c\lambda^2}
\left(\frac{2m-r}{m}\right)^{j+1}
\left(\frac{2m-r+1}{2m-r}\right)^{2j}
e^{-c^2\gamma_j
\left(\frac{2m-r+1}{2m-r}\right)^{2j}
(2m-r)^2}.
\end{align*}
Noting that
$$\min_{r=1,\ldots,2m-12}
\frac{2m-r}{2m-r+1}
=
\frac{12}{13}
>
\frac12,$$
we obtain
\begin{align*}
I_r
&\leq
\frac{1}{c\lambda^2}
\left(\frac{2m-r}{m}\right)^{j+1}
\left(\frac{2m-r+1}{2m-r}\right)^{2j}
e^{-c^2\gamma_j2^{-2j}(2m-r)^2}
\\
&\leq
\frac{m^{-j-1}}{c\lambda^2}
(2m-r)^{j+1}
2^{2j}
e^{-c^2\gamma_j2^{-2j}(2m-r)^2}.
\end{align*}
Setting \(k=2m-r\), the previous estimate yields
$$I_r
\le
\frac{2^{2j}}{c\lambda^2m^{j+1}}
k^{j+1}
e^{-c^2\gamma_j2^{-2j}k^2}.$$

Then, applying the ratio test, we have
$$\lim_{k\to \infty}
\left(1+\frac{1}{k}\right)^{j+1}
\exp\left(
-c^2\gamma_j2^{-2j}(2k+1)
\right)
=
0,$$
therefore, we can conclude that
\begin{align*}
\sum_{r=1}^{2m-12}I_r
\leq
\frac{C(j)}{\lambda^2m^{j+1}}.
\end{align*}

Putting together all the estimates for
\(I_0,\ldots,I_{2m-11}\)
and taking $
y=\frac{3T}{2},$
we have
\begin{align*}
\P\left(
|\tilde{Y}_m^{(j)}(\beta)|
\leq
\frac{3T}{2}
\right)
\leq
T\,C(j)
\left[
1+\frac{1}{\lambda^2m^{j+1}}
\right].
\end{align*}

To bound the second term in \eqref{den:Yj}, we use Chebyshev's inequality:
\begin{align*}
\P\left(|\eta|\geq \frac{T}{2}\right)
\leq
\frac{4\E(|\eta|^2)}{T^2}
=
\frac{4\lambda^2}{3T^2}.
\end{align*}

Thus,
\begin{align*}
\P\left(|P_{\theta,\rho}^{m(j)}(\beta)|\leq T\right)
&\leq
C(j)
\left[
T\left(
1+\frac{1}{\lambda^2 m^{j+1}}
\right)
+
\frac{\lambda^2}{T^2}
\right].
\end{align*}

Choosing
$$\lambda=T^{3/4}m^{-(j+1)/4},$$
we obtain
\begin{align*}
\P\left(|P_{\theta,\rho}^{m(j)}(\beta)|\leq T\right)
&\leq
C(j)
\left(
T+T^{-1/2}m^{-(j+1)/2}
\right).
\end{align*}
\end{proof}

\begin{lemma}\label{2ndmomt}
For all \(m\geq 1\), \(\theta\in[0,2\pi]\), \(\rho\in[-a,a]\), \(\beta\in\mathbb{R}\), and \(j\geq 0\), we have
\begin{align*}
\mathbb{E}\left(P_{\theta,\rho}^{m(j)}(\beta)\right)=0
\qquad\text{and}\qquad
\mathbb{E}\left(\left|P_{\theta,\rho}^{m(j)}(\beta)\right|^2\right)
\leq
\frac{2^{2j}}{2j+1}.
\end{align*}
\end{lemma}

\begin{proof}
Since the random variables \(a_{k,\ell}\) and \(b_{k,\ell}\) are centered, we have $\mathbb{E}\left(P_{\theta,\rho}^{m(j)}(\beta)\right)=0.$ Moreover, since the collection
\(
\{(a_{k,\ell},b_{k,\ell})\}_{k,\ell}
\)
consists of independent random vectors with identity covariance matrix, and
$A_{k,\ell,\theta,\rho,\beta}^2
+
B_{k,\ell,\theta,\rho,\beta}^2
=
1,$ it follows that
\begin{align*}
\mathbb{E}\left(\left|P_{\theta,\rho}^{m(j)}(\beta)\right|^2\right)
&=
\frac{1}{m^2}
\sum_{k,l=0}^{m-1}
\cc^{2j}
\Big(
A_{k,l,\theta,\rho,\beta}^{2}
+
B_{k,l,\theta,\rho,\beta}^{2}
\Big) \\
&=
\frac{1}{m^2}
\sum_{k,l=0}^{m-1}
\cc^{2j} \\
&=
\frac{1}{m^{2(j+1)}}
\sum_{k,l=0}^{m-1}
\left(
k\cos\theta+l\sin\theta
\right)^{2j} \\
&\leq
\frac{2^j}{m^{2(j+1)}}
\sum_{k,l=0}^{m-1}
(k^2+l^2)^j.
\end{align*}

Using $(k^2+l^2)^j \leq
2^{j-1}(k^{2j}+l^{2j}),$
we obtain
\begin{align*}
\frac{1}{m^{2(j+1)}}
\sum_{k,l=0}^{m-1}
(k^2+l^2)^j
&\leq
\frac{2^{j-1}}{m^{2(j+1)}}
\sum_{k,l=0}^{m-1}
(k^{2j}+l^{2j}) \\
&=
\frac{2^j}{m^{2j+1}}
\sum_{k=0}^{m-1}
k^{2j}.
\end{align*}

Since
$$\sum_{k=0}^{m-1}k^{2j}
\leq
\int_0^m x^{2j}\,dx
=
\frac{m^{2j+1}}{2j+1},$$
it follows that
$$\frac{1}{m^{2(j+1)}}
\sum_{k,l=0}^{m-1}
(k^2+l^2)^j
\leq
\frac{2^j}{2j+1}.$$

Therefore,
\begin{equation*}
\mathbb{E}\left(\left|P_{\theta,\rho}^{m(j)}(\beta)\right|^2\right)
\leq
\frac{2^{2j}}{2j+1}.
\end{equation*}
\end{proof}

\subsection{Distribution of the number of zeros of $P_{\theta,\rho}^{m(j)}$}

\begin{lemma}\label{intrs}
Let \(0 \leq \alpha \leq \zeta\). Define the following event on \((\Omega,\mathcal F)\):
\begin{align*}
D_k^{(j)}
=
\Bigl\{
\omega\in\Omega
\,\Big|\,
P_{\theta,\rho}^{m(j)}(\omega,\cdot)
\text{ has at least }k\text{ zeros in }
I=[\alpha,\zeta]
\Bigr\}.
\end{align*}

\begin{description}
\item[i)]
\begin{align*}
{\mathbb P}\left(D_k^{(j)} \cap \{|P_{\theta,\rho}^{m(j)}| \geq T\}\right)
\leq
\frac{2^{2(j+k)}}{T^2(2(j+k)+1)}
\left(
\frac{(\zeta-\alpha)^k}{k!}
\right)^2.
\end{align*}

\item[ii)]
\begin{align*}
\P\left(D_k^{(j)}\right)
\leq
\frac{2^{2(j+k)}}{T^2(2(j+k)+1)}
\left(
\frac{(\zeta-\alpha)^k}{k!}
\right)^2
+
C(j)\left(
T+T^{-1/2}m^{-(j+1)/2}
\right).
\end{align*}
\end{description}
\end{lemma}

\begin{proof}
Case i): On the event \(D_k^{(j)}\), by Rolle's theorem, we can find some points
\begin{equation*}
\zeta \geq t_0 \geq t_1 \geq \cdots \geq t_{k-1} \geq \alpha
\end{equation*}
such that
\begin{equation*}
P_{\theta,\rho}^{m(j+l)}(t_l)=0,
\qquad
l=0,\ldots,k-1.
\end{equation*}

We may then consider the random variable
\begin{align*}
Y_m^{(j)}
=
{\mathbf 1}_{D_k^{(j)}}
\int_{t_0}^{\zeta}
\int_{t_1}^{x_1}
\cdots
\int_{t_{k-1}}^{x_{k-1}}
P_{\theta,\rho}^{m(j+k)}(x_k)
\,dx_k
\cdots
dx_1.
\end{align*}

On the event \(D_k^{(j)}\), the random variables
\(Y_m^{(j)}\) and \(P_{\theta,\rho}^{m(j)}\) are equal, and on the complement of the event \(D_k^{(j)}\), we have \(Y_m^{(j)}=0\). Hence,
\begin{align*}
{\mathbb P}\left(
D_k^{(j)}
\cap
\{|P_{\theta,\rho}^{m(j)}|\geq T\}
\right)
\leq
{\mathbb P}\left(
\{|Y_m^{(j)}|\geq T\}
\right).
\end{align*}

By Markov's inequality,
\begin{align*}
&{\mathbb P}\left(
D_k^{(j)}
\cap
\{|P_{\theta,\rho}^{m(j)}| \geq T\}
\right) \\
&\leq
\frac{1}{T^2}
\,{\mathbb E}\!\left[\left(Y_m^{(j)}\right)^2\right].
\end{align*}

Since \({\mathbf 1}_{D_k^{(j)}} \leq 1\), we obtain
\begin{align*}
&{\mathbb P}\left(
D_k^{(j)}
\cap
\{|P_{\theta,\rho}^{m(j)}| \geq T\}
\right) \\
&\leq
\frac{1}{T^2}
{\mathbb E}\Bigg(
\int_{t_0}^{\zeta}
\int_{t_0}^{\zeta}
\int_{t_1}^{x_1}
\int_{t_1}^{y_1}
\cdots
\int_{t_{k-1}}^{x_{k-1}}
\int_{t_{k-1}}^{y_{k-1}}
P_{\theta,\rho}^{m(j+k)}(x_k)
P_{\theta,\rho}^{m(j+k)}(y_k)
\,dy_k\,dx_k
\cdots
dy_1\,dx_1
\Bigg).
\end{align*}

By Fubini's theorem and the Cauchy--Schwarz inequality,
\begin{align*}
&{\mathbb E}\Big[
P_{\theta,\rho}^{m(j+k)}(x_k)
P_{\theta,\rho}^{m(j+k)}(y_k)
\Big] \\
&\leq
\Big(
{\mathbb E}|P_{\theta,\rho}^{m(j+k)}(x_k)|^2
\Big)^{1/2}
\Big(
{\mathbb E}|P_{\theta,\rho}^{m(j+k)}(y_k)|^2
\Big)^{1/2}.
\end{align*}

Therefore, by Lemma~\ref{2ndmomt},
\begin{align*}
&{\mathbb P}\left(
D_k^{(j)}
\cap
\{|P_{\theta,\rho}^{m(j)}|\geq T\}
\right)\\
&\leq
\frac{2^{2(j+k)}}{T^2(2(j+k)+1)}
\left(
\int_{t_0}^{\zeta}
\int_{t_1}^{x_1}
\cdots
\int_{t_{k-1}}^{x_{k-1}}
dx_k\cdots dx_1
\right)^2.
\end{align*}

Since
\[
\int_{t_0}^{\zeta}
\int_{t_1}^{x_1}
\cdots
\int_{t_{k-1}}^{x_{k-1}}
dx_k\cdots dx_1
=
\frac{(\zeta-t_0)^k}{k!}
\leq
\frac{(\zeta-\alpha)^k}{k!},
\]
we obtain
\begin{align*}
{\mathbb P}\left(
D_k^{(j)}
\cap
\{|P_{\theta,\rho}^{m(j)}|\geq T\}
\right)
\leq
\frac{2^{2(j+k)}}{T^2(2(j+k)+1)}
\left(
\frac{(\zeta-\alpha)^k}{k!}
\right)^2.
\end{align*}
Case ii): For \(T>0\), we may write
\begin{align*}
\P\left(D_k^{(j)}\right)
&\leq
\P\left(
D_k^{(j)}
\cap
\{|P_{\theta,\rho}^{m(j)}(\beta)|\geq T\}
\right)
+
\P\left(
|P_{\theta,\rho}^{m(j)}(\beta)|<T
\right).
\end{align*}

The conclusion follows from part~i) and Lemma~\ref{poly_lem1}.
\end{proof}

\subsubsection{Proof of Proposition \ref{prop-maj-expec-N}}

We are now in a position to prove Proposition~\ref{prop-maj-expec-N}. The results established in the previous subsections provide quantitative estimates for the distributions of \(P_{\theta,\rho}^{m(j)}(\beta)\) and of their zeros, which will allow us to control the moments of
\begin{equation*}
N\!\left(P_{\theta,\rho}^{m},I_{\rho}\right).
\end{equation*}

We first establish the desired moment estimate on intervals of sufficiently small length.

\begin{lemma}\label{expnpp}
Let \(p\geq 1\). Then there exists a constant \(C(p)>0\) such that, for every interval
$I=[\alpha,\zeta]$ satisfying
$$|I|=\zeta-\alpha<\frac{\pi}{\sqrt{2}},$$
we have
\begin{align*}
\sup_m
\sup_{\theta\in[0,2\pi],\,\rho\in[-a,a]}
\E\left(
N\left(P_{\theta,\rho}^{m},I\right)^p
\right)
\leq C(p).
\end{align*}
\end{lemma}

\begin{proof}
Using the identity valid for any positive random variable $Z$,
\begin{align*}
\E(Z^p)
=
p \int_0^{\infty} t^{p-1}\P(Z>t)\,dt,
\end{align*}
and the monotonicity of \(t\mapsto \P(Z>t)\), we obtain
\begin{align*}
\E\left(N\left(P_{\theta,\rho}^m, I \right)^p\right)
&\leq
p
+
p \sum_{k=1}^{\infty}
k^{p-1}
\P\left(N\left(P_{\theta,\rho}^m, I \right)\geq k\right) \\
&=
p
+
p \sum_{k=1}^{\infty}
k^{p-1}
\P(D_k^{(0)}).
\end{align*}

Since
\begin{equation*}
|I|
=
\zeta-\alpha
<
\frac{\pi}{\sqrt{2}}
<
\frac{m\pi}{\sqrt{2}(m-1)},
\end{equation*}
$P_{\theta,\rho}^{m}$ has at most $2m^2-1$ roots in $I$ (see \cite[Thm.~4]{VanH}), 
{unless it is the null function. With the convention that the null function has no zeros, we thus have $N(P_{\theta,\rho}^m,I)\leq 2m^2-1$ in all cases.}
Consequently, for $q\geq 1$ to be chosen later,
\begin{align*}
\E\left(N\left(P_{\theta,\rho}^m, I \right)^p\right)
&\leq
p
+
p \sum_{k=1}^{q}
k^{p-1}
\P(D_k^{(0)})
+
p \sum_{k=q+1}^{2m^2}
k^{p-1}
\P(D_k^{(0)}).
\end{align*}
Moreover, by Rolle's theorem, on the event \(D_k^{(0)}\), the polynomial
$P_{\theta,\rho}^{m(q)}$ has at least $k-q$ roots in the interval $I$.
Therefore,

\begin{align*}
\E\left(N\left(P_{\theta,\rho}^m, I \right)^p\right)
&\leq
p
+
p\sum_{k=1}^{q}k^{p-1}
+
p\sum_{j=1}^{2m^2-q}(j+q)^{p-1}\P(D_j^{(q)})\\
&\leq
p(1+q^p)
+
p(2q)^{p-1}
\sum_{j=1}^{2m^2-q}
j^{p-1}\P(D_j^{(q)}).
\end{align*}
In the last inequality we have used $k+q \leq 2kq$, $k\geq 1,\quad q\geq 1$.

Applying Lemma~\ref{intrs} with \(T=k^{-\gamma}\), we obtain
\begin{align*}
&\E\left(N\left(P_{\theta,\rho}^m, I \right)^p\right)\\
&\leq p(1+q^p)
+p(2q)^{p-1}
\sum_{k=1}^{2m^2-q}
k^{p-1}
\Bigg(
k^{2\gamma}
\left(
\frac{(\zeta-\alpha)^k}{k!}
\right)^2
\frac{2^{2(q+k)}}{2(q+k)+1}\\
&\hspace{4cm}
+
C(q)
\left(
k^{-\gamma}
+
k^{\gamma/2}m^{-(q+1)/2}
\right)
\Bigg).
\end{align*}

We choose \(\gamma=p+1\), so that $p-1-\gamma=-2.$
Next, we choose $q=[2(2p+\gamma)-1]+1,$
so that $$2p+\gamma-\frac{q+1}{2}\leq 0.$$
Therefore,
\begin{align*}
m^{-\frac{q+1}{2}}
\sum_{k=1}^{2m^2}
k^{p-1+\gamma/2}
&\leq
\frac{(2m^2)^{p+\gamma/2}}{p+\gamma/2}
m^{-\frac{q+1}{2}}\\
&=
\frac{2^{p+\gamma/2}}{p+\gamma/2}
m^{-\frac{q+1}{2}+2p+\gamma}\\
&\leq
\frac{2^{p+\gamma/2}}{p+\gamma/2}.
\end{align*}

Moreover,
\begin{align*}
&\E\left(N\left(P_{\theta,\rho}^m, I \right)^p\right)\\
&\leq
p(1+q^p)
+
p(2q)^{p-1}
\sum_{k=1}^{\infty}
k^{p-1}
\Bigg(
k^{2\gamma}
\left(
\frac{(\zeta-\alpha)^k}{k!}
\right)^2
\frac{2^{2(q+k)}}{2(q+k)+1}\\
&\hspace{4cm}
+
C(q)k^{-\gamma}
\Bigg)
+
C(q)\frac{2^{p+\gamma/2}}{p+\gamma/2}.
\end{align*}

Thus,
\begin{align*}
\sup_m
\sup_{\theta\in[0,2\pi],\,\rho\in[-a,a]}
\E\left(
N\left(P_{\theta,\rho}^m, I \right)^p
\right)
<\infty.
\end{align*}
\end{proof}

Using the fact that every interval $I\subset[-a,a]$ can be decomposed into at most $M(a)$ subintervals $I_1,\ldots,I_{M(a)}$ of length smaller than $\pi/\sqrt{2}$, together with
\begin{equation*}
N\left(P_{\theta,\rho}^m,I\right)
\leq
\sum_{r=1}^{M(a)}
N\left(P_{\theta,\rho}^m,I_r\right),
\end{equation*}
and Jensen's inequality
\begin{equation*}
\left(\sum_{r=1}^{M(a)}x_r\right)^p
\leq
M(a)^{p-1}
\sum_{r=1}^{M(a)}x_r^p,
\end{equation*}
for any sequence $(x_r)_{r=1}^{M(a)}$ of nonnegative numbers, Proposition~\ref{prop-maj-expec-N} follows from Lemma~\ref{expnpp}.

\section{Stein's method}\label{appendix-stein}\label{appendix-Stein}

In this section, we summarize some results on Stein's method for proving Donsker-type theorems (see \cite{Decreus20}). We also refer the reader to our previous work (see Section~1.2 of the supplementary material in \cite{lclp}). First, the result is stated in a general context. Second, the main difference with the previous references is a deeper simplification of the proof of the Lemma 4.7 of   \cite{Decreus20} by using Cochran Theorem.

Throughout this section, $E$ denotes a linear space endowed with a norm $\|\cdot\|_E$.

In what follows, we adopt the following notations:
\begin{itemize}
    \item For any integer $d \geq 1$, we denote by $\left(\vec{e}_j^E\right)_{j=1}^d$ a family of vectors in $E$ and for any integer $n \geq 1$, $\left(c_{i,j}\right)_{1 \leq i \leq n,\ 1 \leq j \leq d}$ denotes a family of real coefficients.
    
    \item We let $\underline{X} = (X_i)_{i=1}^n$ be a family of i.i.d. centered real random variables with unit variance and finite fourth moment.

    \item For any family of real random variables $\underline{Z}=(Z_i)_{i=1}^n$, we denote by $Z_n^c$ the random vector with values in $E$ defined by
    \begin{align}
    Z_n^c := \sum_{j=1}^d \sum_{i=1}^n c_{i,j} Z_i \vec{e}_j^E.
    \end{align}
\end{itemize}

\begin{lemma}\label{d-dim-stein-lemma-summarized}
Let $(E,\|\cdot\|_E)$ be a linear space. 

Let $d, n \geq 1$ be integers. Let ${\mathcal E}=(\vec{e}_j^E)_{j=1}^d$ be a family of unit vectors in $E$, let $(c_{i,j})_{1 \leq i \leq n,\ 1 \leq j \leq d}$ be a family of real coefficients, $\underline{X}$ and let $(G_i)_{i=1}^n$ be a standard Gaussian vector in $\mathbb{R}^n$.

We assume that
\begin{align*}
\text{(H1)}\quad 
&\|c\|_{\infty,{\mathcal E}} := 
\sup_{1 \leq i \leq n} 
\left\|\sum_{j=1}^d c_{i,j}\vec{e}_j^E \right\|_E 
< \frac{1}{2}, \\[0.6em]
\text{(H2)}\quad 
&\exists\, C_{E,X,{\mathcal E}}>0 \ \text{such that} \\[0.3em]
&\sup_n \, \mathbb{E}\!\left[ \|X_n^{c}\|_E^4 \right]^{1/4} \le C_{E,X,{\mathcal E}}, \\[0.4em]
&\sup_n \, \mathbb{E}\!\left[ \|G_n^{c}\|_E^4 \right]^{1/4} \le C_{E,X,{\mathcal E}}.
\end{align*}

Then, the following inequality holds:
\begin{align}
    \sup_{f \in \mathcal{L}ip_{1,1}(E,\|\cdot\|_E)} \left| \mathbb{E}\left[ f(X_n^{c}) \right] - \mathbb{E}\left[ f(G_n^{c}) \right] \right| 
    \leq C \, d \, \|c\|_{\infty,{\mathcal E}} \big|\ln \left( \|c\|_{\infty,{\mathcal E}} \right)\big|
 \left((1+C_{E,X,{\mathcal E}}) \mathbb{E} \left[ |X_1|^4 \right]^{3/4}+ C_{E,X,{\mathcal E}}^2\right).
\end{align}
\end{lemma}

\begin{proof}
We follow the same strategy as in Section~4 of \cite{Decreus20}. The Ornstein--Uhlenbeck semigroup on ${\mathbb R}^n$, endowed with the uniform norm $\|\cdot\|_{\infty}$, is defined by
\begin{align}
    P_t^n(f)(\underline{x}) := \mathbb{E}\left[ f\left( \underline{x} e^{-t} + \beta_t \, \underline{G} \right) \right],
    \quad \forall\, \underline{x} \in \mathbb{R}^n,\ t\geq 0.
\end{align}
for any function $f$ such that there exists $c>0$ satisfying
\[
\sup_{\underline{x}\in \mathbb{R}^n}
e^{-c\|\underline{x}\|_{\infty}}
|f(\underline{x})|
<\infty.
\]
Here, $\underline{G} \sim \mathcal{N}(0, I_n)$ denotes a standard Gaussian vector in $\mathbb{R}^n$, and
$$\beta_t := \sqrt{1 - e^{-2t}}.$$
The semigroup $P$ is ergodic and regularizing. More precisely, for any measurable function $f$ satisfying the above growth condition, the function $P_t^n(f)$ is of class $C^{\infty}$ for all $t > 0$.

This regularizing property is discussed, for instance, in \cite{Decreus20} and \cite{shh}.

We now recall the standard formulas for the first and second derivatives of the Ornstein--Uhlenbeck semigroup.
\begin{align*}
    \partial_{x_i} P^n_t(f)(\underline{x}) 
    &= \frac{e^{-t}}{\beta_t} \, \mathbb{E}\left[ G_i \, f\left( \underline{x} e^{-t} + \beta_t \, \underline{G} \right) \right], \\
    \partial_{x_i x_i}^2 P^n_t(f)(\underline{x}) 
    &= \frac{e^{-3t/2}}{\beta_{t/2}^2} \, \mathbb{E}\left[ G_i \tilde{G}_i \, f\left( \underline{x} e^{-t} + e^{-t/2} \beta_{t/2} \, \underline{G} + \beta_{t/2} \, \tilde{\underline{G}} \right) \right],
\end{align*}
where $\tilde{\underline{G}}$ is an independent copy of $\underline{G}$.

Furthermore, the semigroup satisfies
\begin{align*}
    \frac{\partial}{\partial t} P_t^n(f) = \mathcal{L} \left( P_t^n(f) \right),
\end{align*}
where the generator $\mathcal{L}$ is given by
\begin{align*}
    \mathcal{L}(g)(\underline{x}) = - \langle \underline{x}, \nabla g(\underline{x}) \rangle_{\mathbb{R}^n} 
    + \operatorname{Tr} \left( \nabla^2 g(\underline{x}) \right),
\end{align*}
for any function \( g \in C^2(\mathbb{R}^n;\mathbb{R}) \).
Using that the random variables $(X_i)$ are centered, together with a Taylor expansion, and following the proof of Lemma 4.3 in \cite{Decreus20}, we obtain that there exists a universal constant $C > 0$ such that
\begin{align*}
    \left| \mathbb{E} \left[ \mathcal{L}(g)(\underline{X}) \right] \right| 
    &\leq C \sum_{i=1}^n \mathbb{E} \left[ \left| \frac{\partial^2 g}{\partial x_i^2}(\underline{X}^{-i}) - \frac{\partial^2 g}{\partial x_i^2}(\underline{X}) \right| \right] \\
    &\quad + C \sum_{i=1}^n \int_0^1 \mathbb{E} \left[ \left| X_i^2 
    \left( \frac{\partial^2 g}{\partial x_i^2} \big( \underline{X}^{-i} + \theta(\underline{X} - \underline{X}^{-i}) \big) 
    - \frac{\partial^2 g}{\partial x_i^2}(\underline{X}) \right)\right| \right] d\theta,
\end{align*}

where $\underline{x}^{-i} := \left( (1 - \delta_{j,i}) x_j \right)_{1 \leq j \leq n}$, 
and $\delta_{j,i}$ denotes the Kronecker delta.

Now let $g = P_t^n(f^{c})$, where
\[
    f^{c}(\underline{x}) 
    := f\left( 
        \sum_{j=1}^d \sum_{i=1}^n c_{i,j} x_i \, \vec{e}_j^E
    \right),
    \qquad \text{for all } f \in \mathcal{L}ip_{1,1}(E,\|\cdot\|_E).
\]

In this case, we have:
\begin{align*}
    \frac{\partial^2 P_t^n(f^{c})}{\partial x_i^2} 
    \left( \underline{x}^{-i} + \theta (\underline{x} - \underline{x}^{-i}) \right) 
    &= \frac{e^{-3t/2}}{\beta_{t/2}^2} \, 
    \mathbb{E} \Bigg[ 
        G_i \tilde{G}_i \, 
        f\Bigg( 
            \sum_{j=1}^d \sum_{\hat{i}=1}^n c_{\hat{i},j} \Big( 
                e^{-t} x_{\hat{i}} \left[ 1 - (1 - \theta)\delta_{i,\hat{i}} \right] \\
    &\hspace{5cm} 
                + e^{-t/2} \beta_{t/2} G_{\hat{i}} 
                + \beta_{t/2} \tilde{G}_{\hat{i}} 
            \Big) \vec{e}_j^E
        \Bigg) 
    \Bigg].
\end{align*}

Let \( V_j := \sum_{\hat{i}=1}^n c_{\hat{i},j} G_{\hat{i}} \) and 
\( \tilde{V}_j := \sum_{\hat{i}=1}^n c_{\hat{i},j} \tilde{G}_{\hat{i}} \), 
for all \( j = 1, \ldots, d \). Taking the conditional expectation with respect to 
\[
\sigma\left( V_j, \tilde{V}_j : j = 1, \ldots, d \right),
\]
we obtain
\begin{align}
    \frac{\partial^2 P_t^n(f^{c})}{\partial x_i^2}
    \left( \underline{x}^{-i} + \theta (\underline{x} - \underline{x}^{-i}) \right)
    &= \frac{e^{-3t/2}}{\beta_{t/2}^2} \,
    \mathbb{E} \Bigg[ G_i^{\mathrm{ec}} \tilde{G}_i^{\mathrm{ec}} \,
    f\Bigg( \sum_{j=1}^d \sum_{\hat{i}=1}^n c_{\hat{i},j} \Big( \notag \\
    &\quad e^{-t} x_{\hat{i}} \left( 1 - (1 - \theta)\delta_{i,\hat{i}} \right)
    + e^{-t/2} \beta_{t/2} G_{\hat{i}} + \beta_{t/2} \tilde{G}_{\hat{i}} \Big) \vec{e}_j^E
    \Bigg) \Bigg],
\end{align}

where
\begin{align*}
    G_i^{\mathrm{ec}} &:= \mathbb{E} \left[ G_i \,\middle|\, \sigma\left( V_j : j = 1, \ldots, d \right) \right], \\
    \tilde{G}_i^{\mathrm{ec}} &:= \mathbb{E} \left[ \tilde{G}_i \,\middle|\, \sigma\left( \tilde{V}_j : j = 1, \ldots, d \right) \right].
\end{align*}
For any \( f \in \mathcal{L}ip_{1,1}(E, \| \cdot \|_E) \), we have
\begin{align*}
&\left| 
f\left( 
    \sum_{j=1}^d 
        \left( 
            \sum_{\hat{i}=1}^n 
                c_{\hat{i},j} 
                \left\{ 
                    x_{\hat{i}} \big[ 1 - (1 - \theta)\delta_{i,\hat{i}} \big] 
                    + y_{\hat{i}} 
                \right\} 
        \right) 
    \vec{e}_j^E 
\right)
-
f\left( 
    \sum_{j=1}^d 
        \left( 
            \sum_{\hat{i}=1}^n 
                c_{\hat{i},j} 
                \left\{ 
                    x_{\hat{i}} (1 - \delta_{i,\hat{i}}) 
                    + y_{\hat{i}} 
                \right\} 
        \right) 
    \vec{e}_j^E 
\right)
\right| \\
&\quad \leq 
\left( 
    1 
    + \|x^{c}\|_E 
    + \|y^{c}\|_E 
    + |x_i| \left\| \sum_{j=1}^d c_{i,j} \vec{e}_j^E \right\|_E 
\right)
\, |x_i| \left\| \sum_{j=1}^d c_{i,j} \vec{e}_j^E \right\|_E .
\end{align*}

Thus, using the definition of $\|c\|_{\infty,\mathcal{E}}$ and the assumption $\|c\|_{\infty,\mathcal{E}}<1/2$, we obtain
\begin{align*}
   &\left|
   \frac{\partial^2 P_t^n(f^{c})}{\partial x_i^2 }
      (\underline{x}^{-i} + \theta (\underline{x} - \underline{x}^{-i})) 
   - \frac{\partial^2 }{\partial x_i^2 } P_t^n(f^{c})(\underline{x})
   \right| \\
   &\leq 
   \frac{e^{-5t/2}}{\beta_{t/2}^2} \,
   {\mathbb E}\bigg[
      |G_i^{\mathrm{ec}}| \, |\tilde{G}_i^{\mathrm{ec}}| \,
      \Big(
          1
          + e^{-t} \|x^{c}\|_{E} 
          + e^{-t/2} \beta_{t/2} \|G^{c}\|_{E} 
          + \beta_{t/2} \|\tilde{G}^{c}\|_{E} 
          + |x_i|
      \Big)
      |x_i| \, \|c\|_{\infty,\mathcal{E}}
   \bigg].
\end{align*}

Assumption Hypothesis \ref{hypxp-1} implies that
\begin{align*}
   {\mathbb E}\left[ \|G^{c}\|_{E}^2 \right]
   \leq
   {\mathbb E}\left[ \|G^{c}\|_{E}^4 \right]^{1/2}
   \leq C_{E,X,{\mathcal E}}^2.
\end{align*}
Thus, with $C$ a universal constant,
\begin{align*}
   &\left|
      \frac{\partial^2 P_t^n(f^{c})}{\partial x_i^2 }
        (\underline{x}^{-i} + \theta (\underline{x} - \underline{x}^{-i}))
      - \frac{\partial^2}{\partial x_i^2} P_t^n(f^{c})(\underline{x})
   \right| \\
   &\leq 
   C \frac{e^{-5t/2}}{\beta_{t/2}^2} \,
   \|c\|_{\infty,{\mathcal E}} \, |x_i|\,
   {\mathbb E}\left[ |G_i^{\mathrm{ec}}|^2 \right]
   \left( 
      1 + C_{E,X,{\mathcal E}} + \|x^{c}\|_{E} + |x_i|
   \right).
\end{align*}

Then,
\begin{align*}
     \left| {\mathcal L}(P_t^n(f^{c}))(\underline{x}) \right|
     &\leq (1 + C_{E,X,{\mathcal E}}) 
     \sum_{i=1}^n  
     \frac{e^{- \frac{5t}{2}}}{\beta_{\frac{t}{2}}^{2}} 
     \|c\|_{\infty,{\mathcal E}}  
     \left(1 + |x_i|^{3}\right) 
     {\mathbb E}\left[ |G_i^{\mathrm{ec}}|^{2} \right] 
     \left(1 + \|x^{c}\|_{E} + |x_i| \right).
\end{align*}

Using Hölder's inequality with 
$p = \frac{4}{3}$ and $q = 4$, together with Hypothesis \ref{hypxp-1}, we obtain
\begin{align*}
    {\mathbb E} \left[\left( 1 + |X_i|^{3} \right) \|X_n^{c}\|_{E} \right]
    \leq  
    C_{E,X,{\mathcal E}}
    \left( {\mathbb E} \left[ |X_i|^{4} \right] \right)^{3/4},
    \qquad 
    \text{for all } i = 1, \dots, n.
\end{align*}
Therefore,
\begin{align*}
    {\mathbb E} \left[ \left| {\mathcal L}(P_t^n(f^{c}))(\underline{X}) \right| \right]
    \leq 
    C \, 
    \frac{e^{- \frac{5t}{2}}}{\beta_{\frac{t}{2}}^{2}} 
    \|c\|_{\infty,{\mathcal E}} 
    \left(1 + C_{E,X,{\mathcal E}} \right)
    \left( {\mathbb E}\left[ |X_1|^{4} \right] \right)^{3/4}
    \sum_{i=1}^n {\mathbb E}\left[ |G_i^{\mathrm{ec}}|^{2} \right].
\end{align*}

The  Lemma 4.7 of   \cite{Decreus20} is replaced by the following observation.

Let $\underline{G}^{\mathrm{ec}} = \left( G_i^{\mathrm{ec}} \right)_{i=1}^n$. 
The distribution of the random vector $\underline{G}^{\mathrm{ec}}$ is the same as that of 
$\Pi^{c}(\underline{G})$, where $\Pi^{c}$ denotes the orthogonal projection onto 
$H^{c}$, the linear subspace of ${\mathbb R}^n$ generated by
\begin{align*}
  H^{c} 
  = \operatorname{span} \left\{ 
     \sum_{i=1}^n c_{i,j} \, \vec{e}_i^{\,n} 
     \,:\, j = 1, \dots, d 
  \right\}.
\end{align*}

Here, $\left( \vec{e}_i^{\,n} \right)_{i=1}^n$ denotes the canonical basis of ${\mathbb R}^n$.

The dimension of $H^{c}$ is less than or equal to $d$. Therefore, by Cochran's theorem, the distribution of the random variable $\|\Pi^{c}(\underline{G})\|_2^2$
is a chi-squared distribution with degrees of freedom less than or equal to $d$, and
\begin{align*}
    \sum_{i=1}^n {\mathbb E}\left[( G_i^{\mathrm{ec}})^2 \right]
    &= {\mathbb E}\left\| \underline{G}^{\mathrm{ec}} \right\|_2^2  \\
    &= {\mathbb E}\left\| \Pi^{c}(\underline{G}) \right\|_2^2  \\
    &\leq d.
\end{align*}

Thus,
\begin{align*}
    {\mathbb E} \left[ \left| {\mathcal L}(P_t^n(f^{c}))(\underline{X}) \right| \right]
    \leq C d \frac{e^{-\frac{5t}{2}}}{\beta_{\frac{t}{2}}^2} 
    \|c\|_{\infty,{\mathcal E}}
    \left( 1 + C_{E,X,{\mathcal E}} \right)
    {\mathbb E} \left[ |X_1|^4 \right]^{3/4}.
\end{align*}

We now integrate with respect to $t$ from $\tau$ to $\infty$, and obtain
\begin{align*}
    \left|
    {\mathbb E} \left[ P_{\tau}^n (f^{c})(\underline{X}) \right]
    - {\mathbb E} \left[ f^{c}(\underline{G}) \right]
    \right|
    \leq C\, d \,|\ln(\tau)|\,
    \|c\|_{\infty,{\mathcal E}}
    \left( 1 + C_{E,X,{\mathcal E}} \right)
    {\mathbb E} \left[ |X_1|^4 \right]^{3/4}.
\end{align*}

Moreover, since $f^{c}(\underline{x}) = f(x^{c})$ and $\tau > 0$, we have
\begin{align*}
    \left| \mathbb{E} \left[ f^{c}(\underline{X}) \right] - \mathbb{E} \left[ P_{\tau}^n (f^{c})(\underline{X}) \right] \right|
    \leq C \, \mathbb{E} \left[
        \left( \|X^{c}\|_{E} + \|G^{c}\|_{E} \right)
        \left( (1 - e^{-\tau}) \|X^{c}\|_{E} + \beta_{\tau} \|G^{c}\|_{E} \right)
    \right].
\end{align*}

Using the Cauchy--Schwarz inequality and H2,
\begin{align*}
    \left| \mathbb{E} \left[ f^{c}(\underline{X}) \right] - \mathbb{E} \left[ P_{\tau}^n (f^{c})(\underline{X}) \right] \right|
    &\leq C (1 - e^{-\tau}) \, \mathbb{E} \left[ \|X^{c}\|_{E}^2 + \|G^{c}\|_{E}^2 \right] \\
    &\leq C (1 - e^{-\tau}) \, C_{E,X,{\mathcal E}}^2.
\end{align*}

Then, since
\begin{align*}
    \left| \mathbb{E} \left[ f^{c}(\underline{X}) \right] - \mathbb{E} \left[ f^{c}(\underline{G}) \right] \right|
    &\leq \left| \mathbb{E} \left[ f^{c}(\underline{X}) \right] - \mathbb{E} \left[ P_{\tau}^n (f^{c})(\underline{X}) \right] \right| \\
    &\quad + \left| \mathbb{E} \left[ P_{\tau}^n (f^{c})(\underline{X}) \right] - \mathbb{E} \left[ f^{c}(\underline{G}) \right] \right|,
\end{align*}
we conclude that
\begin{align*}
    \left| \mathbb{E} \left[ f^{c}(\underline{X}) \right] - \mathbb{E} \left[ f^{c}(\underline{G}) \right] \right|
    \leq C \left[ (1 - e^{-\tau}) C_{E,X,{\mathcal E}}^2
    + d |\ln(\tau)|\, \|c\|_{\infty,{\mathcal E}}
      \left( 1 + C_{E,X,{\mathcal E}} \right)
    \right]
    \mathbb{E} \left[ |X_1|^4 \right]^{3/4}.
\end{align*}

We choose $\tau$ such that
\[
1 - e^{-\tau} = \|c\|_{\infty,{\mathcal E}},
\]
and derive
\begin{align*}
    \left| \mathbb{E} \left[ f^{c}(\underline{X}) \right] - \mathbb{E} \left[ f^{c}(\underline{G}) \right] \right|
    \leq C d \|c\|_{\infty,{\mathcal E}}
    \big|\ln \left( \|c\|_{\infty,{\mathcal E}} \right)\big|
    \left( (1 + C_{E,X,{\mathcal E}}) \mathbb{E} \left[ |X_1|^4 \right]^{3/4}
    + C_{E,X,{\mathcal E}}^2 \right).
\end{align*}

This concludes the proof of Lemma~\ref{d-dim-stein-lemma-summarized}.
\end{proof}


\bibliographystyle{alpha}
\bibliography{references}
\end{document}